\documentclass[a4paper,12pt,onecolumn,twoside]{article}
\usepackage{amsmath,amsfonts,amssymb,amsthm,mathrsfs,graphicx,fancyhdr}
\usepackage{cases,xcolor,epstopdf}
\usepackage[driverfallback=hypertex,
linktocpage=true,  
colorlinks=true,   %
linkcolor=blue,    
citecolor=red,     
urlcolor=blue,     
pdfborder={0 0 0}, 
]{hyperref}
\def\p{\partial}
\def\ve{\varepsilon}
\def\f{\frac}

\def\na{\nabla}

\def\al{\alpha}

\def\t{\tilde}
\def\vp{\varphi}
\def\O{\Omega}

\def\th{\theta}

\def\g{\gamma}
\def\G{\Gamma}
\def\si{\sigma}
\def\Si{\Sigma}
\def\dl{\delta}

\def\m{\sqrt{\mu_m}}
\def\q{\qquad}

\def\s{q_0^{\f{\gamma-3}{\gamma-1}}}
\def\a{q_0^{-\f 2{\gamma-1}}}
\def\ds{\displaystyle}

\def\rr{\Bbb{R}}
\def\R{\widetilde{R}}
\def\ss{\sum\limits}

\begin{document}

\title{\bf\Large On the instability problem
of a 3-D\\ transonic oblique shock wave  \footnotetext{
* Li Liang and Yin Huicheng
was supported by the NSFC
(No.11025105), Xu Gang was supported
by the NSFC (No.11101190, No.11371189, No.11271164).}}\vspace{1cm}
\author{Li Liang$^{1}$,\qquad Xu Gang$^{2}$, \qquad Yin
Huicheng$^1$\vspace{0.5cm}\\
\small 1. Department of Mathematics and
IMS, Nanjing University, Nanjing 210093, China.\\\vspace{0.5cm}
\small 2.
Faculty of Science, Jiangsu University, Zhenjiang, Jiangsu 212013,
China.\\
}
\date{ }
\maketitle

\vskip 0.2cm

\maketitle
\vskip 0.5 true cm

\centerline {\bf Abstract} \vskip 0.3 true cm
In this paper, we are
concerned with the instability problem of a 3-D transonic oblique shock wave for the
steady supersonic flow past an infinitely long sharp wedge.
The flow is assumed to be isentropic and irrotational.
It was indicated in pages 317 of [9] that if a  steady supersonic
flow comes from minus infinity and hits a sharp symmetric wedge, then it
follows from the Rankine-Hugoniot conditions and the physical entropy
condition that there possibly
appear a weak shock or a strong shock attached at the edge of the sharp wedge, which
corresponds to a supersonic shock or a transonic shock,
respectively. The question arises which of the two actually occurs.
It has frequently been stated that the strong one is unstable and that,
therefore, only the weak one could occur. However, a
convincing proof of this instability has apparently never been
given. The aim of this paper is to understand such a longstanding open question. We will show that the attached 3-D transonic
oblique shock problem is overdetermined, which implies that the 3-D transonic shock is {\bf unstable} in general.
\vskip 0.3 true cm

 \vskip 0.3 true cm


\vskip 0.4 true cm \centerline{\bf $\S 1.$ Introduction} \vskip 0.3
true cm

\vskip 0.2cm \arraycolsep3pt
\newtheorem{Lemma}{Lemma}[section]
\newtheorem{Theorem}{Theorem}[section]
\newtheorem{Definition}{Definition}[section]
\newtheorem{Proposition}{Proposition}[section]
\newtheorem{Remark}{Remark}[section]
\newtheorem{Corollary}{Corollary}[section]


\setcounter{equation}{1}

In this paper, we are concerned with the instability problem of a 3-D  transonic oblique shock for the steady supersonic flow past an infinitely long sharp wedge (see Figure 1 below).
As indicated in pages 317 of [9]: if a supersonic steady flow comes from minus infinity and hits a sharp symmetric wedge, then it follows from the Rankine-Hugoniot conditions and the physical entropy condition that there will appear a weak shock or a strong shock attached at the edge of the sharp wedge, which corresponds to a supersonic shock or a transonic shock,
respectively. The question arises which of the two shocks actually occurs.
It has frequently been stated that the strong one is unstable and that, therefore, only the weak one could occur. However, a convincing proof of this instability has apparently never been
given. The aim of this paper is to understand such a longstanding open question.
With respect to the 2-D weak oblique shock, under some different assumptions on the 2-D sharp wedge, the authors in [4, 19, 23, 32] have respectively established the local or global existence and stability of a supersonic shock solution or a weak solution for the perturbed supersonic incoming flow past a 2-D sharp curved wedge. For the 3-D weak oblique shock, Chen S.X. in [5] has shown its local stability. With respect to the 2-D strong oblique shock, under certain pressure condition at infinity in the downstream subsonic region, the authors in [6] and [33] have proved the global existence and stability of a transonic shock for the 2-D potential equation and the 2-D full Euler system respectively, which are contrary with the conjecture on the instability of the transonic oblique shock (this instability conjecture has been mentioned in the above).
In addition, for the 2-D unsteady potential equation, the authors in [12] constructed
a self-similar analytic solution which connects an attached  2-D strong shock and an attached 2-D weak shock when a supersonic flow hits a 2-D sharp wedge.
Note that the realistic world is  three-dimensional. The aim of this paper is to show that the attached 3-D transonic shock problem is overdetermined, which means that the 3-D transonic shock is {\bf unstable} in general and further gives a rather positive illustration on the instability of a 3-D transonic oblique shock. This also indicates that the space dimensions are essential for answering the stability or instability of the transonic oblique shocks.

\vskip 0.3 true cm
\includegraphics[width=14cm,height=7.5cm]{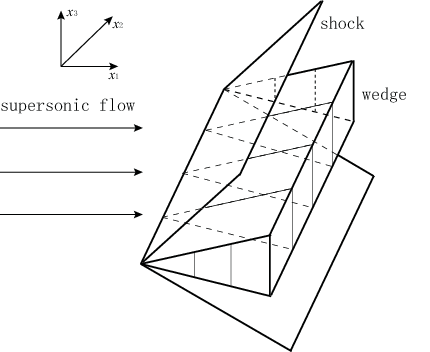}

\centerline{\bf Figure 1. A uniform supersonic flow past a 3-D sharp
wedge}

\vskip 0.3 true cm\hspace{-2cm}
\includegraphics[width=15cm,height=9cm]{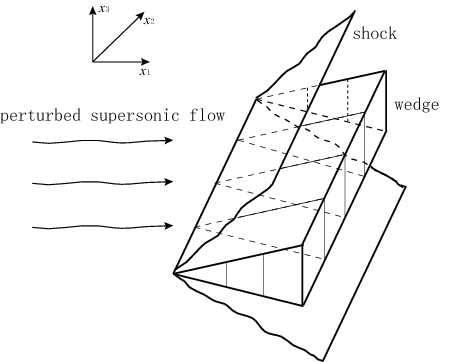}

\centerline{\bf Figure 2. A perturbed supersonic flow past a 3-D sharp
wedge}

\vskip 0.5 true cm
We will assume that the supersonic incoming
flow is of a small perturbation with respect to the constant supersonic state $(\rho_0, q_0, 0, 0)$ and such
a flow hits the sharp 3-D wedge $\{x: x_1\ge 0, x_2\in\Bbb R, -b_0x_1\le x_3\le b_0x_1\}$ along the $x_1$-direction (see Figure 2 above).
Due to the non-interaction property of the transonic oblique shocks on two sides of the wedge, then it suffices
to consider our transonic shock problem only in the upper half-space $x_3\ge 0$ and use a ramp
$\{x: x_1\in\Bbb R, x_2\in\Bbb R, 0\le x_3\le b_0x_1\}$ instead of the wedge (see Figure 3 below).

\vskip 0.3 true cm
\includegraphics[width=13cm,height=7cm]{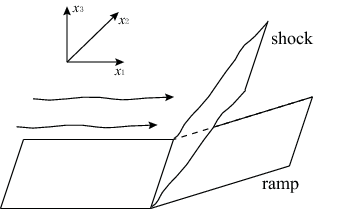}

\centerline{\bf Figure 3. A perturbed supersonic flow past a 3-D ramp}
\vskip 0.6 true cm

The steady and compressible 3-D Euler system is described
by
\begin{equation}
\left\{
\begin{aligned}
&\sum\limits_{j=1}^{3}\p_j(\rho u_j)=0,\\
&\sum\limits_{j=1}^{3}\p_j(\rho u_iu_j)+\p_iP=0,\qquad i=1,2,3,
\end{aligned}
\right.\tag{1.1}
\end{equation}
where $\rho>0$ is the density, $u=(u_1,u_2, u_3)$ is the
velocity, and $P=A\rho^{\g}$ $(1<\g<3)$ is the pressure with $A>0$ a fixed
constant. In addition, $c(\rho)=\sqrt{P'(\rho)}$ is called the local sound speed.

In our paper, we will use the potential equation to describe the
motion of the gas (this model or its variant models have been applied in many other transonic or supersonic
shock problems, one can see [2-3], [13], [17], [25], [34]
and so on). Let $\vp(x)$ be the potential of velocity
$u=(u_1,u_2,u_3)$, i.e., $u_i=\p_i\vp$, then it follows from the
Bernoulli's law that
$$
\f{1}{2}|\na_x\vp|^2+h(\rho)=C_0,\eqno{(1.2)}
$$
here $\na_x=(\p_1,\p_2,\p_3)$, $h(\rho)=\ds\f{c^2(\rho)}{\g-1}$ is the specific enthalpy,
and $C_0=\f{1}{2}q_0^2+h(\rho_0)$ is the
Bernoulli's constant which is determined by the uniform supersonic
incoming flow from the minus infinity with the constant velocity $(q_0,0,0)$ and
the constant density $\rho_0>0$ (see Figure 1 above).

By (1.2) and the implicit function theorem, the density
function $\rho(x)$ can be expressed as
$$\rho=h^{-1}(C_0-\f{1}{2}|\na_x\vp|^2) \equiv
H(\na_x\vp).\eqno{(1.3)}$$
Substituting (1.3) into the mass conservation equation
$\sum\limits_{j=1}^3\p_j(\rho u_j)=0$ in (1.1) yields
$$
\ds\sum_{i=1}^3((\p_i\vp)^2-c^2)\p_i^2\vp
+2\ds\sum_{1\le i<j\le 3}\p_i\vp\p_j\vp\p_{ij}^2\vp=0,\eqno{(1.4)}
$$
where $c=c(H(\na_x\vp))$.

Suppose that the disturbed velocity potentials before and behind the possible attached shock
front $x_3=\chi(x_1,x_2)$ with $\chi(0,x_2)=0$ are denoted by
$\vp^-(x)$ and $\vp^+(x)$ respectively. In this case, the
system (1.4) can be split into two equations, that is,
$\vp^{\pm}(x)$ satisfy the following equations in the
corresponding domains

\begin{align}
\displaystyle\ds\sum_{i=1}^3((\p_i\vp^-)^2-(c^-)^2)\p_i^2\vp^-+2\ds\sum_{1\le i<j\le 3}\p_i\vp^-\p_j\vp^-\p_{ij}^2\vp^-=0 \notag\\
\displaystyle\qquad \qquad \qquad \qquad \text{in \,$\{x_1>0,
x_3>\chi(x_1,x_2)\}$ or $\{x_1\le 0\}$}\tag{1.5}
\end{align}

and
$$
\ds\sum_{i=1}^3((\p_i\vp^+)^2-(c^+)^2)\p_i^2\vp^++2\ds\sum_{1\le i<j\le 3}\p_i\vp^+\p_j\vp^+\p_{ij}^2\vp^+=0\quad \text{in $\{x_1>0,
x_3<\chi(x_1,x_2)\}$} \eqno(1.6)
$$
with $c^{\pm}=c(\rho^{\pm})=c(H(\na_x\vp^{\pm}))$.

It is easy to verify that (1.5) is strictly hyperbolic  with respect
to $x_1$ for $\p_1\vp^{-}>c^-$ and (1.6) is strictly elliptic for
$|\na_{x}\vp^+|<c^+$.

On the ramp surface $\Sigma: x_3=b_0x_1$, $\vp^+$ satisfies
$$b_0\p_1\vp^+ - \p_3\vp^+=0
\quad \text {on $\Sigma$}.\eqno(1.7)$$

Meanwhile, on the possible transonic shock surface $\Gamma: x_3=\chi(x_1,x_2)$ with $\chi(0,x_2)=0$, the Rankine-Hugoniot condition is
$$
[H\p_1\vp]\p_1\chi+[H\p_2\vp]\p_2\chi-[H\p_3\vp]=0 \quad \text { on $\Gamma$},\eqno(1.8)
$$
here we especially point out that the condition $\chi(0,x_2)=0$ comes from the assumption that the transonic
shock is attached at the edge of ramp.

Moreover, the potential $\vp^+(x)$ is continuous across the shock surface $\Gamma$,
namely,
$$
\vp^+(x_1,x_2,\chi(x_1,x_2))=\vp^-(x_1,x_2,\chi(x_1,x_2)),\eqno(1.9)
$$
which obviously means
$$
\vp^+(0,x_2,0)=\vp^-(0,x_2,0).\eqno(1.10)
$$

On $\Gamma$, it follows from  the physical entropy condition that
$$
\rho^-(x_1,x_2,\chi(x_1,x_2))<\rho^+(x_1,x_2,\chi(x_1,x_2)).\eqno(1.11)
$$

In addition, the stable subsonic velocity field behind $\Gamma$ will admit a determined
state:

$$
|\na_x\vp^+|<c^+,\quad\text{and}\quad \text{$\ds\lim_{x_1\to +\infty}\na_{x}\vp^+(x)$ exists for
$b_0x_1\le x_3\le \chi(x_1,x_2)$}.\eqno(1.12)
$$

Finally, we pose the following perturbed initial conditions with respect to the uniform supersonic
constant flow $(\rho_0, q_0, 0, 0)$
$$
\vp^-(0,x_2,x_3) =\ve \vp_0^-(x_2,x_3),\qquad
\p_1\vp^-(0,x_2,x_3) =q_0+\ve \vp_1^-(x_2,x_3),\eqno(1.13)
$$
where $\ve>0$ a small constant, $\vp_i^-(x_2,x_3)\in
C^\infty(\Bbb R^2)$ ($i=0,1$) are supported in $(0,l)$ with respect to the variable
$x_3$, and $l>0$ is some fixed
positive number, moreover, $\vp_i^-(x_2,x_3)=\vp_i^-(x_2+2\pi, x_3)$ holds for $i=0, 1$.
Here we point out that these assumptions on $\vp_i^-(x_2,x_3)$ ($i=0,1$) do not lose the generality
by the finite propagation speed property for the hyperbolic equation (1.5) (one can see more illustrations
in Remark 1.3 below).

In order to solve the transonic shock problem (1.5)-(1.6) together with (1.7)-(1.13),
we will use the partial hodograph transformation in [22] or [26-27] to fix the free boundary
$\Gamma$. To this end, we set $\Phi(x)=\vp^-(x)-\vp^+(x)$ and then it follows
from a direct computation that the problem (1.6)-(1.12) can
be rewritten as

\begin{equation}
\left\{
\begin{aligned}
&\ss_{i,j=1}^3a_{ij}(\nabla_x\vp^--\nabla_x\Phi)\p_{ij}^2\Phi
=\ss_{i,j=1}^3a_{ij}(\nabla_x\vp^--\nabla_x\Phi)\p_{ij}^2\vp^-,\\
&\Phi(x_1,x_2,\chi(x_1,x_2))=0\qquad \qquad \qquad \qquad\qquad \quad \qquad \qquad \qquad \text{on}~~\G,\\
&\big((\rho_+-\rho_-)\p_{1}\vp^--\rho_+\p_{1}\Phi\big)\p_{1}\chi +\big((\rho_+-\rho_-)\p_{2}\vp^--\rho_+\p_{2}\Phi\big)\p_{2}\chi\\ &\qquad\qquad -\big((\rho_+-\rho_-)\p_{3}\vp^--\rho_+\p_{3}\Phi\big)=0\q\q\q\quad\q\text{on}~~\G,\\
&\p_3\Phi-b_0\p_1\Phi=\p_3\vp^--b_0\p_1\vp^-\qquad \qquad \qquad \qquad \quad\qquad\quad\qquad\text{on}~~\Sigma,\\
&\Phi(x_1,x_2+2\pi,x_3)=\Phi(x),\\
&\ds\lim_{x_1+x_3\rightarrow +\infty}\nabla_x\Phi\qquad  \text{exists},\\
\end{aligned}
\right.\tag{1.14}
\end{equation}
where $\vp^-(x)$ is the potential of the supersonic incoming flow, which can be shown to be extended across the shock
$\G$ (see Lemma 2.4 and Remark 2.1 in $\S 2$ below), and
$$
a_{ii}(\na_x\vp^+)=1-\f{\ds (\p_i\vp^+)^2}{\ds (c^+)^2},~~i=1,2,3,\q a_{ij}(\na_x\vp^+)
=-\f{\ds \p_i\vp^+\p_j\vp^+}{\ds (c^+)^2},~~1\le i\neq j\le3.
$$

As in [26-27], we introduce the following partial hodograph
transformation to fix the shock surface $\Gamma$
\begin{equation}
\left\{
\begin{aligned}
&y_1=\f{\ds \Phi(x)}{\ds q_0},\\
&y_2=x_2,\\
&y_3=\f{\ds b_0\Phi(x)}{\ds q_0}+x_3-b_0x_1,\\
\end{aligned}
\right.\tag{1.15}
\end{equation}

In this case, the shock surface $\Gamma $ is changed into $y_1=0$.
Suppose that the inverse transformation of (1.15) is denoted by
\begin{equation}
\left\{
\begin{aligned}
&x_1=u(y),\\
&x_2=y_2,\\
&x_3=y_3-b_0y_1+b_0u(y),\\
\end{aligned}
\right.\tag{1.16}
\end{equation}
where the definition domain of $u(y)$ is the open domain
$Q=\{y\in\Bbb{R}^3: y_1>0, y_2\in \Bbb{R}, y_3>b_0y_1\}$. With respect to the
validity of the invertibility for the transformation (1.15),
one can see the detailed illustrations in $\S 3$ below. In addition, it follows from
(1.9) and (1.16) that
$$
u(0,y_2,0)=0.\eqno{(1.17)}
$$

By (1.14) and (1.16)-(1.17) together with a direct computation, we have
\begin{equation}
\left\{
\begin{aligned}
&L(u, \na_y u, \na_y^2u)\\
&\quad \equiv\ss_{1\le i\le j\le 3} A_{ij}(u, \na_y u)\p_{y_iy_j}^2u+\ds\f{1}{q_0}\ss_{1\le i\le j\le 3} a_{ij}(\na_x\vp^--\na_x\Phi)\p_{x_ix_j}^2\vp^-=0\quad\qquad ~~\text{in}~~Q,\\
&G_1(u, \na_y u)=0\qquad\qquad\qquad\qquad\qquad\qquad\qquad\qquad\qquad\qquad\q\text{on}\qquad y_3=b_0y_1,\\
&G_2(u, \na_y u)=0\qquad\qquad\qquad\qquad\qquad\qquad\qquad\qquad\qquad\qquad\q\text{on}\qquad y_1=0,\\
&u(0,y_2,0)=0,\\
&u(y_1,y_2+2\pi,y_3)=u(y),\\
&\ds\lim_{y_1+y_3\rightarrow +\infty}|\nabla_y u|\quad\text{exists},
\end{aligned}
\right.\tag{1.18}
\end{equation}
where the concrete expressions of $A_{ij}(u, \na_y u)$, $G_1(u, \na_y u)$ and $G_2(u, \na_y u)$ will be given in $\S 3$ below.

Therefore, solving the transonic shock problem (1.5)-(1.6) together with (1.7)-(1.13) is completely
equivalent to solving (1.18). However, unfortunately, (1.18) is an overdetermined problem due to the restriction
$u(0, y_2, 0)=0$ for all $y_2\in\Bbb R$. More precisely,
the following problem can be shown to be uniquely solvable for any fixed $y_2^0\in\Bbb R$
\begin{equation}
\left\{
\begin{aligned}
&L(u, \na_y u, \na_y^2u)=0\qquad\qquad\qquad\qquad\qquad\qquad\qquad\qquad\qquad\qquad\qquad
~~\text{in}\quad ~~Q,\\
&G_1(u, \na_y u)=0\qquad\qquad\qquad\qquad\qquad\qquad\qquad\qquad\qquad\q\text{on}\qquad y_3=b_0y_1,\\
&G_2(u, \na_y u)=0\qquad\qquad\qquad\qquad\qquad\qquad\qquad\qquad\qquad\q\text{on}\qquad y_1=0,\\
&u(0,y_2^0,0)=0,\\
&u(y_1,y_2+2\pi,y_3)=u(y),\\
&\ds\lim_{y_1+y_3\rightarrow +\infty}\nabla_y u\qquad\text{exists}.
\end{aligned}
\right.\tag{1.19}
\end{equation}

Here we emphasize that the difference between (1.18) and (1.19) is: only $u(0,y_2^0,0)=0$
holds for some fixed point $(0, y_2^0, 0)$ in (1.19) other than $u(0,y_2,0)=0$ holds in (1.18) for all $y_2\in\Bbb R$.

We now state our main result in this paper.

\vskip 0.2 true cm

{\bf Theorem 1.1.} {Assume that $b_0>0$ is a small constant, namely, the 3-D ramp is sharp,
then for suitably large supersonic incoming speed $q_0$, the nonlinear problem (1.19) admits a unique smooth solution $u(y)$ in $Q$, which illustrates that (1.18) is overdetermined.}
\vskip 0.2 true cm
{\bf Remark 1.1.} {\it The detailed descriptions on the regularities of $u(y)$ in Theorem 1.1 will be given in Theorem 3.1 of $\S 3$ below.}
\vskip 0.2 true cm
{\bf Remark 1.2.} {\it By the overdetermination of the transonic shock problem
(1.5)-(1.6) together with (1.7)-(1.13) in Theorem 1.1, we know that the transonic oblique shock
is unstable in general. If one could find another point $(0, y_2^1, 0)\not=(0, y_2^0, 0)$ such that $u(0,y_2^1,0)\not=0$ holds for the solution $u$ to (1.19), then the conjecture of the instability for the attached transonic oblique shock is verified in case of the potential flow equation.}
\vskip 0.2 true cm
{\bf Remark 1.3.} {\it Although we pose some restrictions on the perturbed initial data $\vp^-_i(x_2,x_3)$
($i=0,1$) in (1.13), this does not lose the generality. Indeed, if $\vp^-_i(x_2,x_3)\in C_0^{\infty}(\Bbb R\times\Bbb R^+)$, then we can take the smooth initial data $\vp^-_{i,L}(x_2,x_3)$ with a period$-L$ for the variable $x_2$ instead of $\vp^-_i(x_2,x_3)$ in (1.13), where $L>1$ is any fixed constant and $\vp^-_{i,L}(x_2,x_3)=\vp^-_{i}(x_2,x_3)$ holds for $x_2\in [0, L]$. In this case, the related problem (1.19) on $u_L(y)$ can be solved by Theorem 1.1. Moreover, it follows from the proof procedure of Theorem 1.1 that all the $u_L(y)$ for $L\ge 1$ are uniformly bounded for $y\in [0, \infty)\times K\times[0, \infty)$, here $K$ is any fixed compact set in $\Bbb R$. Subsequently, letting $L\to\infty$, then (1.19) can be solved for the given initial data $\vp^-_i(x_2,x_3)$ ($i=0,1$).}

Since the oblique shocks and the conic shocks are two kinds of basic attached shocks for the supersonic flows
past the sharp wedges or sharp cones, we now comment on some interesting and systematic results on the attached conic shocks. It was indicated in pages 317-318 and 414 of [9] that if a uniform supersonic steady
flow  hits a sharp cone in direction of its axis, then it follows from the Rankine-Hugoniot conditions and the physical entropy condition that there possibly occur a weak conic shock (see Figure 4 below) and a strong conic shock (see Figure 5 below) attached at the tip of the cone (this physical phenomena is completely similar to that for the steady supersonic flow past a sharp wedge). For the potential equation, under various assumptions on the supersonic incoming flows and the sharp vertex angles of the conic bodies, the authors have established  the local or global existence and stability of the weak conic shocks or strong conic shocks, one can see [7-11], [17-18], [21], [25-27] and the references therein.
For the full Euler system, because of the essential influences of the rotations,
the authors in [30] and [28] have shown the nonexistence of the global weak solution
with only one stable weak conic shock and the instability of a global transonic conic shock
for the steady supersonic flow past a sharp conic body, respectively. Therefore, these results have
given a basic answer for the global stability or instability of weak and strong
conic shocks.

\includegraphics[width=12cm,height=8cm]{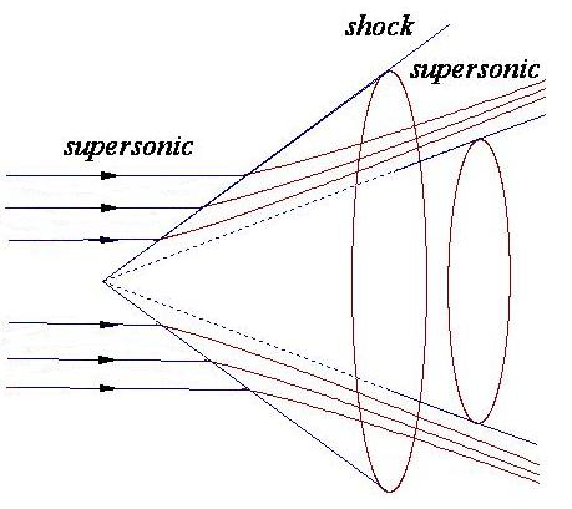}

\centerline{\bf Figure 4. A supersonic shock for the supersonic flow past a 3-D sharp cone}

\vskip 1.2true cm

\includegraphics[width=12cm,height=9cm]{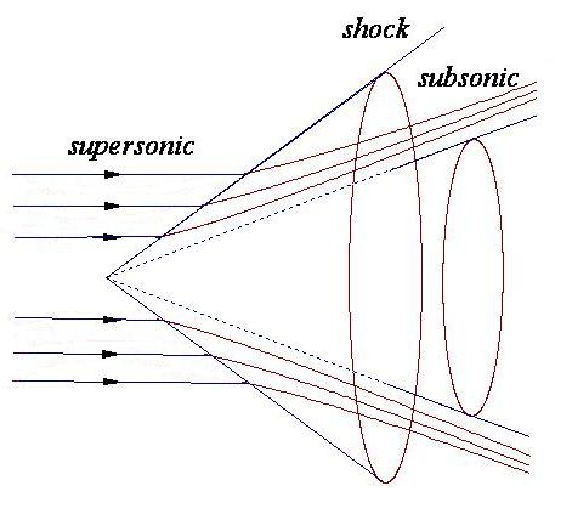}

\centerline{\bf Figure 5. A transonic shock for the supersonic flow past a 3-D sharp cone}

\vskip 0.6 true cm

We now mention some transonic shock problems studied recently in
[2-3], [6], [12-13], [26-27], [33-34] and the references therein. In
these papers, the  considered domains are either 2-D polygons or 3-D conic bodies. For the 2-D polygon
domains (see [2-3], [6], [12-13], [33-34]), it follows from the maximum principle and the barrier function method
for the second order
elliptic equations in the 2-D irregular regions that one can obtain at
least $C^{1,\al} (0<\al<1)$ regularities of the corresponding shock curves and the downstream subsonic solutions.
The $C^{1,\al}$ regularity is
crucial in studying the free boundary problem on the second
order nonlinear elliptic equations whose coefficients contain the
gradients of solution. For the 3-D conic domains (see [26-27], where the maximum principle
can not be used directly), by utilizing the Sturm-Liouville
theorem and separation variable method,  we can write out the expression
of the solution to the linearized elliptic equations and the corresponding
boundary conditions, subsequently we can obtain the regularity, existence and a priori
estimates of the solution to the
nonlinear problem in the conic domain and the suitable weighted H\"older space
with two different weights near the conic point and at infinity.
However, it seems rather difficult for us to choose a weighted H\"older  space as in [26-27]
to deal with the corresponding linearized equation of (1.19) in the 3-D unbounded wedge domain.
The reason is that: we can not expect the solution $u$ of (1.19) to  satisfy $u(0, y_2, 0)\equiv 0$
for all $y_2\in\Bbb R$, thus such properties of $|u(y)|\le Cy_1^{\dl_0}$ ($\dl_0>0$) near the edge $y_1=y_3=0$
and $|u(y)|\le Cy_1^{-\dl_1}$ ($\dl_1>0$) for sufficiently large $y_1>0$
can not hold simultaneously. Note that such kind of weighted space in [26-27] is crucial in
deriving the solvability of the related linearized
potential equation in the unbounded conic domain by the separation variable method.
Therefore, in this paper we should use some other ingredients to overcome this difficulty
so that our problem (1.19) in the unbounded wedge region can be treated.

Next we comment on the proof of Theorem 1.1. To solve (1.19), we will linearize
the nonlinear problem on $u$ by use of the largeness of $q_0$ and
the detailed properties on the background solution, here the
so-called background solution is referred as one to the problem
(1.19) when the uniform supersonic steady flow $(\rho_0, q_0, 0, 0)$ hits the ramp $\{x: x_1>0, x_2\in\Bbb R,  0<x_3<b_0x_1\}$
along $x_1-$direction. By the
linearization, we essentially obtain the Laplacian equation $\Delta
v=f$ in $\Bbb R^3$ with two Neumann boundary conditions on two
different planes in an angular region, a vanishing condition of the first
order derivatives $Dv$ at infinity and a restriction condition $v(0,0,0)=0$
(one can see (4.1) in $\S 4$ below). To study the solvability,
regularity of $v$ and derive the a priori estimates of $v$ in the unbounded wedge region, at first we will
restrict our linearized problem in a bounded wedge domain in addition a Neumann-type boundary condition
on the cut-off surface (see (4.4) of $\S 4$). In this case,
by use of Sturm-Liouville theorem, the
separation variable method, we can derive the concrete expression of the solution
$v_L$ to the cut-off problem (4.4). It follows from the detailed estimates on the related
eigenvalues and eigenfunctions that we can get the existence and
$C^{1,\delta} (0<\delta<1)$ regularity  of $v_L$ up to the
boundaries (including the two boundaries of the angular domain) in Lemma 5.2 of $\S 5$. Based on these
crucial estimates and the scaling techniques for the linear elliptic equations, we can obtain the global estimates
of $v$ in the whole wedge domain by taking the limit $L\to\infty$ for $v_L$. Finally,
by taking a suitable iteration scheme and using the
largeness of $q_0$ and the uniform estimates on the solution to the linearized problem,
we can complete the proof on Theorem 1.1.

Our paper is organized as follows. In $\S 2$, at first, we
give some useful information on the background
solution for large $q_0$, which essentially corresponds to a 2-D transonic oblique shock solution for the uniform supersonic flow
past a 2-D sharp wedge. Secondly, we will define some weighted H\"older spaces which
will be used in subsequent sections. Thirdly, we list or derive some basic properties of the modified Bessel
functions of the first and second kind of order $\nu$ ($\nu\in\Bbb R$) so that one can use the separation
method to study our problems in subsequent $\S 3$-$\S 6$. Fourthly, a global solvability on the problem
(1.5) with (1.13) near the shock $\G$ is given.
In $\S 3$, we will reformulate the problem (1.5)-(1.6) together with (1.7)-(1.13) into
(1.18) meanwhile the detailed expressions of the coefficients in (1.18) can be given. Moreover,
a more precise description on Theorem 1.1 in the weighted H$\ddot o$lder space will be given
in Theorem 3.1. In $\S 4$, the linearized equation and boundary
conditions of (1.19) are given in (4.1), subsequently, a cut-off problem (4.4) with a suitable
Neumann boundary condition on the cut-off surface $\sqrt{y_1^2+y_3^2}=L$ is studied in details,
where the solvability of (4.4) and the rough regularity of the solution $v_L$ to (4.4) in related weighted
H\"older space are shown. In $\S 5$, the higher regularities of $v_L$ in (4.4) are obtained by the
classical Schauder estimate and the regularity theory of solutions to the second order elliptic equations
in a 3-D bounded angular region. Moreover, the global solvability  and estimates of the solution to (4.1) in the unbounded
angular domain $Q$ are established. In $\S 6$, the uniqueness of solution $u$ to (4.1) is proved by the separation variable
method other than by the usual  maximum principle for the second order elliptic equations since it seems that there is no maximum principle
for the problem (4.1) due to the 3-D unbounded angular region and the Neumann boundary conditions (note that $u$ and $\na_y u$ are
actually unknown on the edge $\{y_1=y_3=0\}$ of $Q$). Based on
the estimates in $\S 4$-$\S 6$, Theorem 3.1 and further Theorem 1.1 can be shown in $\S 7$. In addition, some
complicated and useful computations are carried out in the Appendix.

In what follows, we will use the following conventions:

For large $q_0$, $O(q_0^{-\nu})$ $(\nu>0)$ denotes a bounded quantity
such that $|O(q_0^{-\nu})|\le Cq_0^{-\nu}$, where $C>0$ is a generic positive constant.

The Gamma function $\Gamma(a)$
for $a>0$ and  the Beta function $B(a,b)$ ($a,b>0$) are respectively defined as

\begin{align*}
\displaystyle &\Gamma(a)=\int_0^{+\infty} t^{a-1}e^{-t} dt\qquad \text{for $a>0$},\\
\displaystyle &B(a,b)=\int_0^1 t^{a-1} (1-t)^{b-1} dt\qquad \text{for $a>0, b>0$}.
\end{align*}

\vskip 0.4 true cm \centerline{\bf $\S 2.$ Some preliminaries} \vskip 0.4 true cm

At first, we study the background solution to (1.5)-(1.6) together with (1.7)-(1.13)
and derive some useful properties of the transonic oblique shock for the uniform supersonic incoming flow past a sharp ramp.
Since $u_2=0$ always holds in the background solution, it is only required to consider a 2-D  transonic oblique shock
problem temporarily.

\includegraphics[width=12cm,height=7cm]{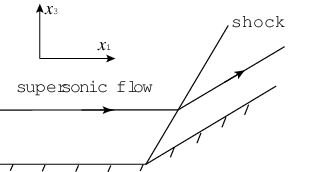}

\centerline{\bf Figure 6. A uniform supersonic flow past a 2-D sharp
ramp}

\vskip 0.6 true cm

Suppose that there is a uniform 2-D supersonic flow $(q_0, 0)$ with constant density $\rho_0>0$
which comes from minus infinity, and the flow hits the 2-D sharp ramp
in the $x_1-$direction (see the Figure 6 above). The ramp boundary is described by $x_3=b_0x_1$
($b_0>0$), then as indicated in pages 317 of [9], there exists a critical value $b^*$  such that there will
appear a transonic shock $x_3=s_0x_1$ $(s_0>b_0)$ attached at
the edge of ramp for $b_0<b^*$. Moreover, it follows from Rankine-Hugoniot conditions and
the boundary condition on the ramp that the constant downstream subsonic flow
$(\rho_0^+, u_{10}^+, u_{30}^+)$ satisfies
\begin{equation}
\left\{
\begin{aligned}
&s_0(\rho_0^+u_{10}^+-\rho_0q_0)-\rho_0^+u_{30}^+=0,\\
&u_{10}^+-q_0+s_0u_{3+}=0,\\
&\f{1}{2}((u_{10}^+)^2+(u_{30}^+)^2)+h(\rho_0^+)\equiv C_0=\ds\f{q_0^2}{2}+h(\rho_0),\\
&u_{30}^+=b_0u_{10}^+\\
\end{aligned}
\right.\tag{2.1}
\end{equation}
with
$$
(u_{10}^+)^2+(u_{30}^+)^2<c^2(\rho_0^+).\eqno{(2.2)}
$$

In addition, the following physical entropy condition holds
$$
\rho_0<\rho_0^+.\eqno{(2.3)}
$$

With respect to the properties of the downstream subsonic flow $(\rho_0^+, u_{10}^+, u_{30}^+)$ and
the slope $s_0$ of the transonic oblique shock, for large $q_0$, we have

\vskip 0.2 true cm

{\bf Lemma 2.1.} {\it If $q_0$ is large and $b_0>0$ is fixed,  then one has for $1<\gamma<3$

(i) $s_0=\ds \f{1}{b_0\rho_0}\bigl(\f{\gamma-1}
{2A\gamma}\bigr)^{\f{1}{\gamma-1}}q_0^{\f{2}{\gamma-1}}\bigl(1
+O(q_0^{-\f{2}{\g-1}})+O(q_0^{-2})\bigr)$.

(ii) $u_{10}^+=O(q_0^{\f{\g-3}{\g-1}})$.

(iii) $u_{30}^+=O(q_0^{\f{\g-3}{\g-1}})$.

(iv) $\rho_0^+=\biggl(\f{\gamma-1}
{2A\gamma}\biggr)^{\f{1}{\gamma-1}}q_0^{\f{2}{\g-1}}\bigl(1+O(q_0^{-2})+O(q_0^{-\f{2}{\g-1}})\bigr)$.

(v)
$c^2(\rho_0^+)=\f{\g-1}{2}q_0^2\bigl(1+O(q_0^{-2})+O(q_0^{-\f{4}{\g-1}})\bigr)$.

(vi) $(q_0^+)^2-c^2(\rho_0^+)=-\ds\f{\g-1}{2}q_0^2\bigl(1+O(q_0^{-2})+O(q_0^{-\f{4}{\g-1}})\bigr)$,
here and below $(q_0^+)^2=(u_{10}^+)^2+(u_{30}^+)^2$.}

{\bf Proof.} (i) It follows from (2.1) that
\begin{equation}
\left\{
\begin{aligned}
&u_{10}^+=q_0-\f{\ds s_0^2q_0(\rho_0^+-\rho_0)}{\ds (1+s_0^2)\rho_0^+},\\
&u_{30}^+=\f{\ds s_0q_0(\rho_0^+-\rho_0)}{\ds (1+s_0^2)\rho_0^+},\\
&h(\rho_0^+)-h(\rho_0)-\f{\ds s_0^2q_0^2((\rho_0^+)^2-\rho_0^2)}{\ds 2(1+s_0^2)(\rho_0^+)^2}=0.
\end{aligned}
\right.\tag{2.4}
\end{equation}

From the third equation in (2.4), we have
$$
\f{A\gamma}{\gamma-1}((\rho_0^{+})^{\gamma-1}-\rho_0^{\gamma-1})
=\f{s_0^2 q_0^2}{2(1+s_0^2)}\bigl(1-(\f{\rho_0}{\rho_0^{+}})^2\bigr).
$$
Denoting by $\alpha=\ds\f{\rho_0^{+}}{\rho_0}$, then one has
$$
\alpha^{\gamma-1}=1+\f{\rho_0^{1-\gamma}(\gamma-1)q_0^2}{2A\gamma}
\bigl(1-\f{1}{1+s_0^2}\bigr)(1-\f{1}{\alpha^2}).\eqno{(2.5)}
$$
Therefore, for large $q_0$ and $\al>1$, one derives
$\alpha=\biggl(\ds\f{\rho_0^{1-\gamma}(\gamma-1)}
{2A\gamma}\biggr)^{\f{1}{\gamma-1}}q_0^{\f{2}{\gamma-1}}\bigl(1
+O(q_0^{-2})\bigr)$
and
\begin{equation}
\left\{
\begin{aligned}
&u_{30}^{+}=\ds\f{s_0q_0}{1+s_0^2}(1-\f{1}{\alpha}),\\
&u_{10}^+=\ds\f{q_0}{1+s_0^2}+\f{s_0^2q_0}{(1+s_0^2)\al}.
\end{aligned}
\right.\tag{2.6}
\end{equation}

Furthermore, by
$$
u_{30}^{+}=b_0u_{10}^+, \eqno{(2.7)}
$$
we arrive at
$$
s_0=\ds\f{\al}{b_0}\biggl(1-\f{1}{\al}-\f{b_0}{s_0}\biggr)=\f{1}{b_0\rho_0}\biggl(\f{\gamma-1}
{2A\gamma}\biggr)^{\f{1}{\gamma-1}}q_0^{\f{2}{\gamma-1}}\bigl(1
+O(q_0^{-\f{2}{\g-1}})+O(q_0^{-2})\bigr),\eqno{(2.8)}
$$
which leads to (i) of Lemma 2.1.

(ii) and (iii) come from (2.6) and (2.8) directly.

(iv)-(vi) come from the system (2.1) and
(i)-(iii).\qed

Next, we introduce some weighted H$\ddot{o}$lder spaces which are motivated by the Chapter 6 of [15] and [14].
These spaces are also applied in [2], [6], [26-27], [33] and so on.

Let $D\subset \Bbb R^3$ be an open set including the $x_2-$axis, for $x, y\in D$, we define $r_x^2=x_1^2+x_3^2$ and $r_{x,y}=\min(r_x,r_y)$. For $m\in\Bbb N\cup\{0\}$, $0<\alpha<1$, $k,l\in \Bbb{R}$ and $u\in C_{loc}^{m,\al}(\bar D\backslash \{(0,x_2,0):x_2\in \Bbb R\})$, we define

\begin{align*}
\displaystyle&[u]_{m,0;D}^{(k,l)}=\max\bigl\{\ds\sup_{
0<r_x<1}\ds\sum_{|\beta|=m} |r_x^{\max(k+m,0)}D^{\beta}u(x)|,
\ds\sup_{r_x>1}\ds\sum_{|\beta|=m}|r_x^{l+m}D^{\beta}u(x)|\bigr\},\\
\displaystyle&[u]_{m,\alpha;D}^{(k,l)}=\max\bigl\{\ds\sup_{0<r_{x,y}<1}
\ds\sum_{|\beta|=m}r_{x,y}^{\max(k+m+\alpha,0)}\f{|D^{\beta}u(x)-D^{\beta}u(y)|}{|x-y|^{\alpha}},\\
\displaystyle&\qquad\qquad \qquad \qquad \qquad  \ds\sup_{r_{x,y}>1}
\ds\sum_{|\beta|=m}r_{x,y}^{l+m+\alpha}\f{|D^{\beta}u(x)-D^{\beta}u(y)|}{|x-y|^{\alpha}}\bigr\},\\
\displaystyle&\|u\|_{m,0;D}^{(k,l)}=\ds\sum_{j=0}^m[u]_{j,0;D}^{(k,l)},\\
\displaystyle&\|u\|_{m,\alpha;D}^{(k,l)}=\|u\|_{m,0;D}^{(k,l)}+[u]_{m,\alpha;D}^{(k,l)},
\end{align*}

and the related function space is defined as
$$
H_{m,\alpha}^{(k,l)}(D)=\{u\in
C_{loc}^{m,\alpha}(\bar D\backslash\{(0,x_2,0):x_2\in \Bbb R\}):
\|u\|_{m,\alpha}^{(k,l)}<+\infty\}.
$$

Let $E=D\bigcap \{(x_1,x_2,x_3):x_1^2+x_3^2<1, x_2\in \Bbb R\}$,  which is a domain
near the $x_2-$axis. For $m\in\Bbb N\cup\{0\}$, $0<\alpha<1$, $k\in \Bbb{R}$ and $u\in C_{loc}^{m,\al}(\bar
E\backslash \{(0,x_2,0):x_2\in \Bbb R\})$, we define
\begin{align*}
\displaystyle&[u]_{m,0;E}^{(k,\star)}=\ds\sup\ds\sum_{|\beta|=m} |r_x^{\max(k+m,0)}D^{\beta}u(x)|,\\
\displaystyle&[u]_{m,\alpha;E}^{(k,\star)}=\ds\sup
\ds\sum_{|\beta|=m}r_{x,y}^{\max(k+m+\alpha,0)}\f{|D^{\beta}u(x)-D^{\beta}u(y)|}{|x-y|^{\alpha}},\\
\displaystyle&\|u\|_{m,0;E}^{(k,\star)}=\ds\sum_{j=0}^m[u]_{j,0;E}^{(k,\star)},\\
\displaystyle&\|u\|_{m,\alpha;E}^{(k,\star)}=\|u\|_{m,0;E}^{(k,\star)}+[u]_{m,\alpha;E}^{(k,\star)},
\end{align*}
and the related function space is defined as
$$
H_{m,\alpha}^{(k,\star)}(E)=\{u\in
C_{loc}^{m,\alpha}(\bar E\backslash\{(0,x_2,0):x_2\in \Bbb R\}):
\|u\|_{m,\alpha}^{(k,\star)}<+\infty\}.
$$

Analogously, set $F=D\bigcap \{(x_1,x_2,x_3):x_1^2+x_3^2>1, x_2\in \Bbb R\}$ which is a domain away from $x_2-$axis. We define for $m\in\Bbb N\cup\{0\}$, $0<\alpha<1$, $l\in \Bbb{R}$ and $u\in C_{loc}^{m,\al}(\bar
F)$,

\begin{align*}
\displaystyle&[u]_{m,0;F}^{(\star,l)}=
\ds\sup\ds\sum_{|\beta|=m}|r_x^{l+m}D^{\beta}u(x)|,\\
\displaystyle&[u]_{m,\alpha;F}^{(\star,l)}=\ds\sup
\ds\sum_{|\beta|=m}r_{x,y}^{l+m+\alpha}\f{|D^{\beta}u(x)-D^{\beta}u(y)|}{|x-y|^{\alpha}}\bigr\},\\
\displaystyle&\|u\|_{m,0;F}^{(\star,l)}=\ds\sum_{j=0}^m[u]_{j,0;F}^{(\star,l)},\\
\displaystyle&\|u\|_{m,\alpha;F}^{(\star,l)}=\|u\|_{m,0;F}^{(\star,l)}+[u]_{m,\alpha;F}^{(\star,l)},
\end{align*}

and the related function space is defined as
$$
H_{m,\alpha}^{(\star,l)}(F)=\{u\in
C_{loc}^{m,\alpha}(\bar F):
\|u\|_{m,\alpha}^{(\star,l)}<+\infty\}.
$$

From the definitions of $H_{m,\alpha}^{(k,l)}(D)$, $H_{m,\alpha}^{(k,\star)}(E)$ and
$H_{m,\alpha}^{(\star,l)}(F)$, one easily knows that the space $H_{m,\alpha}^{(k,l)}(D)$
can be split into the two subspaces $H_{m,\alpha}^{(k,\star)}(E)$ and
$H_{m,\alpha}^{(\star,l)}(F)$.

For the domain $E$ defined above, we set $E_{\sigma}=\{x\in E:r_x>\sigma\}$
for some positive constant $\sigma>0$. The following weighted H\"older space
$H_a^{(b)}(E)$ $(a>0, b\in\Bbb R)$ was introduced in [14]:
$$
H_a^{(b)}(E)=\{u(x)\in C_{loc}^a(E): \ds\sup_{\sigma>0}\sigma^{a+b}\|u\|_{a;E_\sigma}<\infty\},
$$
where $\|\cdot\|_{a;E_\sigma}$ stands for the norm of the H\"older space $C^a(E_{\sigma})$.
In addition, as in [14], we denote by
$$|u|_{a;E}^{(b)}=\ds\sup_{\sigma>0}\sigma^{a+b}\|u\|_{a;E_\sigma}.$$
Then we have
\vskip 0.2 true cm
{\bf Lemma 2.2.} {\it (i) If $a>0$, then $H_a^{(-a)}(E)=C^a(E)$.

(ii) If $a\ge a'\ge0$, $a'+b\ge0$, $b\notin \Bbb N$, then $|u|_{a'; E}^{(b)}\le C|u|_{a; E}^{(b)}$.

(iii) If $0<b<1$ and $0<a<1$, then $|u|_{a;E}^{(b)}\le C||u||_{0, a;E}^{(b,\star)}$.
}

{\bf Proof.} (i) and (ii) can be found in Lemma 2.1 of [14]. We now prove (iii).
Noticing that for any $x\in E_{\sigma}$, one has $\sigma\le r_x\le1$. This derives
$\si^{a+b}|u(x)|\le r_x^{a+b}|u(x)|\le r_x^b|u(x)|$ and
$$
\ds\sup_{\si>0}\big(\si^{a+b}\ds\sup_{x\in E_{\si}}|u(x)|\big)\le \ds\sup_{x\in E}r_x^b|u(x)|.\eqno(2.9)
$$
And similarly, we have
$$
\ds\sup_{\si>0}\big(\si^{a+b}\ds\sup_{x,y\in E_{\si}}\f{\ds |u(x)-u(y)|}{|x-y|^a}\big)\le \ds\sup_{x,y\in E}r_{x,y}^{a+b}\f{\ds |u(x)-u(y)|}{|x-y|^a}.\eqno(2.10)
$$
Therefore it follows from (2.9)-(2.10) that (iii) holds.
\qed

In order to apply the Sturm-Liouville theorem and separation variable method to solve
the linearized problem of (1.19), we require to list or establish  some properties on
the modified Bessel functions $I_\nu(t)$ and $K_\nu(t)$
of the first and second kind of order $\nu$ ($\nu\in\Bbb R$) respectively, where $t\in\Bbb R$, and $I_\nu(t)$ and $K_\nu(t)$
are two linearly independent solutions to the ordinary differential equation $t^2\ds\f{d^2w}{dt^2}+t\ds\f{dw}{dt}+(t^2-\nu^2)w=0$.

\vskip 0.2 true cm

{\bf Lemma 2.3.} {\it For $I_\nu(t)$ and $K_\nu(t)$, we have

(i) $I_\nu(t)$ and $K_\nu(t)$ have the following integral representations:

\begin{eqnarray*}
&&I_\nu(t)=\frac{2 e^t (2t)^\nu}{\sqrt{\pi} \Gamma(\nu+\frac 12)} \int_0^{1}e^{-2tu^2} u^{2\nu}(1-u^2)^{\nu-\frac 12} du
\qquad\text{when\quad $Re\nu>-\frac 12$},\\
&&K_\nu(t)= \frac{\sqrt \pi e^{-t}}{\sqrt{2t} \Gamma(\nu+\frac 12)}\int_0^\infty e^{-u} u^{\nu-\frac 12} (1+\frac{u}{2t})^{\nu-\frac 12} du
\qquad\text{when\quad $Re\nu>-\frac 12$ \quad and\quad  $t>0$.}\\
\end{eqnarray*}

(ii) $I'_\nu(t)=I_{\nu+1}(t)+\ds\frac{\nu}{t}I_{\nu}(t)$ and $K'_\nu(t)=-K_{\nu-1}(t)-\ds\frac{\nu}{t} K_{\nu}(t)$.

Especially,
$I'_0(t)=I_1(t)$ and $K'_0(t)=-K_{-1}(t)=-K_{1}(t)$.

(iii) For any $t>0$, then
\begin{eqnarray*}
&&(a)\quad  I_0(t) \leqslant e^t,\\
&&(b) \quad  K_0(t) \leqslant \frac{\sqrt \pi e^{-t}}{\sqrt{2t}}.
\end{eqnarray*}

(iv) If $\nu>\frac 12$ and $t<1$, then
\begin{eqnarray*}
&&(a)\quad  0<I_\nu(t)\leqslant \frac{e^t  (\frac t2)^\nu }{\Gamma(\nu+1)},\\
&&(b)\quad  0<K_\nu(t) \leqslant \frac{e^{t} \Gamma(\nu) 2^{\nu-1}}{t^\nu}.
\end{eqnarray*}

(v) If $\nu>\frac 12$ and $t\ge 1$, then
$$
0<I_\nu(t)\leqslant  \frac{e^t}{\sqrt{2\pi t}}.
$$

(vi) If $\frac 12<\nu\le M$ and  $t \geqslant 1$, then there exists a constant $C_M>0$ independent of $\nu$ such that
$$
0<K_\nu(t) \leqslant C_M \frac{\sqrt \pi e^{-t}}{\sqrt{2t}}.
$$

(vii) When  $|x|$ is large and $\mu=4\nu^2$, then the following asymptotic expansions hold
\begin{eqnarray*}
&& I'_\nu(t)\sim \frac{e^t}{\sqrt{2\pi t}} \biggl(1-\frac{\mu+3}{8t}+\frac{(\mu-1)(\mu+15)}{2!(8t)^2}-\frac{(\mu-1)(\mu-9)(\mu+35)}{3!(8t)^3}+\cdots\biggr),\\
&& K'_\nu(t) \sim -\sqrt{\frac{\pi}{2t}}e^{-t} \biggl(1+\frac{\mu+3}{8t}+\frac{(\mu-1)(\mu+15)}{2!(8t)^2}+\frac{(\mu-1)(\mu-9)(\mu+35)}{3!(8t)^3}+\cdots\biggr).
\end{eqnarray*}

(viii) When $\nu$ is large, the following expansions hold uniformly with respect to $t$
\begin{eqnarray*}
&& I_\nu(\nu t) \sim \frac{1}{\sqrt{2\pi \nu}} \frac{e^{\nu \eta(t)}}{(1+t^2)^{\frac 14}}\biggl(1+\sum_{k=1}^\infty \frac{u_k(\tau(t))}{\nu^k}\biggr),\\
&& K_\nu(\nu t) \sim \sqrt{\frac{\pi}{2\nu}} \frac{e^{-\nu \eta(t)}}{(1+t^2)^{\frac 14}}\biggl(1+\sum_{k=1}^\infty (-1)^k \frac{u_k(\tau(t))}{\nu^k}\biggr),\\
&& I'_\nu(\nu t) \sim \frac{1}{\sqrt{2\pi \nu}} \frac{(1+t^2)^{\frac 14}}{t}e^{\nu \eta(t)}\biggl(1+\sum_{k=1}^\infty \frac{\upsilon_k(\tau(t))}{\nu^k}\biggr),\\
&& K'_\nu(\nu t) \sim -\sqrt{\frac{\pi}{2\nu}} \frac{(1+t^2)^{\frac 14}}{t}e^{-\nu \eta(t)}\biggl(1+\sum_{k=1}^\infty (-1)^k \frac{\upsilon_k(\tau(t))}{\nu^k}\biggr);
\end{eqnarray*}
where $\tau(t)=\ds\frac{1}{\sqrt{1+t^2}}$, $\eta(t)=\sqrt{1+t^2} + \ln {\ds\frac{t}{1+\sqrt{1+t^2}}}$, and
\begin{eqnarray*}
&&u_{k+1}(s)={\frac 12} s^2(1-s^2)u'_k(s)+{\frac 18}\int_0^s(1-5s^2)u_k(s)ds,\qquad k=0,1,\cdots,\\
&&\upsilon_k(s)=u_k(s)+s(s^2-1)(\frac 12 u_{k-1}(s)+su'_{k-1}(s)),\qquad k=1,2,\cdots,\\
&&u_0(s)=1.
\end{eqnarray*}

(ix) For $t_1\leqslant t_2$,  then $e^{\nu \eta(t_1)}e^{-\nu \eta(t_2)}\leqslant e^{-\nu(t_2 -t_1)}$,
where $\eta(t)$ has been  defined  in (viii).
}

\vskip 0.2 true cm
{\bf Proof.} (i)-(ii) can be found in Pages 204-206 and Pages 79 of [24],
and (vii)-(viii) can be seen from Pages 377-378 of [1].

We now show (iii). It follows from (i) that for $t>0$
\begin{eqnarray*}
 I_0(t)&=&\frac{e^t}{\sqrt{\pi} \Gamma(\frac 12)} \int_0^{1}e^{-2ts} (1-s)^{-\frac 12} s^{-\frac 12} ds\\
&\leqslant& \frac{e^t}{\sqrt{\pi} \Gamma(\frac 12)} \int_0^{1} (1-s)^{-\frac 12} s^{-\frac 12} ds\\
& = & e^t
\end{eqnarray*}
and
\begin{eqnarray*}
K_0(t)
&\leqslant& \frac{e^{-t}}{\sqrt{2t}} \int_0^{\infty} e^{-u} u^{-\frac 12} du \\
& =& \frac{\sqrt \pi e^{-t}}{\sqrt{2t}}.
\end{eqnarray*}
Thus, (iii) is proved.

Next, we start to prove (iv). Since $\nu>\frac 12$ and $t<1$, we have that from (i)
\begin{eqnarray*}
I_\nu(t)&\leqslant& \frac{(2t)^\nu e^t}{\sqrt{\pi}\Gamma(\nu+\frac 12)} \int_0^{1}u^{\nu-\frac 12}(1-u)^{\nu-\frac 12} du\\
&=& \frac{(2t)^\nu e^t}{\sqrt{\pi}\Gamma(\nu+\frac 12)} B(\nu+\frac 12, \nu+\frac 12)\\
&=& \ds\frac{e^t  (\frac t2)^\nu }{\Gamma(\nu+1)}.
\end{eqnarray*}
Similarly, we have that from (i)
\begin{eqnarray*}
K_\nu(t)&=& \frac{\sqrt \pi e^{-t}}{(2t)^\nu \Gamma(\nu+\frac 12)}\int_0^\infty e^{-u} u^{\nu-\frac 12} {(2t+u)}^{\nu-\frac 12} du\\
&\leqslant& \frac{\sqrt \pi e^{-t}}{(2t)^\nu \Gamma(\nu+\frac 12)}\int_0^\infty e^{-u} {(2t+u)}^{2\nu-1} du\\
&\leqslant& \frac{\sqrt \pi e^{t}}{(2t)^\nu} \frac{\Gamma(2\nu)}{\Gamma(\nu+\frac 12)}\\
&=& \frac{e^{t} \Gamma(\nu) 2^{\nu-1}}{t^\nu}.
\end{eqnarray*}
Thus, we complete the proof of (iv).

Next, we  prove (v). Since $t\geqslant 1$ and $\nu>\frac 12$, then by (i)
\begin{eqnarray*}
I_\nu(t)&=&\frac{e^t}{\sqrt{2\pi t} \Gamma(\nu+\frac 12)} \int_0^{2t}e^{-s}s^{\nu-\frac 12}(1-\frac{s}{2t})^{\nu-\frac 12} ds\\
&\leqslant& \frac{e^t}{\sqrt{2 \pi t} \Gamma(\nu+\frac 12)} \int_0^{2t}e^{-s}s^{\nu-\frac 12} ds\\
&\leqslant& \frac{e^t}{\sqrt{2 \pi t} \Gamma(\nu+\frac 12)} \Gamma(\nu+\frac 12)\\
&=& \frac{e^t}{\sqrt{2\pi t}}.
\end{eqnarray*}

We now show (vi). Due to $\frac 12<\nu\le M$ and  $t \geqslant 1$, we have from (i)
\begin{eqnarray*}
K_\nu(t)&=& \frac{\sqrt \pi e^{-t}}{\sqrt{2t} \Gamma(\nu+\frac 12)}\biggl(\int_0^{2t} e^{-u} u^{\nu-\frac 12} (1+\frac{u}{2t})^{\nu-\frac 12} du+\int_{2t}^{\infty} e^{-u} u^{\nu-\frac 12} (1+\frac{u}{2t})^{\nu-\frac 12} du\biggr)\\
&\leqslant& \frac{\sqrt \pi e^{-t}}{\sqrt{2t} \Gamma(\nu+\frac 12)} \biggl(2^{\nu-\frac 12} \int_0^{2t} e^{-u} u^{\nu-\frac 12} du+ (2t)^{\frac 12-\nu}\int_{2t}^{\infty} e^{-u} u^{\nu-\frac 12} (2t+u)^{\nu-\frac 12} du\biggr)\\
&\leqslant& \frac{\sqrt \pi e^{-t}}{\sqrt{2t} \Gamma(\nu+\frac 12)} \bigl(2^{\nu-\frac 12} \Gamma(\nu+\frac 12)+ t^{\frac 12-\nu} \Gamma(2\nu)\bigr)\\
& \leqslant& C_M \frac{\sqrt \pi e^{-t}}{\sqrt{2t}}.
\end{eqnarray*}

Finally, we prove (ix). By $\eta(t)=\sqrt{1+t^2} + \ln {\ds\frac{t}{1+\sqrt{1+t^2}}}$, then one has
$\eta'(t)=\ds\frac{\sqrt{1+t^2}}{t}.$
Thus, there exist a number $\xi \in (t_1, t_2)$ such that
$$
e^{\nu \eta(t_1)}e^{-\nu \eta(t_2)}=e^{-\nu\frac{\sqrt{1+\xi^2}}{\xi}(t_2-t_1)}
\leqslant e^{-\nu(t_2 -t_1)}.
$$
Collecting all the analysis above, we complete the proof of Lemma 2.3.\qed

Finally, we give an existence result of the supersonic solution to (1.5) and (1.13) in
the  domain which is larger than that of left hand side of the
shock surface $\Gamma$.

\vskip 0.2 true cm

{\bf Lemma 2.4.} {\it The equation (1.5) with the initial data
(1.13) has a $C^{\infty}$ solution $\vp^-(x)$ in the domain
$\O_-=\{x:  x_1\ge 0, x_2\in\Bbb R, x_3\ge\ds\f{s_0+b_0}{2}x_1\}$.  Moreover,
$\vp^-(x)-q_0x_1\in C^{\infty}(\O_-)$, $\vp^-(x)=\vp^-(x_1, x_2+2\pi, x_3)$, and there exists a
positive constant $C_k$ independent of $\ve$ such that
$$
\|\vp^-(x)-q_0x_1\|_{C^k(\O_-)}\leq C_k\ve \eqno(2.11)
$$
for any fixed $k\in\Bbb N$.}

{\bf Proof.}  We note that the equation (1.5) is quasi-linear
strictly hyperbolic with respect to the $x_1-$ direction for the
supersonic flow $\p_1\vp^->c^-$, furthermore, the initial condition
(1.13) is of a small perturbation. Thus, in
terms of the finite propagation property of the wave equation, the periodic property of the initial data
$(\vp_0(x_2,x_3), \vp_1(x_2,x_3))$ with respect to the variable $x_2$ and
the Picard iteration (or one can see
[16]), we know that Lemma 2.4
holds.\qed

{\bf Remark 2.1.} {\it By (2.11) and the standard extension theorem (see Theorem 7.25 of [15]),
we can extend the smooth function $\vp^-(x)$ in $\O_-$ into the whole domain $\O=\{x:  x_1\ge 0, x_2\in\Bbb R, x_3\ge b_0x_1\}$
such that the extension function $\t\vp^-(x)\in C^{\infty}(\bar\O)$ satisfies $\t\vp^-(x)=\vp^-(x)$ for $x\in\O_-$
and $\|\t\vp^-(x)-q_0x_1\|_{C^k(\O)}\leq C_k\ve$. For convenience, $\t\vp^-(x)$ will still be denoted by $\vp^-(x)$ later.
Here one should notice that $\t\vp^-(x)$ is not a solution to (1.5) in $\O\setminus\O_-$ in general case.}

\vskip 0.8 true cm
\centerline{\bf $\S 3$. Reformulation on (1.6)-(1.12) and  detailed descriptions on Theorem 1.1}
\vskip 0.5 true cm

By the notations in  (1.14)-(1.16) of $\S 1$,  it follows from Lemma 2.1 that the
function $\Phi(x)$ corresponding to the background solution is
$$
\Phi_0(x)=(q_0-u_{10}^+)x_1-u_{30}^+x_3=q_0x_1+O(\s)(x_1+x_3).\eqno(3.1)
$$
In this case, $\p_{x_1}\Phi_0(x)=q_0+O(q_0^{\f{\gamma-3}{\gamma-1}})>0$ and
$\p_{x_3}\Phi_0(x)=O(q_0^{\f{\gamma-3}{\gamma-1}})$
holds for large $q_0$ and $1<\g<3$. Thus, the transformation (1.15) is inverse since
$\p_{x_1}\Phi(x)$ and $\p_{x_3}\Phi(x)$ will be of the small perturbations of $\p_{x_1}\Phi_0(x)$
and $\p_{x_3}\Phi_0(x)$ respectively.
In addition, the corresponding unknown function $u(y)$ in
(1.16) for $\Phi_0(x)$ can be expressed as
$$
u_0(y)=y_1+O(\a)(y_1+y_3).\eqno(3.2)
$$
To solve (1.5)-(1.6) together with (1.7)-(1.13), we suffice to study the following problem
(one can also see (1.18) in $\S 1$)
$$
\begin{cases}
&L(u, \na_y u, \na_y^2u)\\
&\qquad \equiv \ss_{1\le i\le j\le 3} A_{ij}(u, \na_y u)\p_{y_iy_j}^2u+\ds\f{1}{q_0}\ss_{1\le i\le j\le 3} a_{ij}(\na_x\vp^--\na_x\Phi)\p_{x_ix_j}^2\vp^-=0\quad\qquad ~~\text{in}~~Q,\\
&G_1(u, \na_y u)=0\qquad\qquad\qquad\qquad\qquad\qquad\qquad\qquad\qquad\q\text{on}\qquad y_3=b_0y_1,\\
&G_2(u, \na_y u)=0\qquad\qquad\qquad\qquad\qquad\qquad\qquad\qquad\qquad\q\text{on}\qquad y_1=0,\\
&u(0,y_2,0)=0,\\
&u(y_1,y_2+2\pi,y_3)=u(y),\\
&\ds\lim_{y_1+y_3\rightarrow +\infty}\nabla_yu\quad\text{exists},
\end{cases}\eqno(3.3)
$$
where
\begin{eqnarray*}
&&G_1(u, \na_y u)\equiv -\big(1+b_0^2+\f{b_0}{q_0}(\p_{x_3}\vp^--b_0\p_{x_1}\vp^-)\big)\p_{y_3}u-\f{1}{q_0} (\p_{x_3}\vp^--b_0\p_{x_1}\vp^-)\p_{y_1}u-b_0,\\
&&G_2(u, \na_y u)\equiv \f{1}{q_0}(1-\f{\rho_-}{\rho_+})\p_{x_1}\vp^-\p_{y_1}u
+b_0\big(\f{1}{q_0}(1-\f{\rho_-}{\rho_+})\p_{x_1}\vp^--2\big)\p_{y_3}u \\
&&\q +\f{b_0}{q_0}(1-\f{\rho_-}{\rho_+})\p_{x_1}\vp^-\p_{y_1}u\p_{y_3}u
-(\p_{y_2}u)^2+\big(\f{b_0}{q_0}(1-\f{\rho_-}{\rho_+})(b_0\p_{x_1}\vp^--\p_{x_3}\vp^-)\\
&&\q -(1+b_0^2)\big)(\p_{y_3}u)^2-\f{1}{q_0}(1-\f{\rho_-}{\rho_+})\p_{x_2}\vp^-\p_{y_1}u\p_{y_2}u
-\f{b_0}{q_0}(1-\f{\rho_-}{\rho_+})\p_{x_2}\vp^-\p_{y_2}u\p_{y_3}u\\
&&\q -\f{1}{q_0}(1-\f{\rho_-}{\rho_+})\p_{x_3}\vp^-\p_{y_1}u\p_{y_3}u-1
\end{eqnarray*}
and
\begin{eqnarray*}
&&A_{11}=\f{1}{(\p_{y_1}u+b_0\p_{y_3}u)^3}\bigg(a_{11}(1+b_0\p_{y_3}u)^2+a_{22}(\p_{y_2}u)^2+a_{33}(\p_{y_3}u)^2 -2a_{12}\p_{y_2}u(1+b_0\p_{y_3}u)\\
&&\q\q-2a_{13}\p_{y_3}u(1+b_0\p_{y_3}u)+2a_{23}\p_{y_2}u\p_{y_3}u\bigg),\\
&&A_{22}=\f{a_{22}}{\p_{y_1}u+b_0\p_{y_3}u},\\
&&A_{33}=\f{1}{(\p_{y_1}u+b_0\p_{y_3}u)^3}\bigg(a_{11}b_0^2(1-\p_{y_1}u)^2+a_{22}b_0^2(\p_{y_2}u)^2+a_{33}(\p_{y_1}u)^2 -2a_{12}b_0^2\p_{y_2}u(1-\p_{y_1}u)\\
&&\q\q+2a_{13}b_0\p_{y_1}u(1-\p_{y_1}u)-2a_{23}b_0\p_{y_1}u\p_{y_2}u\bigg),\\
&&A_{12}=A_{21}=\f{1}{(\p_{y_1}u+b_0\p_{y_3}u)^2}\bigg(-a_{22}\p_{y_2}u+a_{12}(1+b_0\p_{y_3}u)-a_{23}\p_{y_3}u\bigg),\\
&&A_{13}=A_{31}=\f{1}{(\p_{y_1}u+b_0\p_{y_3}u)^3}\bigg(a_{11}b_0(1+b_0\p_{y_3}u)(1-\p_{y_1}u) +a_{22}b_0(\p_{y_2}u)^2 -a_{33}\p_{y_1}u\p_{y_3}u\\
&&\q\q +a_{12}b_0\p_{y_2}u(\p_{y_1}u-b_0\p_{y_3}u-2)+a_{13}(\p_{y_1}u+2b_0\p_{y_1}u\p_{y_3}u-b_0\p_{y_3}u) \\ &&\q\q-a_{23}\p_{y_2}u(\p_{y_1}u+b_0\p_{y_3}u)\bigg),\\
&&A_{23}=A_{32}=\f{1}{(\p_{y_1}u+b_0\p_{y_3}u)^2}\bigg(-a_{22}b_0\p_{y_2}u+a_{12}b_0(1-\p_{y_1}u)+a_{23}\p_{y_1}u\bigg).
\end{eqnarray*}
With respect to more precise properties of $G_1(u, \na_y u)$, $G_2(u, \na_y u)$ and $A_{ij}$, one can be referred in $\S 7$ below.
In addition, one should note that $\vp^-(x)$ in (3.3) has become a
function $\vp^-(u(y),y_2,y_3-b_0y_1+b_0u(y))$  depending on the unknown solution $u(y)$, and $\na_x\Phi(x)$
in (3.3) is also a function on $\na_yu$ by the transformations (1.15)-(1.16).
On the other hand, we especially point out that  the condition $u(0, y_2, 0)\equiv 0$ for all $y_2\in\Bbb R$ in (3.3)
comes from the attached shock property. Next we will show that (3.3) is overdetermined since (3.3)
can be solved as long as the condition $u(0, y_2, 0)\equiv 0$
in (3.3) is replaced by $u(0, y_2^0, 0)=0$ for any fixed $y_2^0$. Without loss of generality,
we assume $y_2^0=0$ and consider the following problem instead of (3.3)

\begin{equation}
\left\{
\begin{aligned}
&L(u, \na_y u, \na_y^2u)=0\qquad\qquad\qquad\qquad\qquad\qquad\qquad\qquad\qquad\quad ~~\text{in}~~Q,\\
&G_1(u, \na_y u)=0\qquad\qquad\qquad\qquad\qquad\qquad\qquad\qquad\q\text{on}\qquad y_3=b_0y_1,\\
&G_2(u, \na_y u)=0\qquad\qquad\qquad\qquad\qquad\qquad\qquad\qquad\q\text{on}\qquad y_1=0,\\
&u(0,0,0)=0,\\
&u(y_1,y_2+2\pi,y_3)=u(y),\\
&\ds\lim_{y_1+y_3\rightarrow +\infty}\nabla_yu\quad\text{exists}.
\end{aligned}
\right.\tag{3.4}
\end{equation}

With respect to the problem (3.4), we have

\vskip 0.2 true cm

{\bf Theorem 3.1.} {\it There exist some positive constants $\ve_0>0$, $0<\dl<1$, $0<\dl_0<1$,
and $C>0$ such that for any $\ve\in (0, \ve_0)$, the problem (3.4)
has a unique solution $u\in C^{6,\al}(Q)$ which fulfills
the following estimate
$$
\|u-u_0\|_{6,\alpha;Q}^{(-1-\dl,-\delta_0)}\le
C\ve,
$$
where $0<\al<1$.}

\vskip 0.2 true cm
{\bf Remark 3.1.} {\it  By Theorem 3.1, we have obtained the $C^{6,\al}-$regularities of
$u$ in the domain $Q$ due to the high regularity assumption on the supersonic incoming flow. On the
other hand, as in [27], by use of the separation variable method in $\S 4$ below, it seems that the
$C^{6,\al}-$regularities of $u$ in the interior of $Q$ are required to guarantee the convergence
of the appropriate solution sequence in the $C_{loc}^2-$space for the boundary value problem
of the Laplacian equation with two Neumann boundary conditions and a vanishing
condition of the first order derivatives at infinity. However, such high regularities
($C^{6,\al}(Q)-$regularities) are essentially unnecessary (only $C^{2,\al}(Q)-$regularities should
be enough) if we establish the $L_{loc}^2-$convergence of the appropriate solution sequence
instead of the $C_{loc}^2$ convergence in any compact subdomain of $Q$ and combine with some interior estimate techniques on
the second order linear elliptic equations in Chapter 6 of [15]. Since the requirements
on the higher order interior regularities are not essential for our problem, we omit
the related argument procedure on the reduction from $C^{6,\al}(Q)-$regularities to
$C^{2,\al}(Q)-$regularities in Theorem 3.1.}

\vskip 0.5 true cm
\centerline{\bf $\S 4$. On the linearization
of (3.4) and its related cut-off problem}
\vskip 0.5 true cm

In order to solve the  nonlinear problem (3.4), we first consider
its linearized case, which corresponds to an Neumann
boundary problem of a second order elliptic equation in an unbounded
angular domain. It will be seen that in terms of the smallness of $\ds\f{1}{q_0}$ and
Lemma 2.1, by a tedious but direct computation (see $\S 7$ below), the linearized
problem of (3.4) can be essentially expressed as

\begin{equation}
\left\{
\begin{aligned}
&\Delta\dot u=\dot f\q\q \text{in}\q\quad Q,\\
&\ds\f{\p \dot u}{\p n}=\dot g_2\q\q\text{on}\q \Si_1:y_3=b_0y_1,\\
&\ds\f{\p \dot u}{\p n}=\dot g_1\q\q\text{on}\q \Si_2:y_1=0,\\
&\dot u(y_1,y_2+2\pi,y_3)=\dot{u}(y),\\
&\dot{u}(0,0,0)=0,\\
&\ds\lim_{y_1+y_3\rightarrow\infty}\nabla\dot{u}=0,
\end{aligned}
\right.\tag{4.1}
\end{equation}
where $\dot f\in H_{4,\al}^{(1-\dl,2-\dl_0)}(Q)$ and $\dot g_i\in H_{5,\al}^{(-\dl,1-\dl_0)}(Q)$ for $0<\al<1$
and $i=1,2$.

Introducing the following cylindrical coordinate transformation
$$
y_1=r\cos\th,\q y_2=y_2,\q y_3=r\sin\th,
$$
where $r=\sqrt{y_1^2+y_3^2}$,  $\th\in[\th_0,\f{\pi}{2}]$, and $\th_0=arctan b_0$. Then (4.1) can be changed as

\begin{equation}
\left\{
\begin{aligned}
&\Delta\dot u=\p_r^2\dot{u}+r^{-2}\p_{\th}^2\dot{u}+\p_{y_2}^2\dot{u}+r^{-1}\p_r\dot{u}=\dot f\q\q \text{in}\q Q,\\
&-\p_\th\dot u=r\dot g_1\q\q\text{on}\q \Si_1:\th=\th_0,\\
&\p_\th\dot u=r\dot g_2\q\q\text{on}\q \Si_2:\th=\f{\pi}{2},\\
&\dot u(r,\th, y_2+2\pi)=\dot{u}(r,\th,y_2),\\
&\dot{u}(0,\th, 0)=0,\\
&\ds\lim_{r\rightarrow\infty}\nabla\dot{u}=0,
\end{aligned}
\right.\tag{4.2}
\end{equation}
where $Q=\{(r,\th,y_2): r\in \Bbb R^+,\th\in(\th_0,\f{\pi}{2}), y_2\in \Bbb T\}$
under the  cylindrical coordinate transformation.

Let
$h(r,\th,y_2)=-\ds\f{r}{2(\f{\pi}{2}-\th_0)}(\dot g_1+\dot g_2)(\th-\th_0)^2+r\dot g_1\th$
and $v(r,\th,y_2)=\dot u(r,\th,y_2)-h(r,\th,y_2)$, then the problem (4.2) can be changed as
\begin{equation}
\left\{
\begin{aligned}
&\Delta v=\p_r^2v+r^{-2}\p_{\th}^2v+\p_{y_2}^2v+r^{-1}\p_rv=f\equiv \dot f-\Delta h\q\q \text{in}\q Q,\\
&\p_\th v=0\qquad \qquad\qquad\qquad\q\q\q\text{on}\q \Si_1,\\
&\p_\th v=0\qquad\qquad\qquad\qquad\q\q\q\text{on}\q \Si_2,\\
&v(r,\th, y_2+2\pi)=v(r,\th, y_2),\\
&v(0,\th,0)=0,\\
&\ds\lim_{r\rightarrow\infty}\nabla v=0,
\end{aligned}
\right.\tag{4.3}
\end{equation}
where $f\in H_{4,\al}^{(1-\dl,2-\dl_0)}(Q)$.

In order to solve the unbounded domain problem (4.3),
we will  consider the following cut-off problem
in the bounded domain:
\begin{equation}
\left\{
\begin{aligned}
&\Delta v_L=\p_r^2v_L+r^{-2}\p_{\th}^2v_L+\p_{y_2}^2v_L+r^{-1}\p_rv_L=f_L\\
&\qquad \quad\qquad \quad\qquad \quad\text{in}~ Q_L=\{(r,\th,y_2):0<r<L,(r,\th,y_2)\in Q\},\\
&\p_\th v_L=0\q\q\text{on}\q \Si_{1,L}=\{(r,y_2):0<r<L,(r,y_2)\in \Si_1\},\\
&\p_\th v_L=0\q\q\text{on}\q \Si_{2,L}=\{(r,y_2):0<r<L,(r,y_2)\in \Si_2\},\\
&\p_rv_L=c_L\q\q\text{on}\q \Si_{3,L}=\{(r,\th,y_2):r=L,\th_0\le\th\le\f{\pi}{2},y_2\in \Bbb T\},\\
&v_L(r,\th,y_2+2\pi)=v_L(r,\th,y_2),\\
&v_L(0,\th,0)=0,
\end{aligned}
\right.\tag{4.4}
\end{equation}
where $L\ge 4$, $f_L$ is the restriction of $f$ on $Q_L$, which obviously obeys
$$
||f_L||_{4,\al;Q_L}^{(1-\dl,2-\dl_0)}\le  ||f||_{4,\al;Q}^{(1-\dl,2-\dl_0)}.\eqno(4.5)
$$
In addition, in order to guarantee the solvability of (4.4), we require to choose the constant
$c_L$ in (4.4) such that
$$
\int_{Q_L}f_Ldy=\int_{\Si_{3,L}}c_L dS.\eqno(4.6)
$$
From (4.6), we can arrive at
$$
|c_L|\le CL^{\dl_0-1}||f_L||_{4,\al;Q_L}^{(1-\dl,2-\dl_0)}\le CL^{\dl_0-1}||f||_{4,\al;Q}^{(1-\dl,2-\dl_0)}\rightarrow0
\qquad \text{as $L\rightarrow\infty$}.\eqno(4.7)
$$
With respect to the linear problem (4.4), we have

\vskip 0.2 true cm

{\bf Proposition 4.1.} {\it There exists a unique solution $v_L\in C^2(\bar{Q}_L\backslash\{(0,y_2,0)
:y_2\in \Bbb R\})$ to (4.4) such that
$$
||v_L||_{0,0;Q_L}^{(-1-\dl,-\dl_0)}\le C||f_L||_{4,\al;Q_L}^{(1-\dl,2-\dl_0)},\eqno(4.8)
$$
where the constant $C>0$ is independent of $L$.}

{\bf Proof.}  We will divide the proof of Proposition 4.1  into the following three steps.

\vskip 0.2 true cm

{\bf Step 1. Existence of a formal solution $v_L$ to (4.4)}
\vskip 0.2 true cm

We will use the separation variable method to solve (4.4). To this end, as in [29], we first focus on the corresponding homogeneous equation
of (4.4). Consider the nontrivial solutions to the following problem
\begin{equation}
\left\{
\begin{aligned}
&\Delta v=\p_r^2 v+r^{-2}\p_{\theta}^2 v+\p_{y_2}^2 v+r^{-1}\p_{r} v=0,\\
&\text{$\p_{\theta} v=0$ \qquad  on \qquad $\Sigma_{1,L} \cup \Sigma_{2,L}$},\\
&v(r,\theta,y_2)=v(r,\theta,y_2+2\pi).
\end{aligned}
\right.\tag{4.9}
\end{equation}

Set $v(r,\theta,y_2)=R(r) \Theta(\theta) Y(y_2)$, then we have
\begin{equation}
\left\{
\begin{aligned}
&Y^{''}(y_2)+\lambda Y(y_2)=0,\\
&Y(y_2)=Y(y_2+2\pi),
\end{aligned}
\right.\tag{4.10}
\end{equation}
and
\begin{equation}
\left\{
\begin{aligned}
&\Theta^{''}(\theta)+\mu \Theta(\theta)=0,\\
&\Theta^{'}(\theta_0)=\Theta^{'}(\frac{\pi}2)=0,
\end{aligned}
\right.\tag{4.11}
\end{equation}
and
\begin{equation}
\left\{
\begin{aligned}
&r^2 R^{''}(r)+r R^{'}(r)-(\mu+\lambda r^2)R(r)=0,\\
&R^{'}(L)=0,\qquad \text{$R(0)$\quad is  bounded},
\end{aligned}
\right.\notag
\end{equation}
where $\lambda \in \rr$ and $\mu \in \rr$.

We can get that the eigenvalues of (4.10) and (4.11) are $\lambda_n=n^2 (n=0,1,\cdots)$ and $\mu_m=(\ds\frac{m\pi}{\frac{\pi}2-\theta_0})^2 (m=0,1,\cdots)$, whose corresponding eigenfunctions are $\{\sin(ny_2),\cos(ny_2)\}_{n=0}^{\infty}$ and $\{\cos \sqrt{\mu_m}(\theta-\theta_0)\}_{m=0}^{\infty}$ respectively.

We now solve equation (4.4) by use of the eigenfunction expansion method in terms of the complete
orthogonal basis $\{\cos \sqrt{\mu_m}(\theta-\theta_0)\sin(ny_2),\cos \sqrt{\mu_m}(\theta-\theta_0)\cos(ny_2)\}_{m,n=0}^{\infty}$.

Let

\begin{align}
v_L(r,&\theta,y_2) \notag\\
&=R_{00}(r)+\sum_{m=1}^{\infty} R_{m0}(r) \cos \sqrt{\mu_m}(\theta-\theta_0)+\sum_{n=1}^{\infty}\Big(R_{0n}^{(1)}(r)\sin ny_2+R_{0n}^{(2)}(r)\cos ny_2\Big) \notag\\
&+\sum_{m=1}^{\infty}\sum_{n=1}^{\infty}\Big(R_{mn}^{(1)}(r)\cos \sqrt{\mu_m}(\theta-\theta_0)\sin ny_2+R_{mn}^{(2)}(r)\cos \sqrt{\mu_m}(\theta-\theta_0)\cos ny_2\Big) \tag{4.12}
\end{align}

and

\begin{align}
{f_L}(r,&\theta,y_2) \notag\\
&={f_L}_{00}(r)+\sum_{m=1}^{\infty} {f_L}_{m0}(r) \cos \sqrt{\mu_m}(\theta-\theta_0)+\sum_{n=1}^{\infty}\Big({f_L}_{0n}^{(1)}(r)\sin ny_2+{f_L}_{0n}^{(2)}(r)\cos ny_2\Big) \notag\\
&+\sum_{m=1}^{\infty}\sum_{n=1}^{\infty}\Big({f_L}_{mn}^{(1)}(r)\cos \sqrt{\mu_m}(\theta-\theta_0)\sin ny_2+{f_L}_{mn}^{(2)}(r)\cos \sqrt{\mu_m}(\theta-\theta_0)\cos ny_2\Big),\tag{4.13}
\end{align}

where
\begin{eqnarray*}
&& {f_L}_{00}(r)=\frac{1}{2\pi (\frac{\pi}2-\theta_0)}\int_{\theta_0}^{\frac{\pi}2}\int_{0}^{2\pi}{f_L}(r,\theta,y_2)dy_2 d \theta, \\
&& {f_L}_{m0}(r)=\frac{1}{\pi (\frac{\pi}2-\theta_0)}\int_{\theta_0}^{\frac{\pi}2}\int_{0}^{2\pi}{f_L}(r,\theta,y_2)\cos \sqrt{\mu_m}(\theta-\theta_0)dy_2 d \theta, \\
&& {f_L}_{0n}^{(1)}(r)=\frac{1}{\pi (\frac{\pi}2-\theta_0)}\int_{\theta_0}^{\frac{\pi}2}\int_{0}^{2\pi}{f_L}(r,\theta,y_2)\sin ny_2dy_2 d \theta, \\
&& {f_L}_{0n}^{(2)}(r)=\frac{1}{\pi (\frac{\pi}2-\theta_0)}\int_{\theta_0}^{\frac{\pi}2}\int_{0}^{2\pi}{f_L}(r,\theta,y_2)\cos ny_2dy_2 d \theta, \\
&& {f_L}_{mn}^{(1)}(r)=\frac{2}{\pi (\frac{\pi}2-\theta_0)}\int_{\theta_0}^{\frac{\pi}2}\int_{0}^{2\pi}{f_L}(r,\theta,y_2)\cos \sqrt{\mu_m}(\theta-\theta_0) \sin ny_2 dy_2 d \theta, \\
&& {f_L}_{mn}^{(2)}(r)=\frac{2}{\pi (\frac{\pi}2-\theta_0)}\int_{\theta_0}^{\frac{\pi}2}\int_{0}^{2\pi}{f_L}(r,\theta,y_2)\cos \sqrt{\mu_m}(\theta-\theta_0) \cos ny_2 dy_2 d \theta.
\end{eqnarray*}

Substituting (4.12) and (4.13) into the equation $\Delta v_L={f_L}$ yields
\begin{equation}
\left\{
\begin{aligned}
&R_{00}^{''}(r)+r^{-1}R_{00}^{'}(r)={f_L}_{00}(r),\\
&R_{m0}^{''}(r)-r^{-2}\mu_m R_{m0}(r)+r^{-1}R_{m0}^{'}(r)={f_L}_{m0}(r),\quad m\geqslant 1,\\
&(R_{0n}^{(i)})^{''}(r)-n^{2} (R_{0n}^{(i)})(r)+r^{-1}(R_{0n}^{(i)})^{'}(r)={f_L}_{0n}^{(i)}(r),\quad n\geqslant 1,\qquad i=1,2,\\
&(R_{mn}^{(i)})^{''}(r)-(r^{-2}\mu_m+n^{2}) (R_{mn}^{(i)})(r)+r^{-1}(R_{mn}^{(i)})^{'}(r)={f_L}_{mn}^{(i)}(r),
\quad m,n\geqslant 1,\quad i=1,2.\\
\end{aligned}
\right.\tag{4.14}
\end{equation}
Meanwhile, substituting (4.12) into the condition $v_L(0,\theta,0)=0$ yields
$$
\Big(R_{00}(0)+\sum_{n=1}^{\infty}R_{0n}^{(2)}(0)\Big)+\sum_{m=1}^{\infty}\Big(R_{m0}(0)+\sum_{n=1}^{\infty}R_{mn}^{(2)}(0)\Big)\cos \sqrt{\mu_m}(\theta-\theta_0)=0.
$$
By the orthogonality of $\{\cos \sqrt{\mu_m}(\theta-\theta_0)\}_{m=0}^{\infty}$, then
$$
R_{m0}(0)+\sum_{n=1}^{\infty}R_{mn}^{(2)}(0)=0,\qquad m=0,1,2,\cdots.
\eqno(4.15)
$$
In addition, it follows from $\p_{r} v_L=c_L$ on $\Sigma_{3,L}$ that
\begin{eqnarray*}
&&R_{00}^{'}(L)+\sum_{m=1}^{\infty} R_{m0}^{'}(L) \cos \sqrt{\mu_m}(\theta-\theta_0)+\sum_{n=1}^{\infty}\Big((R_{0n}^{(1)})^{'}(L)\sin ny_2+(R_{0n}^{(2)})^{'}(L)\cos ny_2\Big)\\
&&\quad +\sum_{m=1}^{\infty}\sum_{n=1}^{\infty}\Big((R_{mn}^{(1)})^{'}(L)\cos \sqrt{\mu_m}(\theta-\theta_0)\sin ny_2+(R_{mn}^{(2)})^{'}(L)\cos \sqrt{\mu_m}(\theta-\theta_0)\cos ny_2\Big)\\
&&=c_L.
\end{eqnarray*}
This derives
\begin{equation}
\left\{
\begin{aligned}
&R'_{00}(L)=c_L,~~R_{m0}^{'}(L)=0,~ m=1,2,\cdots,~~\\
&(R_{mn}^{(1)})^{'}(L)=(R_{mn}^{(2)})^{'}(L)=0,~  m=0,1,2,\cdots, n=1,2,\cdots
\end{aligned}
\right.\tag{4.16}
\end{equation}

Collecting (4.14)-(4.16), we obtain the following systems

\begin{equation}
\left\{
\begin{aligned}
&R_{00}^{''}(r)+r^{-1}R_{00}^{'}(r)={f_L}_{00}(r),\\
&R_{00}^{'}(L)=c_L,\\
&R_{00}(0)=-\ds\sum_{n=1}^{\infty}R_{0n}^{(2)}(0),
\end{aligned}
\right.\tag{4.17}
\end{equation}

\begin{equation}
\left\{
\begin{aligned}
&R_{m0}^{''}(r)-r^{-2}\mu_m R_{m0}(r)+r^{-1}R_{m0}^{'}(r)={f_L}_{m0}(r),\quad m=1,2,\cdots,\\
&R_{m0}^{'}(L)=0,\\
&R_{m0}(0)=-\ds\sum_{n=1}^{\infty}R_{mn}^{(2)}(0),\\
&R_{m0}(0) \quad is \quad bounded,
\end{aligned}
\right.\tag{4.18}
\end{equation}

\begin{equation}
\left\{
\begin{aligned}
&(R_{0n}^{(i)})^{''}(r)-n^{2} R_{0n}^{(i)}(r)+r^{-1}(R_{0n}^{(i)})^{'}(r)={f_L}_{0n}^{(i)}(r),\quad n=1,2,\cdots,\\
&(R_{0n}^{(i)})^{'}(L)=0,\\
&R_{0n}^{(i)}(0)\quad is \quad bounded,
\end{aligned}
\right.\tag{4.19}
\end{equation}
and
\begin{equation}
\left\{
\begin{aligned}
&(R_{mn}^{(i)})^{''}(r)-(r^{-2}\mu_m+n^{2}) R_{mn}^{(i)}(r)+r^{-1}(R_{mn}^{(i)})^{'}(r)={f_L}_{mn}^{(i)}(r),\quad m,n=1,2,\cdots,\\
&(R_{mn}^{(i)})^{''}(L)=0,\\
&R_{mn}^{(i)}(0) \quad is \quad bounded.
\end{aligned}
\right.\tag{4.20}
\end{equation}

Notice that the general solutions to the equations in (4.17)-(4.20) are respectively
\begin{equation}
\left\{
\begin{aligned}
&R_{00}(r)=C_{00}^1+C_{00}^2 \ln r+\int_0^r \eta^{-1}\int_0^{\eta} \xi {f_L}_{00}(\xi)d \xi d \eta,\\
&R_{m0}(r)=C_m^1 r^{\sqrt{\mu_m}} +C_m^2 r^{-\sqrt{\mu_m}}+r^{-\sqrt{\mu_m}}\int_r^L \eta^{2\sqrt{\mu_m}-1}\int_{\eta}^L \xi^{1-\sqrt{\mu_m}}{f_L}_{m0}(\xi) d \xi d\eta,\\
&\qquad \quad\qquad \quad \qquad \quad \qquad \quad \qquad \quad \qquad \quad \qquad \quad \qquad \quad  m=1,2,\cdots,\\
&R_{0n}^{(i)}(r)=C_{0n}^{i1}I_{0}(nr)+C_{0n}^{i2}K_{0}(nr)-I_{0}(nr)\int_r^L sK_{0}(ns){f_L}_{0n}^{(i)}(s) d s\\
& \qquad \qquad -K_{0}(nr)\int_0^r sI_{0}(ns){f_L}_{0n}^{(i)}(s) d s,\quad i=1,2,\qquad n=1,2,\cdots,\\
&R_{mn}^{(i)}(r)=C_{mn}^{i1}I_{\sqrt{\mu_m}}(nr)+C_{mn}^{i2}K_{\sqrt{\mu_m}}(nr)-I_{\sqrt{\mu_m}}(nr)\int_r^L sK_{\sqrt{\mu_m}}(ns){f_L}_{mn}^{(i)}(s) d s\\
& \qquad \qquad -K_{\sqrt{\mu_m}}(nr)\int_0^r sI_{\sqrt{\mu_m}}(ns){f_L}_{mn}^{(i)}(s) d s,\quad i=1,2,\qquad m,n=1,2,\cdots,
\end{aligned}
\right.\tag{4.21}
\end{equation}

Next, we determine the constants in (4.21) by the boundary conditions in (4.17)-(4.20).
At first, due to $R_{00}^{'}(r)=C_{00}^2 r^{-1}+r^{-1}\int_0^{r} \xi {f_L}_{00}(\xi)d \xi$
and $R_{00}^{'}(L)=c_L$, one has
$$
C_{00}^2=L c_L-\int_0^{L} \xi {f_L}_{00}(\xi)d \xi.\eqno(4.22)
$$
On the other hand, it follows from  the compatibility condition (4.6) that
$$
c_L=\f{1}{2\pi L(\f{\pi}{2}-\th_0)}\int_{Q_L}f_L(y)dy.
$$
This, together with (4.13) and (4.22), yields $C_{00}^2=0$ holds. It is noted that $C_{00}^1=R_{00}(0)=-\ds\sum_{n=1}^{\infty}R_{0n}^{(2)}(0)$,
then the solution of (4.17) can be expressed as
$$
R_{00}(r)=-\ds\sum_{n=1}^{\infty}R_{0n}^{(2)}(0)+\int_0^r \eta^{-1}\int_0^{\eta} \xi {f_L}_{00}(\xi)d \xi d \eta,\eqno(4.23)
$$
here we especially point out that the numbers $R_{0n}^{(2)}(0)$ ($n=1,2, ...$) in (4.23) have not been known yet.

Secondly, by the boundedness of $R_{m0}(0)$ for $m\in\Bbb N$, then
$$
C_{m0}^2=-\int_0^L \eta^{2\sqrt{\mu_m}-1}\int_{\eta}^L \xi^{1-\sqrt{\mu_m}}{f_L}_{m0}(\xi) d \xi d\eta.
$$
In this case,
$$
R_{m0}(r)=C_{m0}^1 r^{\sqrt{\mu_m}} -r^{-\sqrt{\mu_m}}\int_0^r \eta^{2\sqrt{\mu_m}-1}\int_{\eta}^L \xi^{1-\sqrt{\mu_m}}{f_L}_{m0}(\xi) d \xi d\eta
$$
and
$$
R_{m0}^{'}(L)=C_{m0}^1\sqrt{\mu_m} L^{\sqrt{\mu_m}-1}+\sqrt{\mu_m}L^{-\sqrt{\mu_m}-1}\int_0^L \eta^{2\sqrt{\mu_m}-1}\int_{\eta}^L \xi^{1-\sqrt{\mu_m}}{f_L}_{m0}(\xi) d \xi d\eta.
$$
Together with the boundary condition $R_{m0}^{'}(L)=0$ in (4.18), this yields
$$
C_{m0}^1=-L^{-2\sqrt{\mu_m}}\int_0^L \eta^{2\sqrt{\mu_m}-1}\int_{\eta}^L \xi^{1-\sqrt{\mu_m}}{f_L}_{m0}(\xi) d \xi d\eta.
$$
Consequently, the solution of (4.18) has the following expression

\begin{align}
R_{m0}(r)&=-r^{\sqrt{\mu_m}}L^{-2\sqrt{\mu_m}}\int_0^L \eta^{2\sqrt{\mu_m}-1}\int_{\eta}^L \xi^{1-\sqrt{\mu_m}}{f_L}_{m0}(\xi) d \xi d\eta \notag\\
&\qquad -r^{-\sqrt{\mu_m}}\int_0^r \eta^{2\sqrt{\mu_m}-1}\int_{\eta}^L \xi^{1-\sqrt{\mu_m}}{f_L}_{m0}(\xi) d \xi d\eta.\tag{4.24}
\end{align}

Thirdly, we solve (4.19). By the boundedness of $R_{0n}^{(i)}(0)$ for $i=1,2$ and the properties of Bessel functions as $r\rightarrow 0$
in Lemma 2.3, we can get from the expression of $R_{0n}^{(i)}(r)$
$$
C_{0n}^{i2}=\lim_{r\rightarrow 0}\int_0^r sI_{0}(ns){f_L}_{0n}^{(i)}(s) d s=0.
$$
In addition, a simple computation shows
$$
\Big(R_{0n}^{(i)}\Big)^{'}(r)=C_{0n}^{i1}nI_{0}^{'}(nr)-nI_{0}^{'}(nr)\int_r^L sK_{0}(ns){f_L}_{0n}^{(i)}(s) d s-nK_{0}^{'}(nr)\int_0^r sI_{0}(ns){f_L}_{0n}^{(i)}(s) d s.
$$
This, together with the boundary condition $(R_{0n}^{(i)})^{'}(L)=0$, yields
$$
C_{0n}^{i1}=\frac{K_{0}^{'}(nL)}{I_{0}^{'}(nL)}\int_0^L sI_{0}(ns){f_L}_{0n}^{(i)}(s) d s.
$$
Thus, the solution of (4.19) is
\begin{align}
R_{0n}^{(i)}(r)&=\frac{K_{0}^{'}(nL)}{I_{0}^{'}(nL)}I_{0}(nr)\int_0^L sI_{0}(ns){f_L}_{0n}^{(i)}(s) d s-I_{0}(nr)\int_r^L sK_{0}(ns){f_L}_{0n}^{(i)}(s) d s \notag\\
&-K_{0}(nr)\int_0^r sI_{0}(ns){f_L}_{0n}^{(i)}(s) d s.\tag{4.25}
\end{align}

Finally, we solve (4.20). It follows from $I_{\m}(0)=0$, $K_{\m}(0)=\infty$, the boundedness of $R_{mn}^{(i)}(0)$
and the expression of $R_{mn}^{(i)}(r)$ in (4.21) that
$$
C_{mn}^{i2}=0.
$$
Due to $(R_{mn}^{(i)}\big)^{'}(L)=0$ and
$$
\big(R_{mn}^{(i)}\big)^{'}(L)=C_{mn}^{i1}nI_{\m}^{'}(nL)-nK_{\m}^{'}(nL)\int_0^L sI_{\m}(ns){f_L}_{mn}^{(i)}(s) d s,
$$
one has
$$
C_{mn}^{i1}=\frac{K_{\m}^{'}(nL)}{I_{\m}^{'}(nL)}\int_0^L sI_{\m}(ns){f_L}_{mn}^{(i)}(s) d s.
$$
Thus, we can obtain the solution to (4.20) as follows
\begin{align}
R_{mn}^{(i)}(r)&=I_{\m}(nr)\Big(\frac{K_{\m}^{'}(nL)}{I_{\m}^{'}(nL)}\int_0^L sI_{\m}(ns){f_L}_{mn}^{(i)}(s) d s-\int_r^L sK_{\m}(ns){f_L}_{mn}^{(i)}(s) d s\Big) \notag\\
&-K_{\m}(nr)\int_0^r sI_{\m}(ns){f_L}_{mn}^{(i)}(s) d s.\tag{4.26}
\end{align}

Collecting (4.23)-(4.26), the formal solution of (4.4) can be expressed as
$$
v_L(r,\theta,y_2)=R_{00}(r)+I_1+I_2+I_3,\eqno(4.27)
$$
where
\begin{equation}
\left\{
\begin{aligned}
&I_1=\sum\limits_{m=1}^{\infty} R_{m0}(r) \cos \m(\theta-\theta_0),\\
&I_2=\sum\limits_{n=1}^{\infty}\Big(R_{0n}^{(1)}(r)\sin ny_2+R_{0n}^{(2)}(r)\cos ny_2\Big),\\
&I_3=\sum\limits_{m=1}^{\infty}\sum\limits_{n=1}^{\infty}\Big(R_{mn}^{(1)}(r)\cos \m(\theta-\theta_0)\sin ny_2+R_{mn}^{(2)}(r)\cos \m(\theta-\theta_0)\cos ny_2\Big).
\end{aligned}
\right.\tag{4.28}
\end{equation}

{\bf Step 2. The uniform convergence of (4.27)}
\vskip 0.2 true cm
In order to show that $v_L(r,\theta,y_2)$ in (4.27) is a real solution to (4.4), we require to give a more precise estimates on $f_L$.
Since $f_L\in H_{4,\alpha}^{(1-\delta,2-\delta_0)}$ and $f_L(r,\theta,y_2+2\pi)=f_L(r,\theta,y_2)$, which means that $\p_{r,\theta,y_2}^4 f_L(r,\theta,y_2)$ is continuous for the variables $(r,\theta,y_2)\in (0,L] \times [\theta_0,\frac{\pi}2]\times [0,2\pi]$, then we can use the integration by parts to obtain that for $n,m\geqslant 1$,

\begin{align}
{f_L}_{m0}(r)&=-\ds\f{(-1)^m}{\mu_m\pi(\f{\pi}{2}-\th_0)}\int_0^{2\pi}\bigl(\p_{\th}f_L(r,\f{\pi}{2}, y_2)-\p_{\th}f_L(r,\th_0, y_2)\bigr)dy_2 \notag\\
&-\frac{1}{{\mu_m}\pi (\frac{\pi}2-\theta_0)} \int_{\theta_0}^{\frac{\pi}2}\int_{0}^{2\pi}\p_{\theta}^2 {f_L}(r,\theta,y_2) \cos \m(\theta-\theta_0) d\theta dy_2,\tag{4.29}\\
{f_L}_{0n}^{(1)}(r)&=\frac{1}{n^4 \pi (\frac{\pi}2-\theta_0)}\int_{\theta_0}^{\frac{\pi}2} \int_{0}^{2\pi}\p_{y_2}^4 {f_L}(r,\theta,y_2) \sin ny_2 dy_2 d \theta,\tag{4.30}\\
{f_L}_{mn}^{(1)}(r)&=-\ds\f{2(-1)^m}{n^2\mu_m\pi(\f{\pi}{2}-\th_0)}\int_0^{2\pi}\bigl(\p_{\th}\p_{y_2}^2f_L(r,\f{\pi}{2}, y_2)
-\p_{\th}\p_{y_2}^2f_L(r,\th_0, y_2)\bigr)sin y_2dy_2 \notag\\
&+\frac{2}{n^2 {\mu_m}\pi (\frac{\pi}2-\theta_0)} \int_{\theta_0}^{\frac{\pi}2} \int_{0}^{2\pi}\p_{\theta}^2\p_{y_2}^2 {f_L}(r,\theta,y_2) \sin ny_2 \cos \m(\theta-\theta_0)dy_2 d \theta.\tag{4.31}
\end{align}

From (4.29)-(4.31), we can derive that  for $0<r\leqslant 1$
$$
|{f_L}_{00}(r)| \leqslant C\|f\|_{0,\alpha;Q}^{(1-\delta,2-\delta_0)} r^{\delta-1}\eqno(4.32)
$$
and
\begin{equation}
\left\{
\begin{aligned}
&|{f_L}_{m0}(r)| \leqslant C\|f\|_{0,\alpha;Q}^{(1-\delta,2-\delta_0)} r^{\delta-1},\\
&|{f_L}_{m0}(r)| \leqslant C\mu_m^{-\frac 12}\|f\|_{1,\alpha;Q}^{(1-\delta,2-\delta_0)} r^{\delta-1},\\
&|{f_L}_{m0}(r)| \leqslant C\mu_m^{-1}\|f\|_{2,\alpha;Q}^{(1-\delta,2-\delta_0)} r^{\delta-1}
\end{aligned}
\right.\tag{4.33}
\end{equation}

and
$$
|{f_L}_{0n}^{(1)}(r)| \leqslant Cn^{-k}\|f\|_{k,\alpha;Q}^{(1-\delta,2-\delta_0)} r^{\delta-1-k},\quad
k=0,1,2,3,\eqno(4.34)
$$
and
\begin{equation}
\left\{
\begin{aligned}
&|{f_L}_{mn}^{(1)}(r)| \leqslant Cn^{-k}\|f\|_{k,\alpha;Q}^{(1-\delta,2-\delta_0)} r^{\delta-1-k},\quad k=0,1,2,\\
&|{f_L}_{mn}^{(1)}(r)| \leqslant Cn^{-1}\mu_m^{-1}\|f\|_{3,\alpha;Q}^{(1-\delta,2-\delta_0)} r^{\delta-2},\\
&|{f_L}_{mn}^{(1)}(r)| \leqslant Cn^{-2}\mu_m^{-1}\|f\|_{4,\alpha;Q}^{(1-\delta,2-\delta_0)} r^{\delta-3};
\end{aligned}
\right.\tag{4.35}
\end{equation}
for $r>1$
$$
|{f_L}_{00}(r)| \leqslant C\|f\|_{0,\alpha;Q}^{(1-\delta,2-\delta_0)} r^{\delta_0-1}\eqno(4.36)
$$
and
\begin{equation}
\left\{
\begin{aligned}
&|{f_L}_{m0}(r)|  \leqslant C\|f\|_{0,\alpha;Q}^{(1-\delta,2-\delta_0)} r^{\delta_0-2},\\
&|{f_L}_{m0}(r)| \leqslant C \mu_m^{-\frac 12}\|f\|_{1,\alpha;Q}^{(1-\delta,2-\delta_0)} r^{\delta_0-2},\\
&|{f_L}_{m0}(r)| \leqslant C\mu_m^{-1}\|f\|_{2,\alpha;Q}^{(1-\delta,2-\delta_0)} r^{\delta_0-2}
\end{aligned}
\right.\tag{4.37}
\end{equation}
and
$$
|{f_L}_{0n}^{(1)}(r)| \leqslant Cn^{-k}\|f\|_{k,\alpha;Q}^{(1-\delta,2-\delta_0)} r^{\delta_0-2-k},\quad k=0,1,2,3,
\eqno(4.38)
$$
and
\begin{equation}
\left\{
\begin{aligned}
&|{f_L}_{mn}^{(1)}(r)| \leqslant Cn^{-k}\|f\|_{k,\alpha;Q}^{(1-\delta,2-\delta_0)} r^{\delta_0-2-k},\quad k=0,1,2,\\
&|{f_L}_{mn}^{(1)}(r)| \leqslant Cn^{-1}\mu_m^{-1}\|f\|_{3,\alpha;Q}^{(1-\delta,2-\delta_0)} r^{\delta_0-3},\\
&|{f_L}_{mn}^{(1)}(r)| \leqslant Cn^{-2}\mu_m^{-1}\|f\|_{4,\alpha;Q}^{(1-\delta,2-\delta_0)} r^{\delta_0-4}.
\end{aligned}
\right.\tag{4.39}
\end{equation}

We now start to show that the series in (4.28) are convergent for $(r,\theta,y_2)\in (0,L]\times [\theta_0,\frac{\pi}{2}]\times [0,2\pi]$.
In fact, by Lemma A.1-Lemma A.4 in Appendix
(which are based on (4.32)-(4.39)), one has

\begin{align}
&\|I_1\|_{0,\alpha;Q_L}^{(0,-\delta_0)}\leqslant  C \|f\|_{0,\alpha;Q}^{(1-\delta,2-\delta_0)}.\tag{4.40}\\
&\|I_2\|_{0,0;Q_L}^{(0,0)}\leqslant C \|f\|_{1,\alpha;Q}^{(1-\delta,2-\delta_0)},\tag{4.41}\\
&\|R_{00}\|_{0,0;Q_L}^{(0,-\delta_0)} \leqslant C \|f\|_{1,\alpha;Q}^{(1-\delta,2-\delta_0)},\tag{4.42}\\
&\|I_3\|_{0,0;Q_L}^{(0,0)} \leqslant C\|f\|_{1,\alpha;Q}^{(1-\delta,2-\delta_0)},\tag{4.43}
\end{align}
where $C>0$  is independent of $L$.

Thus, combining (4.40)-(4.43) yields the uniform convergence of (4.27) for $(r,\theta,y_2)\in (0,L]\times [\theta_0,\frac{\pi}{2}]\times [0,2\pi]$. Moreover,
$$
\|v_L\|_{0,0;Q_L}^{(0,-\delta_0)} \leqslant C\|f\|_{1,\alpha;Q}^{(1-\delta,2-\delta_0)},\eqno{(4.44)}
$$
where $C>0$ is independent of $L$.

\vskip 0.2 true cm
{\bf Step 3. The convergence of $\nabla_{r,\th,y_2} v_L$ and $\nabla_{r,\th,y_2}^2 v_L$}
\vskip 0.2 true cm

We only give the proof on the convergence of $\p_{y_2}^2v_L$ since the other cases can be treated analogously.
It follows from (4.27) and a direct computation that
\begin{align}
\p_{y_2}^2 &v_L(r,\theta,y_2)
=-\sum_{n=1}^{\infty} n^2 \Big(R_{0n}^{(1)}(r)\sin ny_2+R_{0n}^{(2)}(r)\cos ny_2\Big) \notag\\
&\quad -\sum_{m=1}^{\infty}\sum_{n=1}^{\infty} n^2 \Big(R_{mn}^{(1)}(r)\cos \m(\theta-\theta_0)\sin ny_2 +R_{mn}^{(2)}(r)\cos \m(\theta-\theta_0)\cos ny_2\Big).\tag{4.45}
\end{align}

Thus, we have
$$
|\p_{y_2}^2 v_L(r,\theta,y_2)|
\leqslant \sum_{n=1}^{\infty} n^2 \Big(|R_{0n}^{(1)}(r)|+|R_{0n}^{(2)}(r)|\Big)
+\sum_{m=1}^{\infty}\sum_{n=1}^{\infty} n^2 \Big(|R_{mn}^{(1)}(r)|+|R_{mn}^{(2)})(r)|\Big).\eqno{(4.46)}
$$

By Lemma A.5-Lemma A.6 in Appendix, we have that
\begin{equation*}
\sum\limits_{n=1}^{\infty} n^2 |R_{0n}^{(1)}(r)| \leqslant
\begin{cases}
&C\|f\|_{3,\alpha;Q}^{(1-\delta,2-\delta_0)}r^{\dl-\f52},\quad 0<r\le 1,\\
&C\|f\|_{3,\alpha;Q}^{(1-\delta,2-\delta_0)}, \quad 1<r\le L
\end{cases}
\end{equation*}

and
\begin{equation*}
\sum\limits_{m=1}^{\infty}\sum\limits_{n=1}^{\infty} n^2 |R_{mn}^{(1)}(r)|\le
\begin{cases}
C \|f\|_{4,\alpha;Q}^{(1-\delta,2-\delta_0)}r^{\min\{\delta,\frac 12\}-\frac 52}, \qquad r\leqslant 1,\\
C \|f\|_{4,\alpha;Q}^{(1-\delta,2-\delta_0)}, \qquad \qquad\quad  1<r\leqslant L,
\end{cases}
\end{equation*}
where the generic $C>0$ is independent of $L$.

Similarly, we also have
\begin{equation*}
\sum\limits_{n=1}^{\infty} n^2 |R_{0n}^{(2)}(r)| \leqslant
\begin{cases}
&C\|f\|_{3,\alpha;Q}^{(1-\delta,2-\delta_0)}r^{\dl-\f52},\quad 0<r\le 1,\\
&C\|f\|_{3,\alpha;Q}^{(1-\delta,2-\delta_0)}, \quad 1<r\le L
\end{cases}
\end{equation*}
and
\begin{equation*}
\sum\limits_{m=1}^{\infty}\sum\limits_{n=1}^{\infty} n^2 |R_{mn}^{(2)}(r)|\le
\begin{cases}
&C \|f\|_{4,\alpha;Q}^{(1-\delta,2-\delta_0)}r^{\min\{\delta,\frac 12\}-\frac 52}, \qquad r\leqslant 1,\\
&C \|f\|_{4,\alpha;Q}^{(1-\delta,2-\delta_0)}, \qquad \qquad\quad  1<r\leqslant L.
\end{cases}
\end{equation*}
Thus, the series in (4.46) are convergent for any $(r,\th,y_2)\in (0,L]\times [\theta_0,\pi] \times [0,2\pi]$.

Collecting Step 1-Step 3, we complete the proof of Proposition 4.1.\qed

\vskip 0.4 true cm
\centerline{\bf $\S 5.$ Higher regularities and existence of the solution to (4.2) }
\vskip 0.2 true cm

In this section, based on Proposition 4.1, we will establish the higher regularities of the solution $v_L$ to (4.4) in the domain $Q_{\f{L}{2}}$,
subsequently, we show the solvability of (4.2) in the whole domain $Q$.

\vskip 0.2 true cm

{\bf Lemma 5.1.} {\it Suppose that $v_L$ is a solution to (4.4), then the following estimate holds
$$
||v_L||_{6,\al;Q_{\f{L}{2}}}^{(0,-\dl_0)}\le C||f_L||_{4,\al;Q_{\f{L}{2}}}^{(1-\dl,2-\dl_0)}.\eqno{(5.1)}
$$
}

{\bf Proof.}  We now apply the scaling technique to establish (5.1).

At first, set $y=Lz$, $\tilde{v}(z)=v_L(y)$ and $\tilde{f}(z)=f_L(y)$, then it follows from (4.4) and a
direct computation that $\tilde{v}(z)$ satisfies
$$
\Delta \tilde{v}(z)=L^2\tilde{f}(z)\qquad\q in\q Q_1,
$$
where $\tilde{v}(z)$ and $\tilde{f}(z)$ are $\ds\f{2\pi}{L}-$periodic with respect to $z_2$.

Denoting $Q(r_1,r_2)=\{y: y\in Q_{r_2}, r_1<\sqrt{y_1^2+y_3^2}<r_2\}$, and applying the standard Schauder interior estimate and boundary estimate
(see Chapter 6 of [15]), one has
$$
||\tilde{v}(z)||_{6,\al;Q(\f13,\f23)}\le C\big(||\tilde{v}(z)||_{0;Q(\f14,\f34)}+||L^2\tilde{f}||_{4,\al;Q(\f14,\f34)} \big).\eqno{(5.2)}
$$
Going back to the function $v_L(y)$, we have that for any $y\in Q(\f{L}{3},\f{2L}{3})$
$$
\ss_{m=0}^6L^m|D_y^mv_L|\le C\big(\|v_L\|_{0;Q(\f{L}{4},\f{3L}{4})}+\ss_{m=0}^4L^{m+2}\|D_y^mf_L\|_{0;Q(\f{L}{4},\f{3L}{4})}
+L^{6+\al}[D_y^4f_L]_{0,\al;Q(\f{L}{4},\f{3L}{4})}\big).\eqno{(5.3)}
$$
Noticing that for any $y\in Q(\f{L}{4},\f{3L}{4})$, we have $r_y\sim L$. Then multiplying $L^{-\dl_0}$ on the two hand sides of (5.3)
yields for $L>4$
$$
\ss_{m=0}^6r_y^{m-\dl}|D_y^mv_L(y)|\le C\big(||v_L||_{0;Q(\f{L}{4},\f{3L}{4})}^{(\star,-\dl_0)}
+||f_L||_{4,\al;Q(\f{L}{4},\f{3L}{4})}^{(\star,2-\dl_0)}\big).\eqno{(5.4)}
$$
On the other hand, for $y,z\in Q(\f{L}{3},\f{2L}{3})$, it follows from (5.2) that
\begin{align*}
L^{6+\al}\f{\ds |D^6v_L(y)-D^6v_L(z)|}{\ds |y-z|^{\al}}&\le C\big(\|v_L\|_{0;Q(\f{L}{4},\f{3L}{4})}+\ss_{m=0}^4L^{m+2}
\|D_y^mf_L\|_{0;Q(\f{L}{4},\f{3L}{4})}\\
&\quad +L^{6+\al}[D_y^4f_L]_{0,\al;Q(\f{L}{4},\f{3L}{4})}\big).
\end{align*}

Similar to the proof of (5.4), we have
$$
r_{y,z}^{6+\al-\dl_0}\f{\ds |D^6v_L(y)-D^6v_L(z)|}{\ds |y-z|^{\al}}\le C\big(||v_L||_{0;Q(\f{L}{4},\f{3L}{4})}^{(\star,-\dl_0)}+||f_L||_{4,\al;
Q(\f{L}{4},\f{3L}{4})}^{(\star,2-\dl_0)}\big).\eqno{(5.5)}
$$
Combining  (5.4)-(5.5) with Proposition 4.1 yields
$$
||v_L||_{6,\al;Q(\f{L}{3},\f{2L}{3})}^{(\star,-\dl_0)}\le C||f_L||_{4,\al;Q_L}^{(1-\dl,2-\dl_0)}.\eqno{(5.6)}
$$

We now continue to prove (5.1). For any fixed point $y_0=(y_1^0,y_2^0,y_3^0)\in Q_{\f{L}{2}}$,
we set $r_0=\sqrt{(y_1^0)^2+(y_3^0)^2}$,
$d_0=\mu r_0$
and the cylindrical domain $C_{d_0}(y_0)\equiv B((y_1^0,y_3^0), d_0)\times \Bbb R$, where $0<\mu<1$ is any fixed constant
and $B((y_1^0,y_3^0), d_0)$ stands for a ball centered at $(y_1^0,y_3^0)$ with the radius $d_0$.
We define the map $T: C_{d_0}(y_0)\rightarrow C_1(O)\equiv B((0,0), 1)\times \Bbb R$ by $T(y)=\ds\f{y-y_0}{d_0}$ for $y\in C_{d_0}(y_0)$.

In order to estimate $v_L(y)$ in $Q_{\f{L}{2}}$, we distinguish two cases:

(i) $C_{d_0}(y_0)\subset Q_{\f{L}2}$;

(ii) $C_{d_0}(y_0)\bigcap\p Q_{\f{L}2}\neq\emptyset$.

We now treat these two cases separately. In case (i), we set $\bar{v}(x)=\f{1}{d_0}v_L(y_0+d_0x)$ and $\bar{f}(x)=f(y_0+d_0x)$ for $x\in C_1(O)$. Then it follows that
$$
\Delta\bar{v}(x)=d_0\bar{f}(x),
$$
where $\bar{v}(x)$ and $\bar{f}(x)$ are $\ds\f{2\pi}{d_0}$-periodic with respect to the variable $x_2$.

By the Schauder interior estimate in Chapter 6 of [15], one has
$$
||\bar{v}||_{6,\al;C_{\f23}(O)}\le C\big(||\bar{v}||_{0;C_1(O)}+||d_0\bar{f}||_{4,\al;C_1(O)}\big),\eqno{(5.7)}
$$
where $C>0$ depends only on $\al$.

For $y\in C_{\f{2d_0}{3}}(y_0)$, then $\ds(\f{1}{\mu}-\f23)d_0\le r_y=\sqrt{y_1^2+y_3^2}\le (\f23+\f{1}{\mu})d_0$ holds, and (5.7) means that
\begin{align}
\ss_{m=0}^6d_0^{m-1}|D_y^mv_L(y)|\le& C\big(d_0^{-1}\|v_L\|_{0;C_{d_0}(y_0)}+\ss_{m=0}^4d_0^{m+1}\|D_y^mf_L\|_{0;C_{d_0}(y_0)} \notag\\
&+d_0^{5+\al}[D_y^4f_L]_{0,\al;C_{d_0}(y_0)}\big).\tag{5.8}
\end{align}

If $r_0\ge1$, then multiplying $d_0^{1-\dl_0}$ on the two hand sides of (5.8) to obtain
\begin{align}
\ss_{m=0}^6d_0^{m-\dl_0}|D_y^mv_L(y)|\le &C\big(d_0^{-\dl_0}\|v_L\|_{0;C_{d_0}(y_0)}+\ss_{m=0}^4d_0^{m+2-\dl_0}\|D_y^mf_L\|_{0;C_{d_0}(y_0)}\notag\\
&+d_0^{6+\al-\dl_0}[D_y^4f_L]_{0,\al;C_{d_0}(y_0)}\big).\tag{5.9}
\end{align}

If $r_0<1$, then multiplying $d_0$ on the two hand sides of (5.8) yields
$$
\ss_{m=0}^6d_0^{m}|D_y^mv_L(y)|\le C\big(\|v_L\|_{0;C_{d_0}(y_0)}+\ss_{m=0}^4d_0^{m+2}\|D_y^mf_L\|_{0;C_{d_0}(y_0)}
+d_0^{6+\al}[D_y^4f_L]_{0,\al;C_{d_0}(y_0)}\big).\eqno(5.10)
$$
Noticing $r_y\sim d_0$ for any $y\in C_{\f{2d_0}{3}}(y_0)$, then it follows from (5.9)-(5.10) and Proposition 4.1 that
$$
||v_L||_{6;C_{\f{2d_0}{3}}(y_0)}^{(0,-\dl_0)}\le C_{\mu}\big(||v_L||_{0;Q_L}^{(0,-\dl_0)}+||f_L||_{4,\al;Q_L}^{(-1-\dl,2-\dl_0)}\big)\le C_{\mu}||f_L||_{4,\al;Q_L}^{(-1-\dl,2-\dl_0)}.\eqno(5.11)
$$
On the other hand, if $x_0\in C_{\f{d_0}{2}}(y_0)$, then we can derive $r_{x_0}\sim d_0$ and $r_{x_0,y_0}=\min(r_{x_0},r_{y_0})\sim d_0$, and further obtain by (5.7) that
\begin{align}
&r_{x_0,y_0}^{6+\al-\dl_0}\f{\ds |D^6v(x_0)-D^6v(y_0)|}{\ds |x_0-y_0|^{\al}}\le C_{\mu}\big(||v_L||_{0;Q_L}^{(0,-\dl_0)}+||f_L||_{4,\al;Q_L}^{(-1-\dl,2-\dl_0)}\big),\q r_0\ge1,\tag{5.12}\\
&r_{x_0,y_0}^{6+\al-\dl_0}\f{\ds |D^6v(x_0)-D^6v(y_0)|}{\ds |x_0-y_0|^{\al}}\le C_{\mu}\big(||v_L||_{0;Q_L}^{(0,-\dl_0)}+||f_L||_{4,\al;Q_L}^{(-1-\dl,2-\dl_0)}\big),\q r_0<1.\tag{5.13}
\end{align}
If $x_0\notin C_{\f{d_0}{2}}(y_0)$ but $x_0\in C_{d_0}(y_0)$, then we have
\begin{align}
&r_{x_0,y_0}^{6+\al-\dl_0}\f{\ds |D^6v(x_0)-D^6v(y_0)|}{\ds |x_0-y_0|^{\al}}\le C_{\mu} \ds\sup_{y\in C_{d_0}(y_0)}r_y^{6-\dl_0}|D_y^6v_L(y)|,\q r_0\ge 1,\tag{5.14}\\
&r_{x_0,y_0}^{6+\al}\f{\ds |D^6v(x_0)-D^6v(y_0)|}{\ds |x_0-y_0|^{\al}}\le C_{\mu} \ds\sup_{y\in C_{d_0}(y_0)}r_y^{6}|D_y^6v_L(y)|,\q r_0<1.\tag{5.15}
\end{align}

In case (ii), set $\bar{v}(y)=\f{1}{d_0}v_L(y_0+d_0y)$ and $\bar{f}(y)=f_L(y_0+d_0y)$
for $y\in M\equiv T\big(C_{d_0}(y_0)\bigcap Q_L\big)$. As in case (i), since we already have shown $v_L\in C^{6,\al}(\overline{Q(\f{L}{3},\f{2L}{3})})$, then
it follows from the Schauder boundary estimate in Chapter 6 of [15] that
$$
||\bar{v}||_{6,\al;C_{\f12}(O)\bigcap M}\le C\big(||\bar{v}||_{0;M}+||d_0\bar{f}||_{4,\al;M}\big),\eqno(5.16)
$$
where $C$ depends only on $\al$. Thus similar to the proof in case (i), one can obtain the similar estimates as in (5.11)-(5.15). Combining all estimates (i) and case (ii), we can derive that (5.1) holds. Thus, the proof of Lemma 5.1 is complete.\qed

Next, we focus on improving the regularities of $v_L$ near $r=0$.

\vskip 0.2 true cm

{\bf Lemma 5.2.} {\it Suppose that $v_L$ is the solution to (4.4), then the following estimate holds
$$
||v_L||_{6,\al;Q_{\f{L}{2}}}^{(-1-\dl,-\dl_0)}\le C||f_L||_{4,\al;Q_L}^{(1-\dl,2-\dl_0)}.\eqno(5.17)
$$
}

{\bf Proof.} From Lemma 3.1 of [20], one has
$$
|v_L|_{2+\al;Q_1}^{(-1-\dl)}\le C\big(|f_L|_{\al;Q_1}^{(1-\dl)}+|v_L|_{0;Q_1}\big).\eqno(5.18)
$$
Then it follows from (4.8), (5.18) and Lemma 2.2 that
$$
|v_L|_{2+\al;Q_1}^{(-1-\dl)}\le C||f_L||_{4,\al;Q_L}^{(1-\dl,2-\dl_0)}.\eqno(5.19)
$$
On the other hand, it follows from (i)-(ii) in Lemma 2.2 and (5.19) that
\begin{align}
& |v_L|_{1+\dl;\overline{Q_1}}\le |v_L|_{2+\al;Q_1}^{(-1-\dl)}\le C ||f_L||_{4,\al;Q_L}^{(1-\dl,2-\dl_0)},\tag{5.20}\\
& |v_L|_{2;Q_1}^{(-1-\dl)}\le C|v_L|_{2+\al;Q_1}^{(-1-\dl)}\le C ||f_L||_{4,\al;Q_L}^{(1-\dl,2-\dl_0)}.\tag{5.21}
\end{align}
It follows from the boundary condition in (4.4) and the regularity of $v_L$ in (5.20) that
$$
\p_{y_1}v_L(0,y_2,0)=\p_{y_3}v_L(0,y_2,0)=0. \eqno(5.22)
$$
From (5.20), we have for any $y\in Q_1$
$$
\ds\sup_{y\in Q_1}\f{\ds |\p_{y_1}v_L(y)-\p_{y_1}v_L(0,y_2,0)|}{\ds |y_1^2+y_3^2|^{\f{\dl}2}}\le C ||f_L||_{4,\al;Q_L}^{(1-\dl,2-\dl_0)},
$$
which means
$$
\ds\sup_{y\in Q_1}|r_y^{-\dl}\p_{y_1}v_L|\le C ||f_L||_{4,\al;Q_L}^{(1-\dl,2-\dl_0)}.\eqno(5.23)
$$
Similarly,
$$
\ds\sup_{y\in Q_1}|r_y^{-\dl}\p_{y_3}v_L|\le C ||f_L||_{4,\al;Q_L}^{(1-\dl,2-\dl_0)}.\eqno(5.24)
$$

In addition, it follows from (5.21) that
$$
\ds\sup_{\si>0}\big(\si^{1-\dl}\ds\sup_{y\in Q_{1\si}}|D^2v_L|\big)\le C ||f_L||_{4,\al;Q_L}^{(1-\dl,2-\dl_0)}.\eqno(5.25)
$$
Choosing $\si=\ds\f{r_y}{2}$, then by (5.25) we have for any $y\in Q_1$
$$
\big(\f{r_y}{2}\big)^{1-\dl}\ds\sup_{z\in Q_{1\si}}|D^2v_L(z)|\le C ||f_L||_{4,\al;Q_L}^{(1-\dl,2-\dl_0)}.
$$
This, together with $r_y>\si$, yields
$$
\ds\sup_{y\in Q_1}r_y^{1-\dl}|D^2v_L(y)|\le C ||f_L||_{4,\al;Q_L}^{(1-\dl,2-\dl_0)}.\eqno(5.26)
$$

Based on (5.26), we then derive the higher estimates near the $y_2-$axis. Let $w_2(y)=\p_{y_2}^2v_L(y)$,
then $w_2(y)$ satisfies following equation
\begin{equation}
\left\{
\begin{aligned}
&\Delta w_2=\p_{y_2}^2f_L\qquad\,\, \q \text{in}\q Q_L,\\
&\p_{n}w_2=0\q\qquad\qquad \, \text{on}\q\th=\th_0 \quad~\text{and}\quad ~\th=\pi,\\
&w_2(y)=\p_{y_2}^2v_L(y)\quad\q \text{on}\q\sqrt{y_1^2+y_2^2}=1.
\end{aligned}
\right.\tag{5.27}
\end{equation}

It follows from (5.27) and the scaling method used in Lemma 5.1 that
$$
||w_2||_{4,\al;Q_1}^{(1-\dl,\star)}\le C\big(||w_2||_{0;Q_1}^{(1-\dl,\star)}+||\p_{2}^2f_L||_{2,\al;Q_1}^{(3-\dl,\star)}\big).
$$
This, together with (5.26), yields
$$
||\p_{y_2}^2v_L||_{4,\al;Q_1}^{(1-\dl,\star)}\le C||f_L||_{4,\al;Q_1}^{(1-\dl,\star)}.\eqno(5.28)
$$

Next, we focus on the estimates on  $\p_{y_1}v_L$. Let $w_1=\p_{y_1}v_L$, we now derive the boundary conditions of $w_1$.On $y_1=0$, $w_1=\p_nv_L=0$ holds. On $y_3=b_0y_1$,  one has from (4.4) that
$$
\p_nv_L=\sin\th_0\p_{y_1}v_L-\cos\th_0\p_{y_3}v_L=0.\eqno(5.29)
$$
Taking the tangent derivative $\cos\th_0\p_{y_1}+\sin\th_0\p_{y_3}$ on two hand sides of (5.29) yields
\begin{align}
0&=(\cos\th_0\p_{y_1}+\sin\th_0\p_{y_3})(\sin\th_0\p_{y_1}v_L-\cos\th_0\p_{y_3}v_L)\notag\\
&=\sin\th_0\cos\th_0\p_{y_1}^2v_L+(\sin^2\th_0-\cos^2\th_0)\p_{y_1y_3}^2v_L
-\sin\th_0\cos\th_0\p_{y_3}^2v_L.\tag{5.30}
\end{align}
In addition, it follows from (4.4) that
$$
\p_{y_3}^2v_L=f_L-\p_{y_1}^2v_L-\p_{y_2}^2v_L.\eqno(5.31)
$$
Substituting (5.31) into (5.30) yields
$$
2b_0\p_1w_1+(b_0^2-1)\p_3w_1=b_0f_L-b_0\p_{y_2}^2v_L.\eqno(5.32)
$$
It is easy to verify
$$
(2b_0,b_0^2-1)\cdot \vec{n}=(1+b_0^2)\cos\th_0>0.
$$
Thus, $w_1(y)$ satisfies the following problem with the oblique derivative boundary condition on $\th=\th_0$
\begin{equation}
\left\{
\begin{aligned}
&\Delta w_1=\p_{y_1}f_L  \qquad\qquad \qquad\qquad \qquad\qquad \quad\q \text{in}\q Q_L,\\
&2b_0\p_{y_1}w_1+(b_0^2-1)\p_{y_3}w_1=b_0f_L-b_0\p_{y_2}^2v_L\qquad \q \text{on}\q\th=\th_0 ,\\
&w_1=0\qquad\qquad\qquad\qquad  \qquad \qquad\q\text{on}\q\th=\f{\pi}2,\\
&w_1(y)=\p_{y_1}v_L(y)\qquad \qquad\qquad \quad\q \text{on}\q\sqrt{y_1^2+y_3^2}=1.
\end{aligned}
\right.\tag{5.33}
\end{equation}

Still applying the scaling method as in Lemma 5.1, we have
$$
||w_1||_{5,\al;Q_1}^{(-\dl,\star)}\le C\bigg(||w_1||_{0;Q_1}^{(-\dl,\star)}+||b_0f_L-b_0\p_{y_2}^2v_L||_{4,\al;Q_1}^{(1-\dl,\star)} +||\p_{y_1}f_L||_{3,\al;Q_1}^{(2-\dl,\star)}\bigg).
$$
This, together with (5.23) and (5.28), yields
$$
||\p_{y_1}v_L||_{5,\al;Q_1}^{(-\dl,\star)}\le C||f_L||_{4,\al;Q_L}^{(-1-\dl,2-\dl_0)}.\eqno(5.34)
$$
Similarly, for $\p_{y_2}v_L$, we have same estimates as follows
$$
||\p_{y_2}v_L||_{5,\al;Q_1}^{(-\dl,\star)}\le C||f_L||_{4,\al;Q_L}^{(-1-\dl,2-\dl_0)}.\eqno(5.35)
$$
Thus, combining (5.1), (5.20), (5.28), and (5.34)-(5.35), we obtain (5.17) and complete the proof of Lemma 5.2.
\qed

Finally, we start to prove the existence of the solution to (4.3), which also means the existence of the problem (4.2).
At first, it follows from (4.5) and (5.17) that
$$
||v_L||_{6,\al;Q_{\f{L}{2}}}^{(-1-\dl,-\dl_0)}\le C||f||_{4,\al;Q}^{(1-\dl,2-\dl_0)}.\eqno(5.36)
$$
Suppose that $\tilde{v}_L$ is an extension of $v_L$ in the whole domain $Q$
(one can see Theorem 7.25 of [15]), which satisfies
$$
\tilde{v}_L|_{Q_{\f{L}{2}}}=v_L,\q ||\tilde{v}_L||_{6,\al;Q}^{(-1-\dl,-\dl_0)}\le C||v_L||_{6,\al;Q_{\f{L}{2}}}^{(-1-\dl,-\dl_0)}.\eqno(5.37)
$$
Then it follows from (5.36)-(5.37) that
$$
||\tilde{v}_L||_{6,\al;Q}^{(-1-\dl,-\dl_0)}\le C||f||_{4,\al;Q}^{(1-\dl,2-\dl_0)}.\eqno(5.38)
$$

Let $L\rightarrow+\infty$, by the standard diagonal method, we can extract a convergent sub-sequence $\tilde{v}_{L_n}$
($n\in\Bbb N$) and a function $v\in H_{6,\al;Q}^{(-1-\dl,-\dl_0)}$ such that
$$
||\tilde{v}_{L_n}-v||_{6,\al;Q_N}^{(-1-\dl,-\dl_0)}\rightarrow 0 \qquad as \qquad n\to\infty,
$$
where $Q_N$ is any fixed sub-domain of $Q$.
Moreover, it is easy to know that $v$ is a solution of (4.3).

\vskip 0.4 true cm
\centerline{\bf $\S 6.$  The uniqueness of the solution to (4.3) }
\vskip 0.4 true cm

In this section, we focus on the uniqueness of the solution $v$ to (4.3) since  the existence of solution $v$ to
(4.3) has been shown in $\S 5$. To this end, we will use the separation  method as in $\S 4$ together with
some technical analysis so that the difficulty
induced by the lack of the maximum principle for (4.3) can be overcome.

\vskip 0.2 true cm

{\bf Lemma 6.1.} {\it There exists a unique solution to (4.3) such that $v\in H_{6,\al}^{(-1-\dl,-\dl_0)}(Q)$.}

{\bf Proof.} Suppose that $v_1(y), v_2(y)\in H_{6,\alpha}^{(-1-\delta,-\delta_0)}(Q)$ are different solutions to (4.3). Then for any $\delta_1>\delta_0$, we have $v_1(y), v_2(y)\in H_{6,\alpha}^{(-1-\delta,-\delta_1)}(Q)$. Let
$W(y)=v_1(y)-v_2(y)$, then $W(y)\in H_{6,\alpha}^{(-1-\delta,-\delta_0)}(Q)\cap H_{6,\alpha}^{(-1-\delta,-\delta_1)}(Q)$.
Denote $W_L(y)$ by the restriction of $W(y)$ on the domain $Q_L$, then we have
\begin{equation}
\left\{
\begin{aligned}
&\triangle W_L=0 \qquad\qquad\qquad\qquad\qquad \mbox{in} \quad  Q_L,\\
&\p_n W_L=0 \qquad\qquad\qquad\qquad \mbox{on} \qquad  \theta=\theta_0\quad  \mbox{and}\quad  \theta=\pi,\\
&W_L(y)=W_L(y_1,y_2+2\pi,y_3),\\
&W_L(0,0,0)=0,\\
&\ds\frac{\p W_L}{\p n}(y)=\frac{\p v_1}{\p n}(y)-\frac{\p v_2}{\p n}(y)=g(y)\qquad \mbox{on}\quad \sqrt{y_1^2+y_3^2}=L.
\end{aligned}
\right.\tag{6.1}
\end{equation}

Introducing the cylindrical $(r,\theta,y_2)$ as in (4.2), then (6.1) can be rewritten as
\begin{equation}
\left\{
\begin{aligned}
&\p_r^2 W_L + r^{-2}\p_{\theta}^2 W_L + \p_{y_2}^2 W_L + r^{-1}\p_{r} W_L =0 \qquad \mbox{in} \quad  Q_L,\\
&\p_n W_L=0 \qquad\qquad\qquad\qquad\qquad\qquad\qquad \mbox{on} \qquad  \theta=\theta_0\quad  \mbox{and}\quad \theta=\frac{\pi}{2},\\
&W_L(r,\theta,y_2)=W_L(r,\theta,y_2+2\pi),\\
&W_L(0,\theta,0)=0,\\
&\p_r W_L (L,\theta,y_2)=g(L,\theta,y_2)\quad \mbox{on}~ \Si_L\equiv\{y: \sqrt{y_1^2+y_3^2}=L, y_2\in \Bbb R\},
\end{aligned}
\right.\tag{6.2}
\end{equation}
where $g(L,\theta,y_2)$ satisfies the compatibility conditions
$$
\text{$\p_\theta g(L,\theta,y_2)=\p_\theta^3 g(L,\theta,y_2)=\p_\theta^5 g(L,\theta,y_2)=0$  \quad on\quad  $\theta=\theta_0$\quad  and\quad $\theta=\frac{\pi}{2}$,}\eqno(6.3)
$$
and $\ds\sup_{0\le\al_1+\al_2\le5}|L^{1-\delta_1+\alpha_1+\alpha_2} (\frac 1r \p_\theta)^{\alpha_1} \p_{y_2}^{\alpha_2} g|\leqslant C$,
which come from the regularity of $W(y)\in H_{6,\alpha}^{(-1-\delta,-\delta_1)}$ and
$$
\sum_{0\leqslant \alpha_1+\alpha_2 \leqslant 5} \sup_{\theta\in[\theta_0,\frac{\pi}{2}]\atop y_2\in \rr}
L^{1-\delta_1+\alpha_2} |\p_\theta^{\alpha_1} \p_{y_2}^{\alpha_2} g|=\|g\|_{5,\Si_L}^{(1-\delta_1)} \leqslant CL^{\delta_0-\delta_1} \|g\|_{5,\Si_L}^{(1-\delta_0)},\eqno(6.4)
$$
respectively. Moreover, the solvability condition holds
$$
\int_0^{2\pi} \int_{\theta_0}^{\frac {\pi}2} g(L,\theta,y_2) d\theta dy_2=0.\eqno(6.5)
$$

By $\S 4$, it is known that  the eigenvalues of the corresponding homogeneous problem of (6.2)
are $\lambda_n=n^2 (n=0,1,\cdots)$ and $\mu_m=(\frac{m\pi}{\frac{\pi}2-\theta_0})^2 (m=0,1,\cdots)$,
and the related complete
orthogonal basis of eigenfunction functions is $\{\cos \m(\theta-\theta_0)\sin(ny_2), \cos \m(\theta-\theta_0)\cos(ny_2)
\}_{m,n=0}^{\infty}$. Suppose that the expansion of $g(L,\theta,y_2)$ is
\begin{align}
g(L,\theta,y_2)&=g_{00}(L)+\sum_{m=1}^{\infty} g_{m0}(L) \cos \m(\theta-\theta_0)+\sum_{n=1}^{\infty}\Big(g_{0n}^{(1)}(L)\sin ny_2+g_{0n}^{(2)}(L)\cos ny_2\Big)\notag\\
&+\sum_{m=1}^{\infty}\sum_{n=1}^{\infty}\Big(g_{mn}^{(1)}(L)\cos \m(\theta-\theta_0)\sin ny_2+g_{mn}^{(2)}(L)\cos \m(\theta-\theta_0)\cos ny_2\Big),\tag{6.6}
\end{align}

where
\begin{align*}
& g_{00}(L)= \frac{1}{2\pi (\frac{\pi}2-\theta_0)} \int_{\theta_0}^{\frac{\pi}2} \int_{0}^{2\pi} g(L,\theta,y_2) dy_2 d \theta, \\
& g_{m0}(L) = \frac{1}{\pi (\frac{\pi}2-\theta_0)} \int_{\theta_0}^{\frac{\pi}2} \int_{0}^{2\pi} g(L,\theta,y_2) \cos \m(\theta-\theta_0) dy_2 d \theta, \\
& g_{0n}^{(1)}(L) = \frac{1}{\pi (\frac{\pi}2-\theta_0)} \int_{\theta_0}^{\frac{\pi}2}
\int_{0}^{2\pi} g(L,\theta,y_2) \sin ny_2dy_2 d \theta, \\
&g_{0n}^{(2)}(L) = \frac{1}{\pi (\frac{\pi}2-\theta_0)} \int_{\theta_0}^{\frac{\pi}2}
\int_{0}^{2\pi} g(L,\theta,y_2) \cos ny_2 dy_2 d \theta, \\
& g_{mn}^{(1)}(L) = \frac{2}{\pi (\frac{\pi}2-\theta_0)} \int_{\theta_0}^{\frac{\pi}2} \int_{0}^{2\pi} g(L,\theta,y_2) \cos \m(\theta-\theta_0) \sin ny_2 dy_2 d \theta, \\
& g_{mn}^{(2)}(L) = \frac{2}{\pi (\frac{\pi}2-\theta_0)} \int_{\theta_0}^{\frac{\pi}2} \int_{0}^{2\pi} g(L,\theta,y_2) \cos \m(\theta-\theta_0) \cos ny_2 dy_2 d \theta.
\end{align*}

Let

\begin{align}
W_L (r,\theta,y_2)&= \R_{00}(r)+\sum_{m=1}^{\infty} \R_{m0}(r) \cos \m(\theta-\theta_0) + \sum_{n=1}^{\infty} \Big(\R_{0n}^{(1)}(r) \sin ny_2 + \R_{0n}^{(2)}(r)\cos n y_2 \Big) \notag\\
 & + \sum_{m=1}^{\infty} \sum_{n=1}^{\infty} \Big(\R_{mn}^{(1)}(r) \cos \m(\theta-\theta_0) \sin ny_2 +\R_{mn}^{(2)}(r)\cos \m(\theta-\theta_0) \cos ny_2 \Big).\tag{6.7}
\end{align}

It follows from (6.2)-(6.3) and (6.7) that
\begin{equation}
\left\{
\begin{aligned}
&\R_{00}^{''}(r)+r^{-1} \R_{00}^{'}(r)=0,\\
&\R_{00}^{'}(L)=g_{00}(L)=0,\\
&\R_{00}(0)=-\ds\sum_{n=1}^{\infty} \R_{0n}^{(2)}(0),\\
\end{aligned}
\right.\tag{6.8}
\end{equation}

\begin{equation}
\left\{
\begin{aligned}
&\R_{m0}^{''}(r)-r^{-2} \mu_m \R_{m0}(r)+r^{-1} \R_{m0}^{'}(r)=0,\quad m=1,2,\cdots,\\
&\R_{m0}^{'}(L)=g_{m0}(L),\\
&\text{$\R_{m0}(0)$ \, is \, bounded, and $\R_{m0}(0)+\sum\limits_{n=1}^{\infty}R_{mn}^{(2)}(0)=0$},\\
\end{aligned}
\right.\tag{6.9}
\end{equation}

\begin{equation}
\left\{
\begin{aligned}
&(\R_{0n}^{(i)})^{''}(r)-n^{2} \R_{0n}^{(i)}(r) + r^{-1}(\R_{0n}^{(i)})^{'}(r)=0,\quad n=1,2,\cdots,\\
&(\R_{0n}^{(i)})^{'}(L)=g_{0n}^{(i)}(L),\\
&\R_{0n}^{(i)}(0)\quad \mbox{ is \quad bounded},
\end{aligned}
\right.\tag{6.10}
\end{equation}
and
\begin{equation}
\left\{
\begin{aligned}
&(\R_{mn}^{(i)})^{''}(r)-(r^{-2}\mu_m+n^{2}) \R_{mn}^{(i)}(r)+r^{-1}(\R_{mn}^{(i)})^{'}(r)=0,\quad m,n=1,2,\cdots,\\
&(\R_{mn}^{(i)})^{'}(L)=g_{mn}^{(i)}(L),\\
&\R_{mn}^{(i)}(0)\quad \mbox{ is \quad bounded},
\end{aligned}
\right.\tag{6.11}
\end{equation}

The general solution of (6.8) is
$$
\R_{00}(r)=C_{00}^1+C_{00}^2 \ln r.
$$
By the boundary condition $\R'_{00}(L)=0$, we have $C_{00}^2=0$ and $C_{00}^1=-\ds\sum_{n=1}^{\infty} \R_{0n}^{(2)}(0)$. Thus
$$
\R_{00}(r)=-\sum_{n=1}^{\infty} \R_{0n}^{(2)}(0).\eqno(6.12)
$$
For the equation (6.9), its general solution is
$$
\R_{m0}(r)=C_{m0}^1 r^{\m} + C_{m0}^2 r^{-\m}.
$$
It is noted that the boundedness of $\R_{m0}(0)$ implies $C_{m0}^2=0$ and
$\R'_{m0}(L)=g_{m0}(L)$ derives $C_{m0}^1 =\ds\frac{g_{m0}(L)}{\m} L^{1-\m}$, then
$$
\R_{m0}(r)=\frac{g_{m0}(L)}{\m} L^{1-\m} r^{\m}.\eqno(6.13)
$$
In addition, the general solution of (6.11) is
$$
\R_{mn}^{(i)}(r)=C_{mn}^{i1} I_{\m}(nr) + C_{mn}^{i2} K_{\m}(nr),\qquad i=1,2.
$$
Due to the boundedness of $\R_{mn}^{(i)}(0)$ and $\Bigl(\R_{mn}^{(i)}\Bigr)'(L)=g_{mn}^{(i)}(L)$, we get
$$
\text{$C_{mn}^{i2}=0$,\qquad $C_{mn}^{i1} = \frac{g_{mn}^{(i)}(L)}{nI'_{\m}(nL)}$},
$$
and then
$$
\text{$\R_{mn}^{(i)}(r)=\frac{g_{mn}^{(i)}(L)}{nI'_{\m}(nL)} I_{\m}(nr)$ \quad for \, $m=0,1,\cdots, n=1,2,\cdots$}.\eqno(6.14)
$$
Here we have used $\R_{0n}^{(i)}(0)=\ds\frac{g_{0n}^{(i)}(L)}{nI'_{0}(nL)}$ and $\R_{mn}^{(i)}(0)=0$.

Combining (6.12)-(6.14), we get a formal solution of (6.2) as follows
\begin{align}
W_L (r,\theta,y_2)&= \R_{00}(r)+\sum_{m=1}^{\infty} \R_{m0}(r) \cos \m(\theta-\theta_0) + \sum_{n=1}^{\infty} \Big(\R_{0n}^{(1)}(r) \sin ny_2 + \R_{0n}^{(2)}(r)\cos n y_2 \Big)\notag\\
 & + \sum_{m=1}^{\infty} \sum_{n=1}^{\infty} \Big(\R_{mn}^{(1)}(r) \cos \m(\theta-\theta_0) \sin ny_2 +\R_{mn}^{(2)}(r)\cos \m(\theta-\theta_0) \cos ny_2 \Big),\tag{6.15}
\end{align}
where
\begin{equation}
\left\{
\begin{aligned}
&\R_{00}(r)=-\ds\sum_{n=1}^{\infty} \R_{0n}^{(2)}(0),\\
&\R_{m0}(r)=\frac{g_{m0}(L)}{\m} L^{1-\m} r^{\m},\quad m=1,2,\cdots,\\
&\R_{mn}^{(i)}(r)=\frac{g_{mn}^{(i)}(L)}{nI'_{\m}(nL)} I_{\m}(nr),\quad i=1,2,\quad m=0,1,\cdots, n=1,2,\cdots.
\end{aligned}
\right.\tag{6.16}
\end{equation}

Next, we show that the series in (6.15) is uniformly convergent in $ Q_L$.
To prove the convergence, we require to establish some estimates on the coefficients in (6.6) as in $\S 4$.
It follows from a direct computation that
\begin{align}
g_{m0}(L) & = \frac{1}{\pi (\frac{\pi}2-\theta_0)} \int_{\theta_0}^{\frac{\pi}2} \int_{0}^{2\pi} g(L,\theta,y_2) \cos \m(\theta-\theta_0) dy_2 d \theta \notag\\
&=-\frac{1}{\m \pi (\frac{\pi}2-\theta_0)} \int_{0}^{2\pi} \int_{\theta_0}^{\frac{\pi}2} \p_{\theta} g(L,\theta,y_2) \sin \m(\theta-\theta_0) dy_2 d \theta,\tag{6.17}\\
g_{0n}^{(1)}(L) & = \frac{1}{\pi (\frac{\pi}2-\theta_0)} \int_{\theta_0}^{\frac{\pi}2}
\int_{0}^{2\pi} g(L,\theta,y_2) \sin ny_2dy_2 d \theta \notag\\
&=\frac{1}{n\pi (\frac{\pi}2-\theta_0)} \int_{\theta_0}^{\frac{\pi}2}
\int_{0}^{2\pi} \p_{y_2}g(L,\theta,y_2) \cos ny_2dy_2 d \theta,\tag{6.18}
\end{align}
and
\begin{align}
g_{mn}^{(1)}(L) &= \frac{2}{\pi (\frac{\pi}2-\theta_0)} \int_{\theta_0}^{\frac{\pi}2} \int_{0}^{2\pi} g(L,\theta,y_2) \cos \m(\theta-\theta_0) \sin ny_2 dy_2 d \theta\notag\\
&=-\frac{2}{\m \pi (\frac{\pi}2-\theta_0)} \int_{\theta_0}^{\frac{\pi}2} \int_{0}^{2\pi} \p_{\theta}g(L,\theta,y_2) \sin \m(\theta-\theta_0) \sin ny_2 dy_2 d \theta\notag\\
&=\frac{2}{\mu_m \pi (\frac{\pi}2-\theta_0)} \int_{\theta_0}^{\frac{\pi}2} \int_{0}^{2\pi} \p_{\theta}^2g(L,\theta,y_2) \cos \m(\theta-\theta_0) \sin ny_2 dy_2 d \theta\notag\\
&=\frac{2}{n\mu_m \pi (\frac{\pi}2-\theta_0)} \int_{\theta_0}^{\frac{\pi}2} \int_{0}^{2\pi} \p_{y_2}\p^2_{\theta}g(L,\theta,y_2) \cos \m(\theta-\theta_0) \cos ny_2 dy_2 d \theta.\tag{6.19}
\end{align}

From (6.17)-(6.19), we arrive at
\begin{equation}
\left\{
\begin{aligned}
&|g_{m0}(L)|\leqslant C \mu_m^{-\frac 12} \|g\|_{1,0;\Si_L}^{(1-\delta_1)} L^{\delta_1-1},\\
&|g_{0n}^{(1)}(L)|\leqslant C n^{-1} \|g\|_{1,0;\Si_L}^{(1-\delta_1)} L^{\delta_1-2},\\
&|g_{mn}^{(1)}(L)|\leqslant C n^{-1} \mu_m^{-1} \|g\|_{3,0;\Si_L}^{(1-\delta_1)} L^{\delta_1-2}.
\end{aligned}
\right.\tag{6.20}
\end{equation}

Analogously, (6.20) are also true for $|g_{0n}^{(2)}(L)|$ and $|g_{mn}^{(2)}(L)|$.

We now show that the series in the expression of $W_L (r,\theta,y_2)$ are convergent for any fixed point $(r,\theta,y_2)\in (0,L]\times [\theta_0,\frac{\pi}{2}] \times [0,2\pi]$. At first, by (6.16) and (6.20), we get
\begin{align}
\sum_{m=1}^{\infty}|\R_{m0}(r)| &\leqslant \sum_{m=1}^{\infty}\frac{|g_{m0}(L)|}{\m} L^{1-\m} r^{\m}\notag\\
&\leqslant C\|g\|_{1,0;\Si_L}^{(1-\delta_1)} \sum_{m=1}^{\infty} \mu_m^{-1}L^{\delta_1-\m} r^{\m}\notag\\
&\leqslant C\|g\|_{1,0;\Si_L}^{(1-\delta_1)} r^{\delta_1} \sum_{m=1}^{\infty} \mu_m^{-1}\notag\\
&\leqslant C\|g\|_{1,0;\Si_L}^{(1-\delta_1)} r^{\delta_1},\tag{6.21}
\end{align}
which derives
\begin{equation}
\left\{
\begin{aligned}
&\sum\limits_{m=1}^{\infty}|\R_{m0}(r)|\leqslant C\|g\|_{1,0;\Si_L}^{(1-\delta_1)},\qquad r\leqslant 1;\\
&r^{-\delta_1}\sum\limits_{m=1}^{\infty}|\R_{m0}(r)|\leqslant C\|g\|_{1,0;\Si_L}^{(1-\delta_1)},\qquad r\geqslant 1.
\end{aligned}
\right.\tag{6.22}
\end{equation}

Next, we use the properties of modified Bessel functions in Lemma 2.3 and (6.20) to show the convergence of  the
series $\ds\sum_{m=1}^{\infty}\ds\sum_{n=1}^{\infty}|\R_{mn}^{(1)}(r)|$. To this end, the following three cases
will be considered separately, where $M$ will represent a suitably large fixed integer.
\vskip 0.2 true cm

{\bf Case (a)\quad $m<M$ and $nr<1$}
\begin{align}
\sum_{m<M}\sum_{nr<1}|\R_{mn}^{(1)}(r)| &\leqslant \sum_{m<M}\sum_{nr<1} \frac{|g_{mn}^{(1)}(L)|}{nI'_{\m}(nL)} I_{\m}(nr)\notag\\
&\leqslant \sum_{m<M}\sum_{nr<1} n^{-1} \frac{\sqrt{2\pi nL}}{e^{nL}}\frac{e^{nr}(nr)^{\m}}{2^{\m}\Gamma(\m+1)} |g_{mn}^{(1)}(L)|\notag\\
&\leqslant C \|g\|_{2,0;\Si_L}^{(1-\delta_1)}\sum_{m<M}\sum_{nr<1} n^{-\frac 32} e^{-n}2^{-\m}\mu_m^{-\frac 12} L^{\delta_1-\frac 32}\notag\\
&\leqslant C \|g\|_{2,0;\Si_L}^{(1-\delta_1)}.\tag{6.23}
\end{align}

{\bf Case (b)\quad $m<M$ and $1\leqslant nr \leqslant nL$}
\begin{align}
\sum_{m<M}\sum_{1\leqslant nr \leqslant nL}|\R_{mn}^{(1)}(r)| &\leqslant \sum_{m=1}^{\infty}\sum_{n=1}^{\infty} \frac{|g_{mn}^{(1)}(L)|}{nI'_{\m}(nL)} I_{\m}(nr)\notag\\
&\leqslant C \|g\|_{3,0;\Si_L}^{(1-\delta_1)}\sum_{m<M}\sum_{1\leqslant nr \leqslant nL} n^{-1} \frac{\sqrt{2\pi nL}}{e^{nL}} \frac{e^{nr}}{\sqrt {2\pi nr}} n^{-1} \mu_m^{-1} L^{\delta_1-2}\notag\\
&\leqslant C \|g\|_{3,0;\Si_L}^{(1-\delta_1)} L^{\delta_1-\frac 32} \sum_{m<M}\sum_{1\leqslant nr \leqslant nL} n^{-\frac 32} \mu_m^{-1}\notag\\
&\leqslant C \|g\|_{3,0;\Si_L}^{(1-\delta_1)}.\tag{6.24}
\end{align}

{\bf Case (c)\quad  $m>M$ }

It follows from a direct computation and Lemma 2.3 that
\begin{align}
\sum_{m=1}^{\infty}\sum_{n=1}^{\infty}|\R_{mn}^{(1)}(r)|
&\leqslant \sum_{m=1}^{\infty}\sum_{n=1}^{\infty} \frac{|g_{mn}^{(1)}(L)|}{nI'_{\m}(nL)} I_{\m}(nr)\notag\\
&\leqslant \sum_{m=1}^{\infty}\sum_{n=1}^{\infty} n^{-1}|g_{mn}^{(1)}(L)|
 \frac{z(L)}{(1+z^2(L))^{\frac 14}} \frac{1}{(1+z^2(r))^{\frac 14}} e^{\m \bigl(\eta(r)-\eta(L)\bigr)} \notag\\
&\leqslant C\|g\|_{3,0;\Si_L}^{(1-\delta_1)} L^{\delta_1-\frac 32} \sum_{m=1}^{\infty}\sum_{n=1}^{\infty} n^{-\frac 32} \mu_m^{-\frac 34}\notag\\
&\leqslant C\|g\|_{3,0;\Si_L}^{(1-\delta_1)}.\tag{6.25}
\end{align}
Thus, collecting (6.23)-(6.25) yields
$$
\sum_{m=1}^{\infty}\sum_{n=1}^{\infty}|\R_{mn}^{(1)}(r)| \leqslant C\|g\|_{3,0;\Si_L}^{(1-\delta_1)}.
\eqno(6.26)
$$
Finally, we prove the convergence of the series $\ds\sum_{n=1}^{\infty}|\R_{0n}^{(1)}(r)|$. In fact, we have
\begin{align}
\sum_{n=1}^{\infty}|\R_{0n}^{(1)}(r)|  &\leqslant \sum_{n=1}^{\infty} \frac{|g_{0n}^{(1)}(L)|}{nI'_{0}(nL)} I_{0}(nr)\notag\\
&\leqslant C \|g\|_{1,0;\Si_L}^{(1-\delta_1)} L^{\delta_1-2} \sum_{n=1}^{\infty} n^{-2} \frac{\sqrt{2\pi nL}}{e^{nL}} e^{nr}\notag\\
&\leqslant C \|g\|_{1,0;\Si_L}^{(1-\delta_1)} L^{\delta_1-\frac 32} \sum_{n=1}^{\infty} n^{-\frac 32}\notag\\
&\leqslant C \|g\|_{1,0;\Si_L}^{(1-\delta_1)},\tag{6.27}
\end{align}
which implies
$$
|\R_{00}(r)|\leqslant \sum_{n=1}^{\infty} |\R_{0n}^{(2)}(0)|\leqslant C \|g\|_{1,0;\Si_L}^{(1-\delta_1)}.\eqno(6.28)
$$
By (6.21) and (6.26)-(6.28), we have established the convergence of $W_L (r,\theta,y_2)$ in $ Q_L$ and the following estimates
\begin{equation*}
\left\{
\begin{aligned}
&|W_L (r,\theta,y_2)|
\leqslant C \|g\|_{3,0;\Si_L}^{(1-\delta_1)},\qquad  r \leqslant 1;\\
&r^{-\delta_1}|W_L (r,\theta,y_2)|
\leqslant C \|g\|_{3,0;\Si_L}^{(1-\delta_1)}, \qquad r \geqslant 1.
\end{aligned}
\right.
\end{equation*}

And hence
$$
\|W_L\|_{0,0;Q_L}^{(0,-\dl_1)}\le C\|g\|_{3,0;\Si_L}^{(1-\delta_1)},\eqno(6.29)
$$
where the constant $C>0$ is independent of $L$.

Let $L\rightarrow+\infty$, then one has from (6.4) that
$$
\|g\|_{3,0;\Si_L}^{(1-\delta_1)}\rightarrow 0.\eqno(6.30)
$$
Combining (6.29) with (6.30) yields
$$
\|W\|_{0,0;Q}^{(0,-\dl_1)}=\ds\lim_{L\rightarrow+\infty}\|W_L\|_{0,0;Q_L}^{(0,-\dl_1)}=0.
$$
Thus, the proof of  Lemma 6.1 is complete. \qed

Going back to (4.2), based on Lemma 5.1 and Lemma 6.1 we have
\vskip 0.2 true cm

{\bf Proposition 6.2.} {\it Suppose $\dot f\in H_{4,\al}^{(1-\dl,2-\dl_0)}(Q)$, $\dot g_i\in H_{5,\al}^{(-\dl,1-\dl_0)} \quad
(i=1,2)$, then there exists a unique solution $\dot u$ to (4.2), which satisfies the following estimate
$$
\|\dot{u}\|_{6,\al;Q}^{(-1-\dl,-\dl_0)}\le C\bigg(\|\dot f\|_{4,\al;Q}^{(1-\dl,2-\dl_0)}
+\ss_{i=1}^2\|\dot{g}_i\|_{5,\al;Q}^{(-\dl,1-\dl_0)}\bigg).\eqno(6.31)$$}

\vskip 0.4 true cm
\centerline{\bf $\S 7.$ Proofs of Theorem 3.1 and Theorem 1.1. }
\vskip 0.4 true cm

In this section, first we will use the contraction mapping principle to show Theorem 3.1. To this end, we define
the space $K=\{v\in C(\bar Q): v-u_0\in H_{6,\al}^{(-1-\dl,-\dl_0)}, \|v-u_0\|_{6,\al; Q}^{(-1-\dl,-\dl_0)}\le \ve\}$. Set $u=\dot{u}+u_0$, where $\dot{u}$ is defined as the solution to the following linearized problem which is analogous to the problem (4.1) in $\S 4$
\begin{equation}
\left\{
\begin{aligned}
&\Delta\dot{u}=\dot{F}\qquad \q\text{in}\q Q,\\
&\p_{n}\dot{u}=\dot{G}_1\quad\q\text{on}\q y_3=b_0y_1,\\
&\p_{y_1}\dot{u}=\dot{G}_2\quad\q\text{on}\q y_1=0,\\
&\dot{u}(0,0,0)=0,\\
&\ds\lim_{y_1+y_3\rightarrow\infty}|\nabla_y\dot{u}|=0,
\end{aligned}
\right.\tag{7.1}
\end{equation}
where
\begin{align*}
&\dot{F}(v, \na_yv,\na_y^2v)
=\Delta\dot{v}-L(v, \na_yv,\na_y^2v)\dot{v},\\
&\dot{G}_1(v, \na_yv)=(G_1(u_0,\na_yu_0)u_0-G_1(v,\na_yv)u_0)+\f{\p\dot{v}}{\p n}-G_1(v,\na_yv)\dot{v},\\
&\dot{G}_2(v, \na_yv)=(G_2(u_0,\na_yu_0)u_0-G_2(v,\na_yv)u_0)+\p_{y_1}\dot{v}-G_2(v,\na_yv)\dot{v},
\end{align*}
with $\dot v=v-u_0$. Denote the mapping $J$ by $J(v)=u$, then we have the following lemma.

\vskip 0.2 true cm

{\bf Lemma 7.1.} {\it Suppose that the positive constants $\al$, $\dl$ and $\dl_0$ $(0<\al, \dl, \dl_0<1)$ are given in Proposition 6.2,
then there exists an $\ve_0$ such that for $\ve\in(0,\ve_0)$, $J$ is a mapping from $K$ to itself.}

{\bf Proof.} Set
$$
\dot{F}(v, \na_yv, \na_y^2v)\equiv \ss_{i=1}^7I_i,\eqno(7.2)
$$
where
\begin{align*}
&I_1=-\f{1}{q_0}\ss_{i,j=1}^3a_{ij}(\na_x\vp^-,\na_yv)\p_{x_ix_j}^2\vp^-,\\
&I_{i+1}=(1-A_{ii}(v,\na_yv))\p_{y_i}^2\dot{v},\q i=1,2,3,\\
&I_5=-2A_{12}(v,\na_yv)\p_{y_1y_2}^2\dot{v},\\
&I_6=-2A_{13}(v,\na_yv)\p_{y_1y_3}^2\dot{v},\\
&I_7=-2A_{23}(v,\na_yv)\p_{y_2y_3}^2\dot{v}.
\end{align*}

We now treat each $I_i$ separately. It follows from $\nabla^2\Phi_0=0$ and Lemma 2.4 that
$$
\|a_{ij}(\na_x\vp^-,\na_yv)\p_{x_ix_j}^2\vp^-\|_{4,\al;Q}^{(1-\dl,2-\dl_0)}
\le\|a_{ij}(\na_x\vp^-,\na_yv)\|_{4,\al;Q}^{(0,0)}\|\p_{x_ix_j}^2\vp^-\|_{6,\al;Q}^{(1-\dl,2-\dl_0)}\le C\ve,
$$
then we have
$$
\|I_1\|_{6,\al;Q}^{(-1-\dl,-\dl_0)}\le Cq_0^{-1}\ve.\eqno(7.3)
$$
To analyze $I_2$, we rewrite $I_2=\ss_{i=1}^6I_2^i$ with
\begin{align*}
&I_2^1=\bigg(1-\f{\ds a_{11}}{\ds \p_{y_1}v+b_0\p_{y_3}v}\bigg)\p_{1}^2\dot{v},\\
&I_2^2=-\f{\ds a_{22}(\p_{y_2}v)^2}{\ds (\p_{y_1}v+b_0\p_{y_3}v)^3}\p_{1}^2\dot{v},\\
&I_2^3=-\f{\ds a_{33}(\p_{y_3}v)^2}{\ds (\p_{y_1}v+b_0\p_{y_3}v)^3}\p_{1}^2\dot{v},\\
&I_2^4=\f{\ds 2a_{12}\p_{y_2}v}{\ds (\p_{y_1}v+b_0\p_{y_3}v)^2}\p_{1}^2\dot{v},\\
&I_2^5=\f{\ds 2a_{13}\p_{y_3}v}{\ds (\p_{y_1}v+b_0\p_{y_3}v)^2}\p_{1}^2\dot{v},\\
&I_2^6=-\f{\ds 2a_{23}\p_{y_2}v\p_{y_3}v}{\ds (\p_{y_1}v+b_0\p_{y_3}v)^2}\p_{1}^2\dot{v}.
\end{align*}

Notice that
\begin{align}
&\|1-a_{11}(\na_x\vp^-,\na_yv)\|_{4,\al;Q}^{(0,0)}=\|\f{\ds (\p_{x_1}\vp^--\p_{x_1}\Phi)^2}{\ds c^2(\na_x\vp^-,\na_yv)}\|_{4,\al;Q}^{(0,0)}\notag\\
\le&\|\f{\ds (\p_{x_1}\vp^--\p_{x_1}\Phi)^2}{\ds c^2(\na_x\vp^-,\na_yv)}-\f{\ds (\p_{x_1}\vp_0^--\p_{x_1}\Phi_0)^2}{\ds c^2(\na_x\vp^-,\na_yv)}\|_{4,\al;Q}^{(0,0)} +\|\f{\ds (\p_{x_1}\vp_0^--\p_{x_1}\Phi_0)^2}{\ds c^2(\na_x\vp^-,\na_yu_0)}\|_{4,\al;Q}^{(0,0)}\notag\\
\le&\|\f{\ds (\p_{x_1}\vp^--\p_{x_1}\Phi)^2}{\ds c^2(\na_x\vp^-,\na_yv)c^2(\na_x\vp_0^-,\na_yu_0)}\|_{4,\al;Q}^{(0,0)}\|c^2(\na_x\vp^-,\na_yv)-c^2(\na_x\vp_0^-, \na_yu_0)\|_{4,\al;Q}^{(-\dl,1-\dl_0)}\notag\\
&\qquad \qquad +\|\f{\ds (\p_{x_1}\vp^--\p_{x_1}\Phi)+(\p_{x_1}\vp_0-\p_{x_1}\Phi_0)}{\ds c^2(\na_x\vp_0^-,\na_yu_0)}\|_{4,\al;Q}^{(0,0)}\bigg(\|\p_{x_1}\vp^--\p_{x_1}\vp_0^-\|_{4,\al;Q}^{(-\dl,1-\dl_0)} \notag\\
&\qquad \qquad +\|\p_{x_1}\Phi-\p_{x_1}\Phi_0\|_{4,\al:Q}^{(-\dl,1-\dl_0)}\bigg)+Cq_0^{-\f{4}{\g-1}}\notag\\
\le&C(q_0^{-1}+q_0^{-\f{2}{\g-1}})(q_0\ve+\ve)+Cq_0^{-\f{4}{\g-1}},\tag{7.4}
\end{align}
and similarly,
$$
\|1-a_{ii}(\na_x\vp^-,\na_yv)\|_{4,\al;Q}^{(0,0)}\le C(q_0^{-1}+q_0^{-\f{2}{\g-1}})(q_0\ve+\ve)+Cq_0^{-\f{4}{\g-1}},\q i=2,3,\eqno(7.5)
$$
$$
\|a_{ij}(\na_x\vp^-,\na_yv)\|_{4,\al;Q}^{(0,0)}\le C(q_0^{-1}+q_0^{-\f{2}{\g-1}})(q_0\ve+\ve)+Cq_0^{-\f{4}{\g-1}}, \q 1\le i\neq j\le3.\eqno(7.6)
$$
In addition, one has
$$
\p_{y_1}u_0=1+O(q_0^{-\f{2}{\g-1}}),~~\p_{y_2}u_0=0,~~\p_{y_3}u_0=O(q_0^{-\f{2}{\g-1}}),~~\nabla_{y_iy_j}^2u_0=0,~~1\le i,j\le 3,
$$
this yields together with (7.4)
$$
\|I_2^1\|_{4,\al;Q}^{(1-\dl,2-\dl_0)}\le C(q_0^{-\f2{\g-1}}+\ve)\ve,
$$
and analogously,
$$
\|I_2^i\|_{4,\al;Q}^{(1-\dl,2-\dl_0)}\le C(q_0^{-\f2{\g-1}}+\ve)\ve,~~i=2,...,6.
$$
Therefore, we have
$$
\|I_2\|_{4,\al;Q}^{(1-\dl,2-\dl_0)}\le C(q_0^{-\f2{\g-1}}+\ve)\ve.\eqno(7.7)
$$
By the same method, we can arrive at
$$
\|I_i\|_{4,\al;Q}^{(1-\dl,2-\dl_0)}\le C(q_0^{-\f2{\g-1}}+\ve)\ve,~~i=3,...,7.\eqno(7.8)
$$
Thus, substituting (7.3) and (7.7)-(7.8) into (7.2) yields
$$
\|\dot{F}(v,\na_yv,\na_y^2v)\|_{4,\al;Q}^{(1-\dl,2-\dl_0)}\le C(q_0^{-\f2{\g-1}}+\ve)\ve.\eqno(7.9)
$$

On the other hand, it follows from a direct computation that
\begin{align*}
\dot{G}_1(v,&\na_yv)=\f{1}{q_0}(\sin\th_0\p_{y_3}u_0+\cos\th_0\p_{y_1}u_0)\big((\p_{x_3}\vp-\p_{x_3}\vp_0)-b_0(\p_{x_1}\vp^--\p_{x_1}\vp_0^-)\big)\\
&+\f{\sin\th_0}{q_0}\big(b_0(q_0-\p_{x_1}\vp^-)+\p_{x_3}\vp^-\big)\p_{y_3}\dot{v}+\f{\cos\th_0}{q_0}\p_{x_3}\vp^-\p_{y_1}\dot{v}
-\f{\sin\th_0}{q_0}(\p_{x_1}\vp^--q_0)\p_{y_1}\dot{v},
\end{align*}
which yields
$$
\|\dot{G}_1(v, \na_y v)\|_{5,\al}^{(-\dl,1-\dl_0)}\le Cq_0^{-1}\ve.\eqno(7.10)
$$
Analogously,
$$
\|\dot{G}_2(v, \na_y v)\|_{5,\al}^{(-\dl,1-\dl_0)}\le Cq_0^{-1}\ve.\eqno(7.11)
$$

For appropriately large $q_0$ and small $\ve>0$, then Proposition 6.2 implies that there exists a unique solution $\dot{u}\in H_{6,\al}^{(-1-\dl,-\dl_0)}$ to (7.1) such that
$$
\|\dot{u}\|_{6,\al;Q}^{(-1-\dl,-\dl_0)}\le C(q_0^{-1}+\ve)\ve\le \ve,
$$
which means that mapping $J$ is from $K$ to itself. \qed

\vskip 0.2 true cm

{\bf Lemma 7.2.} {\it Under the assumptions of Lemma 7.1, the mapping $J$ is a contractible mapping from $K$ to itself.}

{\bf Proof.} Taking $v_1, v_2\in K$. Let $u_i=Jv_i$ and $\dot{u}_i=u_i-u_0$ in $Q$, then we have
\begin{equation}
\left\{
\begin{aligned}
&\Delta(u_2-u_1)=\dot{F}(v_1,\na_yv_1,\na_y^2v_1)-\dot{F}(v_2,\na_yv_2,\na_y^2v_2)\qquad\qquad \q \text{in}\q Q,\\
&\p_n(u_2-u_1)=\dot{G}_1(v_1,\na_yv_1)-\dot{G}_1(v_2,\na_yv_2)\qquad\qquad\qquad\q\text{on}\q y_3=b_0y_1,\\
&\p_{y_1}(u_2-u_1)=\dot{G}_2(v_1,\na_yv_1)-\dot{G}_2(v_2,\na_yv_2)\qquad\qquad\qquad\q\text{on}\q y_1=0,\\
&(u_2-u_1)(0,0,0)=0,\\
&\ds\lim_{y_1+y_3\rightarrow\infty}|\nabla_y(u_2-u_1)|=0.
\end{aligned}
\right.\tag{7.12}
\end{equation}

As in Lemma 7.1, a direct computation yields
\begin{align}
\|\dot{F}(v_2,\na_yv_2, \na_y^2v_2)-&\dot{F}(v_1,\na_yv_1,\na_y^2v_1)\|_{4,\al;Q}^{(1-\dl,2-\dl_0)}\le C\ve\|v_2-v_1\|_{6,\al;Q}^{(-1-\dl,-\dl_0)}\notag\\ &+C(q_0^{-1}+q_0^{-\f{2}{\g-1}}+\ve)\|u_2-u_1\|_{6,\al;Q}^{(-1-\dl,-\dl_0)}\tag{7.13}
\end{align}
and
\begin{align}
&\|\dot{G}_i(v_2,\na_yv_2)-\dot{G}_i(v_1,\na_yv_1)\|_{5,\al;Q}^{(-\dl,1-\dl_0)}\notag\\
\le& C\ve\|v_2-v_1\|_{6,\al;Q}^{(-1-\dl,-\dl_0)} +C(q_0^{-1}+q_0^{-\f{2}{\g-1}}
+\ve)\|u_2-u_1\|_{6,\al;Q}^{(-1-\dl,-\dl_0)},\qquad i=1,2.\tag{7.14}
\end{align}
By Proposition 6.2, we have
$$
\|u_2-u_1\|_{6,\al;Q}^{(-1-\dl,-\dl_0)}\le C\ve\|v_2-v_1\|_{6,\al;Q}^{(-1-\dl,-\dl_0)} +C(q_0^{-1}+q_0^{-\f{2}{\g-1}}+\ve)\|u_2-u_1\|_{6,\al;Q}^{(-1-\dl,-\dl_0)}.\eqno(7.15)
$$
Choosing appropriately large $q_0$ and small $\ve_0$ yields
$$
\|Jv_2-Jv_1\|_{6,\al;Q}^{(-1-\dl,-\dl_0)}\le \f12\|v_2-v_1\|_{6,\al;Q}^{(-1-\dl,-\dl_0)},
$$
which means that $J$ is a contractible mapping.\qed

We now prove Theorem 3.1.

{\bf Proof of Theorem 3.1.}  By Lemma 7.1 and Lemma 7.2, we know that the mapping $u=Jv$ has a unique fixed point
in the space $H_{6,\al}^{(-1-\dl,-\dl_0)}(Q)$,
which implies that Theorem 3.1 is shown. \qed

Based on Theorem 3.1, we can show Theorem 1.1.

{\bf Proof of Theorem 1.1.} By Theorem 3.1, one knows that the problem (1.19) admits a unique solution
$u\in H_{6,\al}^{(-1-\dl,-\dl_0)}(Q)$. Since only the condition $u(0, y_2^0, 0)=0$ other than $u(0, y_2, 0)\equiv 0$
for all $y_2\in\Bbb R$ is applied in order to solve
the 3-D attached strong oblique shock problem (1.18), then (1.18) is obviously overdetermined. \qed

\vskip 0.8 true cm
\centerline{\bf Appendix}
\vskip 0.4 true cm
{\bf Lemma A.1.} {\it For the term $I_1$ defined in (4.28), we have
$$
\|I_1\|_{0,\alpha;Q_L}^{(0,-\delta_0)}\leqslant  C \|f\|_{0,\alpha;Q}^{(1-\delta,2-\delta_0)},\eqno(A.1)
$$
where the generic positive constant $C$ is independent of $L$.}

{\bf Proof.} To estimate $\|I_1\|_{0,\alpha}^{(0,-\delta_0)}$, by the definition of $\|I_1\|_{0,\alpha}^{(0,-\delta_0)}$,
we need to study the two cases
including $0<r\le 1$ and $r\ge 1$ separately.

{\bf Case i. $0<r \leqslant 1$}

By the first inequalities in (4.33) and (4.37), and the fact of $0<\delta<\sqrt{\mu_1}-1$, then we get
\begin{align*}
&\sum_{m=1}^{\infty}|R_{m0}(r)|\\
&\leqslant \sum_{m=1}^{\infty} r^{\m}L^{-2\m} \biggl\{\int_0^1 \eta^{2\m-1} \Bigl(\int_{\eta}^1 \xi^{1-\m}|{f_L}_{m0}(\xi)| d \xi\\
&\quad  + \int_1^L \xi^{1-\m}|{f_L}_{m0}(\xi)| d \xi \Bigr)d\eta
+\int_{1}^L  \eta^{2\m-1}\int_{\eta}^L \xi^{1-\m}|{f_L}_{m0}(\xi)| d \xi d\eta \biggr\}\\
&\quad +r^{-\m}\int_0^r \eta^{2\m-1}\Big(\int_{\eta}^1 \xi^{1-\m}|{f_L}_{m0}(\xi)| d \xi
+\int_1^L \xi^{1-\m}|{f_L}_{m0}(\xi)| d \xi\Big) d\eta\\
&\leqslant C \|f\|_{0,\alpha;Q}^{(1-\delta,2-\delta_0)} \sum_{m=1}^{\infty} \biggl\{r^{\m}L^{-2\m}
\biggl[\int_0^1 \eta^{2\m-1} \Bigl(\int_{\eta}^1 \xi^{\delta-\m} d \xi
+ \int_1^L \xi^{\delta_0-1-\m} d \xi \Bigr)d\eta\\
&\quad +\int_{1}^L  \eta^{2\m-1}\int_{\eta}^L \xi^{\delta_0-1-\m} d \xi d\eta \biggr]
+r^{-\m}\int_0^r \eta^{2\m-1}\Big(\int_{\eta}^1 \xi^{\delta-\m} d \xi\\
&\quad +\int_1^L \xi^{\delta_0-1-\m} d \xi\Big) d\eta \biggr\}\\
\end{align*}
\begin{align}
&\leqslant C \|f\|_{0,\alpha;Q}^{(1-\delta,2-\delta_0)} \sum_{m=1}^{\infty} \biggl\{r^{\m}L^{-2\m}\biggl[\int_0^1 \eta^{2\m-1} \Bigl(\frac{\eta^{1+\delta-\m}-1}{\m-\delta-1}
+ \frac{1-L^{\delta_0-\m}}{\m-\delta_0} \Bigr)d\eta\notag \\
&\quad +\int_{1}^L  \eta^{2\m-1} \frac{\eta^{\delta_0-\m}-L^{\delta_0-\m}}{\m-\delta_0} d\eta \biggr]
+r^{-\m} \int_0^r \eta^{2\m-1}\Big( \frac{\eta^{1+\delta-\m}-1}{\m-\delta-1}\notag\\
&\quad +\frac{1-L^{\delta_0-\m}}{\m-\delta_0} \Big) d\eta \biggr\}\notag\\
&\leqslant C \|f\|_{0,\alpha;Q}^{(1-\delta,2-\delta_0)} \sum_{m=1}^{\infty} \biggl\{ \frac{r^{\m}L^{-2\m}+r^{1+\delta}}{(\m-\delta-1)(\m+\delta+1)}
+ \frac{r^{\m}(L^{-2\m}+1)}{2\m(\m-\delta_0)}\notag\\
&\quad +\frac{r^{\m}L^{\delta_0-\m}}{(\m-\delta_0)(\m+\delta_0)} \biggr\}\notag\\
&\leqslant C\|f\|_{0,\alpha;Q}^{(1-\delta,2-\delta_0)} \sum_{m=1}^{\infty} \biggl\{ \frac{1}{\mu_m-(\delta+1)^2}
+ \frac{1}{\m(\m-\delta_0)}
+\frac{1}{\mu_m-\delta_0^2} \biggr\}\notag\\
&\leqslant C \|f\|_{0,\alpha;Q}^{(1-\delta,2-\delta_0)}.\tag {A.2}
\end{align}

{\bf Case ii. $r \geqslant 1$}

In this case, we have
\begin{align*}
|R_{m0}(r)|
&\leqslant r^{\m} L^{-2\m} \int_0^L \eta^{2\m-1} \int_{\eta}^L  \xi^{1-\m} |f_{m0}(\xi)| d \xi d\eta\\
 &\quad + r^{-\m} \int_0^r \eta^{2\m-1} \int_{\eta}^L  \xi^{1-\m}|f_{m0}(\xi)| d \xi d\eta \\
&\equiv A_{m1} +A_{m2}.
\end{align*}

For $A_{m1}$ with $m\in\Bbb N$, by the same method as in Case i, we can get
\begin{align}
\sum_{m=1}^{\infty} A_{m1}
& \leqslant C \|f\|_{0,\alpha;Q}^{(1-\delta,2-\delta_0)} \sum_{m=1}^{\infty} \biggl(\frac{r^{\m}L^{-2\m}}{\mu_m-(\delta+1)^2}
+\frac{r^{\m}L^{-2\m}}{2\m(\m-\delta_0)}
+\frac{r^{\m}L^{\delta_0-\m}}{\mu_m-\delta_0^2} \biggr)\notag\\
&\leqslant C\|f\|_{0,\alpha;Q}^{(1-\delta,2-\delta_0)} \sum_{m=1}^{\infty} \biggl(\frac{1}{\mu_m-(\delta+1)^2}
+\frac{1}{\m(\m-\delta_0)}
+\frac{r^{\dl_0}}{\mu_m-\delta_0^2} \biggr)\notag\\
&\leqslant C \|f\|_{0,\alpha;Q}^{(1-\delta,2-\delta_0)}r^{\dl_0}.\tag{A.3}
\end{align}

For $A_{m2}$ with $m\in\Bbb N$,
\begin{align*}
&\sum_{m=1}^{\infty} A_{m2}\\
&\leqslant  C\|f\|_{0,\alpha;Q}^{(1-\delta,2-\delta_0)} \sum_{m=1}^{\infty} r^{-\m} \biggl\{ \int_0^1 \eta^{2\m-1}
\Big( \int_{\eta}^1 \xi^{1-\m} \xi^{\delta-1} d \xi + \int_1^L \xi^{1-\m} \xi^{\delta_0-2} d \xi\Big)  d\eta \\
& \qquad + \int_1^r \eta^{2\m-1} \int_{\eta}^L  \xi^{1-\m} \xi^{\delta_0-2} d \xi d\eta \biggr\}\\
\end{align*}

\begin{align}
&\leqslant C\|f\|_{0,\alpha;Q}^{(1-\delta,2-\delta_0)} \sum_{m=1}^{\infty} r^{-\m} \biggl\{ \int_0^1 \eta^{2\m-1}
\Big( \frac{\eta^{1+\delta-\m}-1}{\m-\delta-1} + \frac{1-L^{\delta_0-\m}}{\m-\delta_0} \Big) d\eta\notag\\
& \qquad + \int_1^r \eta^{2\m-1} \frac{\eta^{\delta_0-\m}-L^{\delta_0-\m}}{\m-\delta_0}d\eta \biggr\}\notag\\
&\leqslant C\|f\|_{0,\alpha;Q}^{(1-\delta,2-\delta_0)} \sum_{m=1}^{\infty} \biggl\{ \frac{r^{-\m}}{(\m-\delta-1)(\m+\delta+1)}
+ \frac{r^{-\m}}{2\m(\m-\delta_0)}\notag\\
&\qquad  +\frac{r^{\delta_0}}{(\m-\delta_0)(\m+\delta_0)} \biggr\}\notag\\
& \leqslant C\|f\|_{0,\alpha;Q}^{(1-\delta,2-\delta_0)}r^{\delta_0} \sum_{m=1}^{\infty} \biggl(\frac{1}{\mu_m-(\delta+1)^2}
+ \frac{1}{\m(\m-\delta_0)}
+\frac{1}{\mu_m-\delta_0^2} \biggr)\notag\\
&\leqslant C \|f\|_{0,\alpha;Q}^{(1-\delta,2-\delta_0)} r^{\delta_0}.\tag {A.4}
\end{align}
Thus, collecting the estimates (A.2)-(A.4) yields Lemma A.1.\qed

\vskip 0.2 true cm
{\bf Lemma A.2.} {\it For the term $I_2$ defined in (4.28), we have
$$
\|I_2\|_{0,0;Q_L}^{(0,0)}\leqslant C \|f\|_{1,\alpha;Q}^{(1-\delta,2-\delta_0)},
$$
where $C>0$ is independent of $L$.
}

{\bf Proof.} By the expression of $I_2$, it suffices to only treat  $\ds\sum_{n=1}^{\infty}|R_{0n}^{(1)}(r)|$ since   $\ds\sum_{n=1}^{\infty}|R_{0n}^{(2)}(r)|$
can be analogously estimated. We write
$$
R_{0n}^{(1)}(r)=B_{1}^n-B_{2}^n,\eqno(A.5)
$$
where
\begin{align*}
&B_{1}^n=\frac{K_{0}^{'}(nL)}{I_{0}^{'}(nL)}I_{0}(nr)\int_0^L sI_{0}(ns){f_L}_{0n}^{(i)}(s) d s,\\
&B_{2}^n=I_{0}(nr)\int_r^L sK_{0}(ns){f_L}_{0n}^{(i)}(s) d s
+K_{0}(nr)\int_0^r sI_{0}(ns){f_L}_{0n}^{(i)}(s) d s.
\end{align*}

Next we deal with $\ds\sum_{n=1}^{\infty} |B_{1}^n|$
and $\ds\sum_{n=1}^{\infty} |B_{2}^n|$ respectively.

By (iii) and (v) in Lemma 2.3, and the inequalities in (4.34) and (4.38), we obtain
\begin{align*}
&\ds\sum_{n=1}^{\infty} |B_{1}^n|
\leqslant \sum_{n=1}^{\infty} \frac{\sqrt{\frac{\pi}{2nL}}e^{-nL}}{\frac{e^{nL}}{\sqrt{2\pi nL}}} e^{nr}
\biggl(\int_0^{\frac 1n}+\int_{\frac 1n}^1+\int_1^{L}\biggr)s I_{0}(ns) |{f_L}_{0n}^{(1)}(s)| d s\\
&\quad \leqslant C\sum_{n=1}^{\infty} e^{nr-2nL}  \biggl\{\int_0^{\frac 1n}  e^{ns}s^{\delta} ds \|f\|_{0,\alpha;Q}^{(1-\delta,2-\delta_0)}+
\biggl(\int_{\frac 1n}^1 e^{ns}s^{\delta-1} d s\\
&\qquad +\int_1^{L} e^{ns}s^{\delta_0-2} d s \biggr) n^{-1} \|f\|_{1,\alpha;Q}^{(1-\delta,2-\delta_0)} \biggr\}\\
\end{align*}

\begin{align}
&\quad \leqslant C \sum_{n=1}^{\infty} e^{nr-2nL}  \biggl\{n^{-1-\delta}(e-1) \|f\|_{0,\alpha;Q}^{(1-\delta,2-\delta_0)}+
\biggl( n^{-1-\delta}(e^{n}-e)\notag\\
&\qquad +n^{-2} (e^{nL}-e^n)\biggr) \|f\|_{1,\alpha;Q}^{(1-\delta,2-\delta_0)}\biggr\}\notag\\
&\quad \leqslant \sum_{n=1}^{\infty} \biggl\{n^{-1-\delta}e^{-nL} \|f\|_{0,\alpha;Q}^{(1-\delta,2-\delta_0)}+
\bigl( n^{-1-\delta} e^{n(1-L)}+n^{-2}\bigr) \|f\|_{1,\alpha;Q}^{(1-\delta,2-\delta_0)}\biggr)\notag\\
&\quad \leqslant C \|f\|_{1,\alpha;Q}^{(1-\delta,2-\delta_0)}.\tag{A.6}
\end{align}

To estimate $\ds\sum_{n=1}^{\infty} |B_{2}^n|$, we will divide this procedure into the following three cases.

{\bf Case i. $0<nr \leqslant 1$.}

In this case, it follows from (iii) in Lemma 2.3 that
\begin{align}
&\sum_{nr\leqslant 1} |B_{2}^n|\leqslant \sum_{nr\leqslant 1}
\biggl\{e^{nr} \biggl( \int_r^{\frac 1n} +\int_{\frac 1n}^1
+\int_1^{L}\biggr)s K_{0}(ns) |{f_L}_{0n}^{(1)}(s)| d s \notag\\
&\qquad + K_{0}(nr)\int_0^r sI_{0}(ns) |{f_L}_{0n}^{(1)}(s)|ds \biggr\} \notag\\
& \leqslant C \|f\|_{0,\alpha;Q}^{(1-\delta,2-\delta_0)} \sum_{nr\leqslant 1} \biggl\{e^{nr} \biggl(\int_r^{\frac 1n}  \frac{e^{-ns}}{\sqrt{2ns}}s^{\delta} ds +
  \int_{\frac 1n}^1 \frac{e^{-ns}}{\sqrt{2ns}}s^{\delta} d s +\int_1^{L} \frac{e^{-ns}}{\sqrt{2ns}}s^{\delta_0-1} d s \biggr) \notag\\
&\qquad  + \frac{e^{-nr}}{\sqrt{2nr}} \int_0^r e^{ns}s^{\delta} ds \biggr\} \notag\\
& \leqslant C \|f\|_{0,\alpha;Q}^{(1-\delta,2-\delta_0)} \sum_{nr\leqslant 1}  \biggl(\frac{n^{-1-\delta}}{\delta+\frac 12} +
\frac{n^{-1-\min(\frac 12,\delta)}}{e} +\frac{n^{-1-\frac 12}}{e^{n}} +\frac{n^{-1-\delta}}{1+\delta} \biggr) \notag\\
& \leqslant C \|f\|_{0,\alpha;Q}^{(1-\delta,2-\delta_0)}.\tag{A.7}
\end{align}

{\bf Case ii. $1\leqslant nr$ and $r\leqslant 1$}

We get
\begin{align*}
&\sum_{1\le nr\leqslant n} |B_{2}^n|\leqslant \sum_{1\leqslant nr \leqslant n}
\biggl\{e^{nr} \biggl(\int_{r}^1 K_{0}(ns) |{f_L}_{0n}^{(1)}(s)| d s +\int_1^{L} s K_{0}(ns) |{f_L}_{0n}^{(1)}(s)| d s \biggr)\\
&\quad + K_{0}(nr)\biggl(\int_0^{\frac 1n} sI_{0}(ns) |{f_L}_{0n}^{(1)}(s)|ds + \int_{\frac 1n}^r sI_{0}(ns) |{f_L}_{0n}^{(1)}(s)|ds \biggr)\biggr\}\\
&\leqslant C\|f\|_{0,\alpha;Q}^{(1-\delta,2-\delta_0)} \sum_{1\leqslant nr \leqslant n} \biggl\{e^{nr} \biggl( n^{-\frac 12}\int_{r}^1 e^{-ns}s^{\delta-\frac 12} d s
+ n^{-\frac 12}\int_1^{L} e^{-ns}s^{\delta_0-\frac 32} d s \biggr)\\
&\quad + \frac{e^{-nr}}{\sqrt{2nr}} \biggl( \int_0^{\frac 1n} e^{ns}s^{\delta} ds + \int_{\frac 1n}^r e^{ns}s^{\delta}ds \biggr) \biggr\} \\
\end{align*}

\begin{align}
& \leqslant C\|f\|_{0,\alpha;Q}^{(1-\delta,2-\delta_0)} \sum_{1\leqslant nr \leqslant n} \biggl\{ \biggl( n^{-1-\min(\frac 12,\delta)}
+ n^{-\frac 32} \biggr) +  \biggl( e^{-nr} n^{-1-\delta} e
+ r^{\delta-\frac 12}n^{-\frac 32}\biggr) \biggr\} \notag\\
& \leqslant C \|f\|_{0,\alpha;Q}^{(1-\delta,2-\delta_0)}.\tag{A.8}
\end{align}

{\bf Case iii. $ 1\leqslant r \leqslant L$.}

At this time, we obtain
\begin{align}
&\sum_{n=1}^{\infty} |B_{2}^n|\leqslant \sum_{n=1}^{\infty}
\biggl\{e^{nr} \int_{r}^L K_{0}(ns) |{f_L}_{0n}^{(1)}(s)| d s
+ K_{0}(nr)\biggl(\int_0^{\frac 1n} +\int_{\frac 1n}^1+\int_1^r \biggr)
sI_{0}(ns) |{f_L}_{0n}^{(1)}(s)|ds \biggr\} \notag\\
& \leqslant C\|f\|_{0,\alpha;Q}^{(1-\delta,2-\delta_0)} \sum_{n=1}^{\infty} \biggl\{e^{nr} \int_r^{L} \frac{e^{-ns}}{\sqrt{2ns}}s^{\delta_0-1} d s
+ \frac{e^{-nr}}{\sqrt{2nr}}\biggl(\int_0^{\frac 1n} e^{ns}s^{\delta} + \int_{\frac 1n}^1 e^{ns}s^{\delta}ds \notag\\
&\qquad +\int_1^r e^{ns} s^{\delta_0-1}ds \biggr) \biggr\} \notag\\
&\leqslant C\|f\|_{0,\alpha;Q}^{(1-\delta,2-\delta_0)} \sum_{n=1}^{\infty} n^{-\frac 12}
\biggl\{e^{nr}\int_r^{L} e^{-ns}s^{\delta_0-\frac 32} d s + e^{-nr}\biggl( \int_0^{\frac 1n} e^{ns}s^{\delta} ds + \int_{\frac 1n}^1 e^{ns}s^{\delta}ds \notag\\
&\qquad  + \int_1^r e^{ns} s^{\delta_0-1}ds \biggr) \biggr\} \notag\\
& \leqslant C\|f\|_{0,\alpha;Q}^{(1-\delta,2-\delta_0)} \sum_{n=1}^{\infty}
\bigl(n^{-\frac 32}+ n^{-\delta-\frac 32} e^{-n}\bigr) \notag\\
& \leqslant C \|f\|_{0,\alpha;Q}^{(1-\delta,2-\delta_0)}.\tag{A.9}
\end{align}
Consequently, collecting (A.6)-(A.9) yields Lemma A.2. \qed

\vskip 0.2 true cm
{\bf Lemma A.3.} {\it For the term $R_{00}(r)$ defined in (4.27), we have
$$
\|R_{00}\|_{0,0;Q_L}^{(0,-\delta_0)} \leqslant C \|f\|_{1,\alpha;Q}^{(1-\delta,2-\delta_0)},\eqno(A.10)
$$
where $C>0$ is independent of $L$.}

{\bf Proof.} Noting that from (4.23)
$$
|R_{00}(r)|\leqslant \sum_{n=1}^{\infty}|R_{0n}^{(2)}(0)|+\int_0^r \eta^{-1}\int_0^{\eta} \xi |{f_L}_{00}(\xi)|d \xi d \eta\equiv A_{01}+A_{02}.
$$
By $I_0(0)=1$ and Lemma A.2, we know that
\begin{align}
A_{01} &\leqslant \sum_{n=1}^{\infty} \biggl(\biggl|\frac{K_{0}^{'}(nL)}{I_{0}^{'}(nL)}\biggr|\int_0^L sI_{0}(ns)|{f_L}_{0n}^{(2)}(s)| d s + \int_r^L sK_{0}(ns)|{f_L}_{0n}^{(2)}(s)| d s\biggr) \notag\\
& \leqslant C \|f\|_{1,\alpha;Q}^{(1-\delta,2-\delta_0)}.\tag{A.11}
\end{align}

To estimate $A_{02}$, we divide this process into the following two cases.

\vskip 0.2 true cm

{\bf Case i: $0<r\leqslant 1$}

By (4.32), we have
\begin{align}
A_{02} &\leqslant C\|f\|_{0,\alpha;Q}^{(1-\delta,2-\delta_0)} \int_0^r \eta^{-1}\int_0^{\eta} \xi^{\delta}d \xi d \eta \notag\\
&\leqslant C\|f\|_{0,\alpha;Q}^{(1-\delta,2-\delta_0)} \frac{r^{1+\delta}}{(1+\delta)^2} \notag\\
&\leqslant C\|f\|_{0,\alpha;Q}^{(1-\delta,2-\delta_0)}.\tag{A.12}
\end{align}

{\bf Case ii: $1<r\leqslant L$}

By (4.36), we get
\begin{align}
A_{02} & \leqslant \int_0^1 \eta^{-1}\int_0^{\eta} \xi |{f_L}_{00}(\xi)|d \xi d \eta+ \int_1^r \eta^{-1}\biggl(\int_0^1 \xi |{f_L}_{00}(\xi)|d \xi + \int_1^{\eta} \xi |{f_L}_{00}(\xi)|d \xi \biggr)d \eta \notag\\
& \leqslant C\|f\|_{0,\alpha;Q}^{(1-\delta,2-\delta_0)}
\biggl\{\int_0^1 \eta^{-1}\int_0^{\eta} \xi^{\delta} d \xi d \eta +
\int_1^r \eta^{-1}\biggl(\int_0^1 \xi^{\delta} d \xi + \int_1^{\eta} \xi^{\delta_0-1} d \xi \biggr)d \eta \biggr\} \notag\\
&\leqslant C\|f\|_{0,\alpha;Q}^{(1-\delta,2-\delta_0)} \biggl(\frac{1}{(1+\delta)^2} +
\frac{\ln r}{(1+\delta)^2} +  \frac{r^{\delta_0}}{\delta_0^2} \biggr) \notag\\
&\leqslant C\|f\|_{0,\alpha;Q}^{(1-\delta,2-\delta_0)} r^{\delta_0}.\tag {A.13}
\end{align}
Thus, combining (A.11)-(A.13) yields (A.10).\qed

\vskip 0.2 true cm

{\bf Lemma A.4.} {\it For the term $I_3$ defined in (4.28), we have
$$
\|I_3\|_{0,0;Q_L}^{(0,0)} \leqslant C\|f\|_{1,\alpha;Q}^{(1-\delta,2-\delta_0)},\eqno(A.14)
$$
where $C>0$ is independent of $L$.}

{\bf Proof.} Since it follows from (vii)-(viii) in Lemma 2.3 that
the Bessel functions $I_{\m}(x)$ and $K_{\m}(x)$ have different properties for  $m\leqslant M$ and $m>M$
($M\in\Bbb N$ is some suitably large integer), then we require to divide the process of estimating $I_3$
into the following four cases.

\vskip 0.2 true cm

{{\bf Case i.} $m\leqslant M$, and $nr \leqslant 1$ with $0<r \leqslant 1$.}

At this time, by (iv)-(vii) in Lemma 2.3 and the inequalities in (4.35) and (4.39), we derive that
\begin{align*}
&|R_{mn}^{(1)}(r)|\leqslant  I_{\m}(nr)\biggl(\frac{K_{\m}^{'}(nL)}{I_{\m}^{'}(nL)}
\int_0^L sI_{\m}(ns)|{f_L}_{mn}^{(1)}(s)| d s+\int_r^L sK_{\m}(ns)|{f_L}_{mn}^{(1)}(s)| d s\biggr)\\
&\qquad +K_{\m}(nr)\int_0^r sI_{\m}(ns)|{f_L}_{mn}^{(1)}(s)| d s\\
\end{align*}

\begin{align*}
& \leqslant C\|f\|_{0,\alpha;Q}^{(1-\delta,2-\delta_0)}\frac{e^{nr} (\frac {nr}2)^{\m}}{\Gamma({\m}+1)}
\biggl\{ e^{-2nL} \biggl[ \int_0^{\frac 1n} \frac{e^{ns} (\frac {ns}2)^{\m} }{\Gamma({\m}+1)} s^{\delta} ds
+ \int_{\frac 1n}^1 \frac{e^{ns}}{\sqrt{2\pi ns}} s^{\delta} ds\\
&\qquad + \int_1^L \frac{e^{ns}}{\sqrt{2\pi ns}} s^{\delta_0-1} ds \biggr]
+ \biggl[\int_r^{\frac 1n} \frac{e^{ns} \Gamma({\m}) 2^{\m-1}}{(ns)^{\m}} s^{\delta} ds
+ \int_{\frac 1n}^1 \frac{\sqrt \pi e^{-ns}}{\sqrt{2ns}} s^{\delta} ds\\
&\qquad + \int_1^L \frac{\sqrt \pi e^{-ns}}{\sqrt{2ns}} s^{\delta_0-1} ds \biggr] \biggr\}
+ C\|f\|_{0,\alpha;Q}^{(1-\delta,2-\delta_0)}\frac{e^{nr} \Gamma({\m}) 2^{\m-1}}{(nr)^{\m}}
\int_0^r \frac{e^{ns} (\frac {ns}2)^{\m} }{\Gamma({\m}+1)} s^{\delta} ds\\
& \leqslant C \|f\|_{0,\alpha;Q}^{(1-\delta,2-\delta_0)} \biggl\{ \frac 1{e^{2nL} 2^{\m}}
\biggl[ \int_0^{\frac 1n} n^{\m} s^{\delta+\m} ds + n^{-\frac 12}\int_{\frac 1n}^1 e^{ns} s^{\delta-\frac 12} ds\\
&\qquad  +
n^{-\frac 12} \int_1^L e^{ns} s^{\delta_0-\frac 32} ds \biggr]
+ \biggl[ \int_r^{\frac 1n} \frac{e^{ns} s^{\delta}}{\m}ds + n^{-\frac 12}\int_{\frac 1n}^1 \frac{e^{-ns} s^{\delta-\frac 12}}{2^{\m}} ds + n^{-\frac 12}\int_1^L \frac{e^{-ns} s^{\delta_0-\frac 32}}{ 2^{\m}} ds \biggr]\\
&\qquad + \int_0^r \frac{e^{ns} s^{\delta}}{\m}ds \biggr\}\\
& \leqslant C \|f\|_{0,\alpha;Q}^{(1-\delta,2-\delta_0)} \frac {n^{-1-\min{(\frac 12,\delta)}}}{\m},
\end{align*}
which yields
$$
\sum_{m=1}^{M}\sum_{nr\leqslant 1 \atop 0<r \leqslant 1} |R_{mn}^{(1)}(r)| \leqslant  C \|f\|_{0,\alpha;Q}^{(1-\delta,2-\delta_0)}. \eqno(A.15)
$$

{{\bf Case ii.} $m\leqslant M$, and $nr \geqslant 1$ with $0<r \leqslant 1$.}

In this case, we have
\begin{align*}
&|R_{mn}^{(1)}(r)|\leqslant  I_{\m}(nr)\biggl(\frac{K_{\m}^{'}(nL)}{I_{\m}^{'}(nL)}
\int_0^L sI_{\m}(ns)|{f_L}_{mn}^{(1)}(s)| d s+\int_r^L sK_{\m}(ns)|{f_L}_{mn}^{(1)}(s)| d s\biggr)\\
&\qquad +K_{\m}(nr)\int_0^r sI_{\m}(ns)|{f_L}_{mn}^{(1)}(s)| d s\\
& \leqslant C\|f\|_{0,\alpha;Q}^{(1-\delta,2-\delta_0)}\frac{e^{nr}}{\sqrt{2\pi nr}}
 \biggl\{ e^{-2nL} \biggl[ \int_0^{\frac 1n} \frac{e^{ns} (\frac {ns}2)^{\m} }{\Gamma({\m}+1)} s^{\delta} ds
+ \int_{\frac 1n}^1 \frac{e^{ns}}{\sqrt{2\pi ns}} s^{\delta} ds\\
&\qquad + \int_1^L \frac{e^{ns}}{\sqrt{2\pi ns}} s^{\delta_0-1} ds \biggr]
+ \biggl[\int_r^1 \frac{\sqrt \pi e^{-ns}}{\sqrt{2ns}} s^{\delta} ds + \int_1^L \frac{\sqrt \pi e^{-ns}}{\sqrt{2ns}} s^{\delta_0-1} ds \biggr] \biggr\}\\
&\qquad + \|f\|_{0,\alpha;Q}^{(1-\delta,2-\delta_0)}\frac{\sqrt \pi e^{-nr}}{\sqrt{2nr}}
\biggl[ \int_0^{\frac 1n} \frac{e^{ns} (\frac {ns}2)^{\m} }{\Gamma({\m}+1)} s^{\delta} ds +
\int_{\frac 1n}^r \frac{e^{ns}}{\sqrt{2\pi ns}} s^{\delta} ds \biggr] \\
\end{align*}

\begin{align*}
& \leqslant C \|f\|_{0,\alpha;Q}^{(1-\delta,2-\delta_0)} \biggl\{\frac{ n^{-1-\delta} e^{-n}} {2^{\m}\Gamma({\m}+1)}
+ n^{-1-\min{(\frac 12,\delta)}}+ n^{-\frac 32}+ n^{-1-\min{(\frac 12,\delta)}} e^{-nr}+ n^{-\frac 32}e^{-n}\\
& \qquad  +\frac{ n^{-1-\delta}} {2^{\m}\Gamma({\m}+1)}+ n^{-1-\min{(\frac 12,\delta)}} \biggr\}\\
& \leqslant C \|f\|_{0,\alpha;Q}^{(1-\delta,2-\delta_0)}  n^{-1-\min{(\frac 12,\delta)}},
\end{align*}
which derives
$$
\sum_{m=1}^{M}\sum_{nr\geqslant 1 \atop 0<r \leqslant 1} |R_{mn}^{(1)}(r)| \leqslant  C \|f\|_{0,\alpha;Q}^{(1-\delta,2-\delta_0)}. \eqno(A.16)
$$

{{\bf Case iii.} $m\leqslant M$ and $1 \leqslant r \leqslant L$.}

Set
$$
R_{mn}^{(1)}(r)=B_{mn}^1(r)-B_{mn}^2(r)-B_{mn}^3(r),\eqno(A.17)
$$
where
\begin{align*}
&B_{mn}^1(r)=I_{\m}(nr)\Big(\frac{K_{\m}^{'}(nL)}{I_{\m}^{'}(nL)}\int_0^L sI_{\m}(ns){f_L}_{mn}^{(1)}(s) d s\biggr),\\
&B_{mn}^2(r)=I_{\m}(nr)\int_r^L sK_{\m}(ns){f_L}_{mn}^{(1)}(s) d s\\
&B_{mn}^3(r)=K_{\m}(nr)\int_0^r sI_{\m}(ns){f_L}_{mn}^{(1)}(s) d s.
\end{align*}

By the same method as in Case ii and the fact of $1 \leqslant r \leqslant L$, one has
$$
|B_{mn}^1(r)|\leqslant C \|f\|_{0,\alpha;Q}^{(1-\delta,2-\delta_0)} n^{-\frac 32-\min{(\frac 12,\delta)}}.\eqno(A.18)
$$
On the other hand, we have
\begin{align}
&|B_{mn}^2(r)|+|B_{mn}^3(r)|\leqslant C\|f\|_{0,\alpha;Q}^{(1-\delta,2-\delta_0)}
\biggl\{\frac{e^{nr}}{\sqrt{2\pi nr}}\biggl[\int_r^L \frac{\sqrt \pi e^{-ns}}{\sqrt{2ns}} s^{\delta_0-1} ds \biggr] \notag\\
&\qquad +\frac{\sqrt \pi e^{-nr}}{\sqrt{2nr}}
\biggl[ \int_0^{\frac 1n} \frac{e^{ns} (\frac {ns}2)^{\m} }{\Gamma({\m}+1)} s^{\delta} ds +
\int_{\frac 1n}^1 \frac{e^{ns}}{\sqrt{2\pi ns}} s^{\delta} ds + \int_1^L \frac{e^{ns}}{\sqrt{2\pi ns}} s^{\delta_0-1} ds \biggr] \biggr\} \notag\\
& \leqslant C \|f\|_{0,\alpha;Q}^{(1-\delta,2-\delta_0)} \biggl\{ e^{nr}n^{-1}
\int_r^L e^{-ns} s^{\delta_0-\frac 32}ds+ e^{-nr}\biggl(\frac {n^{-\frac 12}\int_0^{\frac 1n}s^{\delta} ds}{2^{\m}\Gamma({\m}+1)} \notag\\
& \qquad  + n^{-1} \int_{\frac 1n}^1 e^{ns} s^{\delta-\frac 12} ds
+ n^{-1}\int_1^r e^{ns} s^{\delta_0-\frac 32} ds \biggr)\biggr\} \notag\\
& \leqslant C \|f\|_{0,\alpha;Q}^{(1-\delta,2-\delta_0)} \biggl(n^{-2} + \frac{n^{-\frac 32-\delta}} {2^{\m}\Gamma({\m}+1)} +  n^{-\frac 32-\min{(\frac 12,\delta)}} \biggr) \notag\\
& \leqslant C \|f\|_{0,\alpha;Q}^{(1-\delta,2-\delta_0)} n^{-\frac 32-\min{(\frac 12,\delta)}}.\tag{A.19}
\end{align}

Combining (A.18) with (A.19) yields
$$
\sum_{m=1}^{M}\sum_{n\geqslant 1 \atop 0<r \leqslant L} |R_{mn}^{(1)}(r)| \leqslant  C \|f\|_{0,\alpha:Q}^{(1-\delta,2-\delta_0)}. \eqno(A.20)
$$

{{\bf Case iv.} $m> M$ and $0 < r \leqslant L$.}

For the natational convenience, we set
$$
z(r)\equiv\frac{nr}{\m}, \quad \eta(r)\equiv\sqrt{1+z^2(r)} + \ln {\frac{z(r)}{1+\sqrt{1+z^2(r)}}},\quad
F(r)\equiv e^{-\m\eta(r)},\quad
\t{F}(r)=-\ds\f{1}{F(r)}.\eqno(A.21)
$$
It is easy to know that $F(r)$ is decreasing and $\t{F}(r)$ is increasing  with respect to $r$.

By (viii) in Lemma 2.3, we obtain
\begin{align}
 I_{\m}(nr)=I_{\m}(\m \frac{nr}{\m}) &\sim \frac{1}{\sqrt{2\pi \m}} \frac{\widetilde{F}(r)}{(1+z^2(r))^{\frac 14}};\tag{A.22}\\
K_{\m}(nr)&\sim \sqrt{\frac{\pi}{2\m}} \frac{F(r)}{(1+z^2(r))^{\frac 14}};\tag{A.23}\\
 I'_{\m}(nr)&\sim \frac{1}{\sqrt{2\pi \m}} \frac{(1+z^2(r))^{\frac 14}}{z(r)} \widetilde{F}(r);\tag{A.24}\\
 K'_{\m}(nr)&\sim -\sqrt{\frac{\pi}{2\m}} \frac{(1+z^2(r))^{\frac 14}}{z(r)} F(r).\tag{A.25}
\end{align}
By (A.17), we have
$$
|R_{mn}^{(1)}(r)|\leqslant  |B_{mn}^1|+|B_{mn}^2|+|B_{mn}^3|.\eqno(A.26)
$$
Next we deal with each $B_{mn}^i$ ($i=1,2,3$) in (A.26) separately.

\vskip 0.2 true cm

{\bf (A) Estimation of  $B_{mn}^1$}

By the inequalities in (4.35), (4.39) and (A.22)-(A.25), we obtain that
\begin{align*}
|B_{mn}^1| &\leqslant I_{\m}(nr) \biggl| \frac{K_{\m}^{'}(nL)}{I_{\m}^{'}(nL)} \biggl|
\int_0^L sI_{\m}(ns) |{f_L}_{mn}^{(1)}(s)| ds \notag\\
& \leqslant C \frac{F(L)}{\t{F}(L)} \frac{\t{F}(r)}{\m} \|f\|_{1,\alpha;Q}^{(1-\delta,2-\delta_0)}n^{-1}
 \biggl\{ \int_0^1 \frac{\t{F}(s)s^{\delta-1}}{(1+z^2(s))^{\frac 14}} ds
  + \int_1^L \frac{\t{F}(s)s^{\delta_0-2}}{(1+z^2(s))^{\frac 14}} ds \biggr\} \notag\\
& \leqslant C \frac{F(L)}{\t{F}(L)} \frac{\t{F}(r)}{\m} \|f\|_{1,\alpha;Q}^{(1-\delta,2-\delta_0)}n^{-1}
 \biggl\{ \int_0^1 \t{F'}(s) \frac{s^{\delta}}{(1+z^2(s))^{\frac 14}} \frac{z(s)}{ns} \frac{1}{\sqrt{1+z^2(s)}} ds \notag\\
  & \qquad +  \int_1^L \t{F'}(s) \frac{s^{\delta_0-1}}{(1+z^2(s))^{\frac 14}} \frac{z(s)}{ns} \frac{1}{\sqrt{1+z^2(s)}}ds \biggr\} \notag\\
  \end{align*}

  \begin{align}
  & \leqslant C \frac{F(L)}{\t{F}(L)} \frac{\t{F}(r)}{\m} \|f\|_{1,\alpha;Q}^{(1-\delta,2-\delta_0)}n^{-1}
 \biggl\{ \int_0^1 \t{F'}(s) \biggl(\frac{s}{z(s)}\biggr)^{\delta} \frac{\mu_m^{-\frac 12}}{(1+z^2(s))^{\frac 34-\frac \delta 2}} ds \notag\\
  & \qquad + \int_1^L \t{F'}(s) \biggl(\frac{s}{z(s)}\biggr)^{\delta_0}
  \frac{\mu_m^{-\frac 12}}{(1+z^2(s))^{\frac 34-\frac {\delta_0} {2}}}ds \biggr\} \notag\\
 & \leqslant C \frac{F(L)\t{F}(r)}{\t{F}(L)} \|f\|_{1,\alpha;Q}^{(1-\delta,2-\delta_0)}
 \biggl(\t{F}(1)n^{-1-\delta}\mu_m^{-1+\frac{\delta}{2}} + \t{F}(L) n^{-1-\delta_0}\mu_m^{-1+\frac{\delta_0}{2}}\biggr) \notag\\
& \leqslant C  \|f\|_{1,\alpha;Q}^{(1-\delta,2-\delta_0)}
 \biggl( n^{-1-\delta}\mu_m^{-1+\frac{\delta}{2}} +  n^{-1-\delta_0}\mu_m^{-1+\frac{\delta_0}{2}}\biggr).\tag{A.27}
\end{align}

{\bf (B) Estimation of  $B_{mn}^2$}

We will treat $B_{mn}^2$ in two cases of $0<r\leqslant 1$ and $r\ge 1$.

\vskip 0.2 true cm

{\bf Case (a) $0<r\leqslant 1$.}
\begin{align}
|B_{mn}^2| &\leqslant I_{\m}(nr) \biggl\{ \int_r^1 sK_{\m}(ns)|{f_L}_{mn}^{(1)}(s)| d s
  +  \int_1^L sK_{\m}(ns)|{f_L}_{mn}^{(1)}(s)| d s \biggr\} \notag\\
&\leqslant C\frac{1}{n\m} \frac{\widetilde{F}(r)}{(1+z^2(r))^{\frac 14}} \|f\|_{1,\alpha;Q}^{(1-\delta,2-\delta_0)}
 \biggl\{ \int_r^1 \frac{F(s)s^{\delta-1}}{(1+z^2(s))^{\frac 14}} ds
  + \int_1^L \frac{F(s)s^{\delta_0-2}}{(1+z^2(s))^{\frac 14}} ds \biggr\} \notag\\
 &\leqslant C\frac{1}{n\m} \frac{\t{F}(r)}{(1+z^2(r))^{\frac 14}} \|f\|_{1,\alpha;Q}^{(1-\delta,2-\delta_0)}
  \biggl\{ \int_r^1 \biggl(-F'(s)\biggr) \biggl(\frac{s}{z(s)}\biggr)^{\delta}
  \frac{\mu_m^{-\frac 12}}{(1+z^2(s))^{\frac 34-\frac \delta 2}} ds \notag\\
  & \qquad + \int_1^L \biggl(-F'(s)\biggr) \biggl(\frac{s}{z(s)}\biggr)^{\delta_0}
  \frac{\mu_m^{-\frac 12}}{(1+z^2(s))^{\frac 34-\frac {\delta_0} {2}}}ds \biggr\} \notag\\
& \leqslant C  \|f\|_{1,\alpha;Q}^{(1-\delta,2-\delta_0)} \t{F}(r)
 \biggl(F(r)n^{-1-\delta}\mu_m^{-1+\frac{\delta}{2}} + F(1) n^{-1-\delta_0}\mu_m^{-1+\frac{\delta_0}{2}}\biggr) \notag\\
 & \leqslant C  \|f\|_{1,\alpha;Q}^{(1-\delta,2-\delta_0)}
 \biggl(n^{-1-\delta}\mu_m^{-1+\frac{\delta}{2}} +  n^{-1-\delta_0}\mu_m^{-1+\frac{\delta_0}{2}}\biggr).\tag{A.28}
\end{align}

{\bf Case (b) $1\leqslant r \leqslant L$}

As in case (a), we can arrive at
\begin{align}
|B_{mn}^2| &\leqslant I_{\m}(nr) \int_r^L sK_{\m}(ns)|{f_L}_{mn}^{(1)}(s)| ds \notag\\
& \leqslant C  \|f\|_{1,\alpha;Q}^{(1-\delta,2-\delta_0)}
n^{-1-\delta_0}\mu_m^{-1+\frac{\delta_0}{2}}.\tag{A.29}
\end{align}

{\bf (C) Estimation of  $B_{mn}^3$}

As in (B) above, we also treat $B_{mn}^3$ in two cases of $0<r\leqslant 1$ and $r\ge 1$ separately.

{\bf Case (a) $0<r\leqslant 1$}
\begin{align}
|B_{mn}^3| &\leqslant K_{\m}(nr) \int_0^r sI_{\m}(ns)|{f_L}_{mn}^{(1)}(s)| ds \notag\\
&\leqslant C\sqrt{\frac{\pi}{2\m}} \frac{F(r)}{(1+z^2(r))^{\frac 14}} \|f\|_{1,\alpha;Q}^{(1-\delta,2-\delta_0)} n^{-1}
\biggl\{ \int_0^r \frac{1}{\sqrt{2\pi \m}} \frac{\t{F}(r)}{(1+z^2(r))^{\frac 14}}s^{\delta-1} ds \biggr\} \notag\\
&\leqslant C\frac{1}{n\m} F(r) \|f\|_{1,\alpha;Q}^{(1-\delta,2-\delta_0)}
\biggl\{ \int_0^r \t{F'}(s) \frac{s^{\delta}}{(1+z^2(s))^{\frac 14}} \frac{z(s)}{ns} \frac{1}{\sqrt{1+z^2(s)}} ds\biggr\} \notag\\
& \leqslant C \frac{1}{n\m} F(r) \|f\|_{1,\alpha;Q}^{(1-\delta,2-\delta_0)} \biggl\{\int_0^r \t{F'}(s) \biggl(\frac{s}{z(s)}\biggr)^{\delta} \frac{\mu_m^{-\frac 12}}{(1+z^2(s))^{\frac 34-\frac \delta 2}} ds \biggr\} \notag\\
&\leqslant C F(r) \|f\|_{1,\alpha;Q}^{(1-\delta,2-\delta_0)}
\t{F}(r)n^{-1-\delta}\mu_m^{-1+\frac{\delta}{2}} \notag\\
& \leqslant C  \|f\|_{1,\alpha;Q}^{(1-\delta,2-\delta_0)}
n^{-1-\delta}\mu_m^{-1+\frac{\delta}{2}}.\tag{A.30}
\end{align}

{\bf Case (b) $1 \leqslant r \leqslant L$}

In this case, we have
\begin{align}
|B_{mn}^3| &\leqslant K_{\m}(nr)  \biggl\{\int_0^1 sI_{\m}(ns)|{f_L}_{mn}^{(1)}(s)| d s
 + \int_0^1 sI_{\m}(ns)|{f_L}_{mn}^{(1)}(s)| ds \biggr\}  \notag\\
& \leqslant C \frac{1}{n\m} F(r) \|f\|_{1,\alpha;Q}^{(1-\delta,2-\delta_0)} \biggl\{\int_0^1 \t{F'}(s) \biggl(\frac{s}{z(s)}\biggr)^{\delta} \frac{\mu_m^{-\frac 12}}{(1+z^2(s))^{\frac 34-\frac \delta 2}} ds \notag\\
&\qquad + \int_1^r \t{F'}(s) \biggl(\frac{s}{z(s)}\biggr)^{\delta_0} \frac{\mu_m^{-\frac 12}}{(1+z^2(s))^{\frac 34-\frac {\delta_0}{2}}} ds \biggr\} \notag\\
& \leqslant C \|f\|_{1,\alpha;Q}^{(1-\delta,2-\delta_0)} F(r)
 \biggl(\t{F}(1)n^{-1-\delta}\mu_m^{-1+\frac{\delta}{2}} + \t{F}(r) n^{-1-\delta_0}\mu_m^{-1+\frac{\delta_0}{2}}\biggr) \notag\\
& \leqslant C  \|f\|_{1,\alpha;Q}^{(1-\delta,2-\delta_0)}
 \biggl(n^{-1-\delta}\mu_m^{-1+\frac{\delta}{2}} +  n^{-1-\delta_0}\mu_m^{-1+\frac{\delta_0}{2}}\biggr).\tag{A.31}
\end{align}
Substituting (A.27)-(A.31) into (A.26) yields
$$
\sum_{m=M}^{\infty}\sum_{n=1}^{\infty} |R_{mn}^{(1)}(r)| \leqslant  C \|f\|_{1,\alpha;Q}^{(1-\delta,2-\delta_0)}. \eqno(A.32)
$$
By (A.15)-(A.16), (A.20) and (A.32), we complete the proof of Lemma A.4.\qed

{\bf Lemma A.5.} {\it For the term  $\ds\sum\limits_{n=1}^{\infty} n^2 |R_{0n}^{(1)}(r)|$ defined in (4.46), we have
\begin{equation}
\ds\sum\limits_{n=1}^{\infty} n^2 |R_{0n}^{(1)}(r)|\le
\begin{cases}
&C \|f\|_{3,\alpha;Q}^{(1-\delta,2-\delta_0)}r^{\dl-\f52},\quad 0<r\le 1,\\
&C \|f\|_{3,\alpha;Q}^{(1-\delta,2-\delta_0)},\quad 1<r\le L,
\end{cases}\tag{A.33}
\end{equation}
where $C>0$ is independent of $L$.}

{\bf Proof.} As in (A.5), we write $R_{0n}^{(1)}(r)=B_{1}^n-B_{2}^n$. First, we estimate
$\ds\sum\limits_{n=1}^{\infty} n^2 |B_{1}^{n}(r)|$.
It follows from (iii)-(vii) of Lemma 2.3, the inequalities in (4.34) and (4.38), and the fact of $L\geqslant 4$ that
\begin{align}
&\ds\sum\limits_{n=1}^{\infty} n^2 |B_{1}^{n}(r)|\leqslant \sum_{n=1}^{\infty} n^2 \frac{\sqrt{\frac{\pi}{2nL}}e^{-nL}}{\frac{e^{nL}}{\sqrt{2\pi nL}}} e^{nr}
 \biggl(\int_0^{\frac 1n}+\int_{\frac 1n}^1+\int_1^{L}\biggr)s I_{0}(ns) |{f_L}_{0n}^{(1)}(s)| ds \notag\\
& \leqslant C\sum_{n=1}^{\infty} n^2  e^{nr-2nL} \biggl\{\int_0^{\frac 1n}  e^{ns}s^{\delta} ds \|f\|_{0,\alpha;Q}^{(1-\delta,2-\delta_0)}+
  \biggl(\int_{\frac 1n}^1 e^{ns}s^{\delta-3} ds \notag \\
&\qquad +\int_1^{L} e^{ns}s^{\delta_0-4} d s \biggr) n^{-3} \|f\|_{3,\alpha;Q}^{(1-\delta,2-\delta_0)} \biggr\} \notag\\
& \leqslant C\sum_{n=1}^{\infty} \biggl\{n^{1-\delta}e^{-nL} \|f\|_{0,\alpha;Q}^{(1-\delta,2-\delta_0)}+
  \bigl( n^{1-\delta} e^{n(1-L)}+n^{-2}\bigr) \|f\|_{3,\alpha;Q}^{(1-\delta,2-\delta_0)}\biggr\} \notag\\
& \leqslant C\|f\|_{3,\alpha;Q}^{(1-\delta,2-\delta_0)}
  \sum_{n=1}^{\infty} \bigl(n^{1-\delta}e^{-2n} + n^{1-\delta}e^{-n} + n^{-2} \bigr) \notag\\
& \leqslant C \|f\|_{3,\alpha;Q}^{(1-\delta,2-\delta_0)}.\tag{A.34}
\end{align}

To estimate $\ds\sum\limits_{n=1}^{\infty} n^2 |B_{2}^{n}(r)|$, we will divide this procedure into the following cases:
\vskip 0.2 true cm

{\bf Case i. $0<r\leqslant 1$}

In this case, we have
\begin{align*}
&\ds\sum\limits_{n=1}^{\infty} n^2 |B_{2}^{n}(r)|
\leqslant \sum_{n=1}^{\infty} n^2
  \biggl\{I_{0}(nr) \biggl(\int_{r}^1 K_{0}(ns) |{f_L}_{0n}^{(1)}(s)| d s +\int_1^{L} s K_{0}(ns) |{f_L}_{0n}^{(1)}(s)| d s \biggr)\notag\\
  &\qquad + K_{0}(nr)\biggl(\int_0^{\frac r2} sI_{0}(ns) |{f_L}_{0n}^{(1)}(s)|ds + \int_{\frac r2}^r sI_{0}(ns) |{f_L}_{0n}^{(1)}(s)|ds \biggr)\biggr\}\notag\\
& \leqslant \sum_{n=1}^{\infty} C n^2 \biggl\{e^{nr} n^{-2} \|f\|_{2,\alpha;Q}^{(1-\delta,2-\delta_0)}
  \biggl( \int_{r}^1 \frac{e^{-ns}}{\sqrt{2ns}}s^{\delta-2} d s +\int_1^{L} \frac{e^{-ns}}{\sqrt{2ns}}s^{\delta_0-3} d s \biggr) \notag\\
  & \qquad + \frac{e^{-nr}}{\sqrt{2nr}}\biggl( \|f\|_{0,\alpha;Q}^{(1-\delta,2-\delta_0)}
\int_0^{\frac r2} e^{ns}s^{\delta}ds + \|f\|_{2,\alpha;Q}^{(1-\delta,2-\delta_0)}n^{-2} \int_{\frac r2}^r e^{ns}s^{\delta-2}ds \biggr) \biggr\} \notag\\
& \leqslant  C \|f\|_{2,\alpha;Q}^{(1-\delta,2-\delta_0)}\sum_{n=1}^{\infty} \biggl\{e^{nr}
  \biggl( n^{-\frac 12}\int_{r}^1 e^{-ns}s^{\delta-\frac 52} d s
  + n^{-\frac 12}\int_1^{L} e^{-ns}s^{\delta_0-\frac 72} d s \biggr) \notag\\
 \end{align*}

  \begin{align}
  &\quad + \frac{e^{-nr}}{\sqrt{2nr}} \biggl( n^2 \int_0^{\frac r2} e^{ns}s^{\delta} ds + \int_{\frac r2}^r e^{ns}s^{\delta-2}ds \biggr) \biggr\} \notag\\
& \leqslant C \|f\|_{2,\alpha;Q}^{(1-\delta,2-\delta_0)} \sum_{n=1}^{\infty} \biggl(n^{-\frac 32} r^{\delta-\frac 52} +e^{-\frac{nr}2}n^{\frac 12} r^{\delta-\frac 12} +r^{\delta-\frac 52}n^{-\frac 32} \biggr) \notag\\
& \leqslant C \|f\|_{2,\alpha;Q}^{(1-\delta,2-\delta_0)} r^{\delta-\frac 52}.\tag{A.35}
\end{align}

{\bf Case ii. $ 1\leqslant r \leqslant L$}

 At this time, we obtain
\begin{align}
&\ds\sum\limits_{n=1}^{\infty} n^2 |B_{2}^{n}(r)|
\leqslant \sum_{n=1}^{\infty} n^2
\biggl\{e^{nr} \int_{r}^L K_{0}(ns) |{f_L}_{0n}^{(1)}(s)| ds \notag \\
&\quad + K_{0}(nr)\biggl(\int_0^{\frac 12} sI_{0}(ns) |{f_L}_{0n}^{(1)}(s)|ds + \int_{\frac 12}^1 sI_{0}(ns) |{f_L}_{0n}^{(1)}(s)|ds +\int_1^r sI_{0}(ns) |{f_L}_{0n}^{(1)}(s)|ds \biggr)\biggr\} \notag\\
& \leqslant C\sum_{n=1}^{\infty} \biggl\{e^{nr} \|f\|_{2,\alpha;Q}^{(1-\delta,2-\delta_0)} n^{-2} \int_r^{L} \frac{e^{-ns}}{\sqrt{2ns}}s^{\delta_0-3} ds \notag\\
& \quad
+ \frac{e^{-nr}}{\sqrt{2nr}}\biggl(\|f\|_{0,\alpha;Q}^{(1-\delta,2-\delta_0)} \int_0^{\frac 12} e^{ns}s^{\delta}+ \|f\|_{2,\alpha}^{(1-\delta,2-\delta_0)} n^{-2} \biggl[ \int_{\frac 12}^1 e^{ns}s^{\delta-2}ds +\int_1^r e^{ns} s^{\delta_0-3}ds\biggl] \biggr) \biggr\} \notag\\
& \leqslant C\|f\|_{2,\alpha;Q}^{(1-\delta,2-\delta_0)} \sum_{n=1}^{\infty} \biggl\{e^{nr}
\bigl( n^{-\frac 32}e^{-nr}\bigr) + e^{-nr} \biggl( n^{\frac 12} e^{\frac n2}+ n^{-\frac 32}e^{n}
+ n^{-\frac32}e^{nr} \biggr) \biggr\} \notag\\
& \leqslant C\|f\|_{2,\alpha;Q}^{(1-\delta,2-\delta_0)} \sum_{n=1}^{\infty} \bigl(n^{\frac 12} e^{-\frac n2}+ n^{-\frac 32}\bigr) \notag\\
& \leqslant C \|f\|_{2,\alpha;Q}^{(1-\delta,2-\delta_0)}.\tag{A.36}
\end{align}
Combining (A.34)-(A.36) yields (A.33). \qed

{\bf Lemma A.6.} {\it For the term $\sum\limits_{m=1}^{\infty}\sum\limits_{n=1}^{\infty} n^2 |R_{mn}^{(1)}(r)|$ defined in (4.46), we have
\begin{equation*}
\sum\limits_{m=1}^{\infty}\sum\limits_{n=1}^{\infty} n^2 |R_{mn}^{(1)}(r)|\le
\begin{cases}
&C \|f\|_{4,\alpha;Q}^{(1-\delta,2-\delta_0)}r^{\min\{\delta,\frac 12\}-\frac 52}, \qquad r\leqslant 1,\\
&C \|f\|_{4,\alpha;Q}^{(1-\delta,2-\delta_0)}, \qquad \qquad\quad  1<r\leqslant L,
\end{cases}\tag{A.37}
\end{equation*}
where $C>0$ is independent of $L$.}

{\bf Proof.}  To prove (A.37), we will divide this procedure into the following four cases.
As in (A.26), one has $|R_{mn}^{(1)}(r)|\leqslant  B_{mn}^1+B_{mn}^2+B_{mn}^3.$

 {{\bf Case i.} \quad $m\leqslant M$, and $n\leqslant \ds\frac 1r $ with $0<r \leqslant 1$}

At this time, we can choose a positive integer $N$ such that  $Nr \leqslant1$ holds.
And by (iv)-(vii) of Lemma 2.3 and the inequalities in (4.35) and (4.39), we have
\begin{align*}
&n^2 |R_{mn}^{(1)}(r)|\leqslant  n^2 I_{\m}(nr)\biggl(\frac{K_{\m}^{'}(nL)}{I_{\m}^{'}(nL)}
\int_0^L sI_{\m}(ns)|{f_L}_{mn}^{(1)}(s)| d s\\
&\qquad +\int_r^L sK_{\m}(ns)|{f_L}_{mn}^{(1)}(s)| d s\biggr)
+n^2 K_{\m}(nr)\int_0^r sI_{\m}(ns)|{f_L}_{mn}^{(1)}(s)| d s\\
& \leqslant C\|f\|_{2,\alpha;Q}^{(1-\delta,2-\delta_0)} \frac{e^{nr} (\frac {nr}2)^{\m}}{\Gamma(\m+1)}
n^2
\biggl\{ e^{-2nL} \biggl[ \int_0^{\frac 1n} \frac{e^{ns} (\frac {ns}2)^{\m}s^{\delta}}{\Gamma({\m}+1)}  ds
+ \int_{\frac 1n}^1 \frac{e^{ns}s^{\delta}}{\sqrt{2\pi ns}}  ds\\
&\qquad  +\ds\f{1}{n^2}\int_1^L \frac{e^{ns} s^{\delta_0-3}}{\sqrt{2\pi ns}} ds \biggr]
+ n^{-2}\biggl[\int_r^{Nr} \frac{e^{ns} \Gamma(\m) 2^{\m-1}}{(ns)^{\m}} s^{\delta-2} ds
+ \int_{Nr}^1 \frac{\sqrt \pi e^{-ns}}{\sqrt{2ns}} s^{\delta-2} ds\\
&\qquad + \int_1^L \frac{\sqrt \pi e^{-ns}}{\sqrt{2ns}} s^{\delta_0-3} ds \biggr] \biggr\}
+ \|f\|_{0,\alpha;Q}^{(1-\delta,2-\delta_0)}\frac{e^{nr} \Gamma(\m) 2^{\m-1}}{(nr)^{\m}}
\int_0^r \frac{e^{ns} (\frac {ns}2)^{\m} }{\Gamma(\m+1)} s^{\delta} ds\\
& \leqslant C \|f\|_{2,\alpha;Q}^{(1-\delta,2-\delta_0)} n^2 \biggl\{ \frac 1{e^{2nL} 2^{\m}\Gamma(\m+1)}
\biggl[ \int_0^{\frac 1n}  s^{\delta} ds + n^{-\frac 12}\int_{\frac 1n}^1 e^{ns} s^{\delta-\frac 12} ds\\
&\qquad +
n^{-\frac 32} \int_1^L e^{ns} s^{\delta_0-\frac 72} ds \biggr]
+ n^{-2}\biggl[ \int_r^{Nr} \frac{e^{ns} s^{\delta-2}}{\m}ds + \int_{Nr}^1 \frac{e^{-ns} s^{\delta-2}}{2^{\m}} ds\\
 &\qquad + n^{-\frac 12}\int_1^L \frac{e^{-ns} s^{\delta_0-\frac 72}}{ 2^{\m}} ds \biggr] + \int_0^r \frac{e^{ns} s^{\delta}}{\m}ds \biggr\}\\
& \leqslant C \|f\|_{2,\alpha;Q}^{(1-\delta,2-\delta_0)} \biggl(
  \frac {n^{1-\min{(\frac 12,\delta)}} e^{-n} }{2^{\m}}+ \frac {n^{-\frac 32} e^{-n} }{2^{\m}}+ \frac {n^{-1}}{\m}r^{\delta-2}
  + \frac {n^{-1}e^{-\frac n2}}{2^{\m}} \biggr).
\end{align*}

Thus we get
\begin{equation*}
\sum_{m=1}^{M}\sum_{n\leqslant\frac 1r \atop 0<r \leqslant 1} n^2|R_{mn}^{(1)}(r)|\leqslant C\|f\|_{2,\alpha;Q}^{(1-\delta,2-\delta_0)}r^{\delta-2}.\tag{A.38}
\end{equation*}
{{\bf Case ii.}\quad  $m\leqslant M$, and $n\geqslant \frac 1r $ with $0<r \leqslant 1$}

In this case, we get
\begin{align*}
& n^2 |R_{mn}^{(1)}(r)|\\
&\leqslant C n^2 \|f\|_{2,\alpha;Q}^{(1-\delta,2-\delta_0)}\frac{e^{nr}}{\sqrt{2\pi nr}}
\biggl\{ e^{-2nL} \biggl[ \int_0^{\frac 1n} \frac{e^{ns} (\frac {ns}2)^{\m}}{\Gamma(\m+1)} s^{\delta} ds
+ \int_{\frac 1n}^1 \frac{e^{ns}}{\sqrt{2\pi ns}} s^{\delta} ds\\
&\qquad  +n^{-2} \int_1^L \frac{e^{ns}}{\sqrt{2\pi ns}} s^{\delta_0-3} ds \biggr]
+ n^{-2} \biggl[\int_r^1 \frac{\sqrt \pi e^{-ns}}{\sqrt{2ns}} s^{\delta-2} ds + \int_1^L \frac{\sqrt \pi e^{-ns}}{\sqrt{2ns}} s^{\delta_0-3} ds \biggr] \biggr\}\\
\end{align*}
\begin{align*}
&\qquad + n^2 \frac{\sqrt \pi e^{-nr}}{\sqrt{2nr}}
\biggl[ \|f\|_{0,\alpha;Q}^{(1-\delta,2-\delta_0)} \int_0^{\frac r2} \frac{e^{ns}}{\sqrt{2\pi ns}} s^{\delta} ds + n^{-2} \|f\|_{2,\alpha}^{(1-\delta,2-\delta_0)}\int_{\frac r2}^r \frac{e^{ns}}{\sqrt{2\pi ns}} s^{\delta-2} ds \biggr] \\
& \leqslant C \|f\|_{2,\alpha;Q}^{(1-\delta,2-\delta_0)} n^2 \biggl\{ e^{nr-2nL}
\biggl[ \frac {\int_0^{\frac 1n}s^{\delta} ds}{2^{\m}\Gamma({\m}+1)}
+ n^{-\frac 12}\int_{\frac 1n}^1 e^{ns} s^{\delta-\frac 12} ds +
n^{-\frac 52} \int_1^L e^{ns} s^{\delta_0-\frac 72}ds \biggr]\\
&\qquad +e^{nr} n^{-\frac 52} \biggl[ \int_r^1e^{-ns} s^{\delta-\frac 52} ds + \int_1^L e^{-ns} s^{\delta_0-\frac 72} ds \biggr] +e^{-nr}n^{-\frac 12}\biggl[\int_0^{\frac r2} e^{ns}s^{\delta-\frac 12} ds\\
&\qquad + n^{-2}\int_{\frac r2}^r e^{ns} s^{\delta-\frac 52} ds \biggr] \biggr\}\\
& \leqslant C \|f\|_{2,\alpha;Q}^{(1-\delta,2-\delta_0)} \biggl(n^{1-\min{(\frac 12,\delta)}}e^{-n} + n^{-\frac 32}+ n^{-\frac 32} r^{\delta-\frac 52} +  n^{\frac 32} e^{-\frac {nr}{2}} r^{\delta+\frac 12}\biggr).
\end{align*}

This, together with $\ds\sum_{ n\geqslant\frac 1r} n^{\frac 32} e^{-\frac {nr}{2}} r^{\delta+\frac 12} \leqslant C r^{\delta-2}$, yields
\begin{equation*}
\sum_{m=1}^{M}\sum_{n\geqslant\frac 1r \atop 0<r \leqslant 1} n^2|R_{mn}^{(1)}(r)|\leqslant C\|f\|_{2,\alpha;Q}^{(1-\delta,2-\delta_0)}r^{\delta-\frac 52}.\tag{A.39}
\end{equation*}

{{\bf Case iii.}\quad $m\leqslant M$ and $1 \leqslant r \leqslant L$}

As in Case (b), one has
$$
n^2 |B_{mn}^1|\leqslant C \|f\|_{2,\alpha;Q}^{(1-\delta,2-\delta_0)} \biggl(\frac{ n^{1-\delta} e^{-n}} {2^{\m}\Gamma({\m}+1)}
+ n^{1-\min{(\frac 12,\delta)}}e^{-n} + n^{-\frac 32}\biggr).
$$

In addition,
\begin{align*}
& n^2(|B_{mn}^2|+|B_{mn}^3|)\\
&\leqslant C n^2
\biggl\{\frac{e^{nr}}{\sqrt{2\pi nr}}\|f\|_{2,\alpha;Q}^{(1-\delta,2-\delta_0)} n^{-2} \biggl[\int_r^L \frac{\sqrt \pi e^{-ns}}{\sqrt{2ns}} s^{\delta_0-3} ds \biggr]\\
&\qquad +\frac{\sqrt \pi e^{-nr}}{\sqrt{2nr}}
\biggl[ \|f\|_{0,\alpha;Q}^{(1-\delta,2-\delta_0)}\int_0^{\frac 12} \frac{e^{ns}}{\sqrt{2\pi ns}} s^{\delta} ds\\
&\qquad  +
\|f\|_{2,\alpha;Q}^{(1-\delta,2-\delta_0)} n^{-2}\biggl(\int_{\frac 12}^1 \frac{e^{ns}}{\sqrt{2\pi ns}} s^{\delta-2} ds + \int_1^L \frac{e^{ns}}{\sqrt{2\pi ns}} s^{\delta_0-3} ds\biggr) \biggr] \biggr\}\\
& \leqslant C \|f\|_{2,\alpha;Q}^{(1-\delta,2-\delta_0)} \biggl\{ e^{nr}n^{-1}
\int_r^L e^{-ns} s^{\delta_0-\frac 72}ds \\
& \qquad + e^{-nr}\biggl(n \int_0^{\frac 12} e^{ns}s^{\delta-\frac 12} ds + n^{-1} \int_{\frac 12}^1 e^{ns} s^{\delta-\frac 52} ds
+ n^{-1}\int_1^r e^{ns} s^{\delta_0-\frac 72} ds \biggr)\biggr\}\\
& \leqslant C \|f\|_{2,\alpha;Q}^{(1-\delta,2-\delta_0)} (n^{-2} + n e^{-n}).
\end{align*}
Hence, we get
\begin{equation*}
\sum_{m=1}^{M}\sum_{n\geqslant1 \atop 1 \leqslant r \leqslant 1} n^2|R_{mn}^{(1)}(r)|\leqslant C\|f\|_{2,\alpha;Q}^{(1-\delta,2-\delta_0)}.
\tag{A.40}
\end{equation*}

{{\bf Case iv.}\quad $m> M$ and $0 < r \leqslant L$}

Assume that $z(r),\eta(r), F(r), \t F(r)$ are the functions defined in (A.21).
Following (ix) in Lemma 2.3, it is easy to know that $F(r_2)\t{F}(r_1)$ $\leqslant e^{-n(r_2 -r_1)}$ holds for $r_1\leqslant r_2$.
This, together with (A.22)-(A.25) and a direct computation, yields
\begin{align}
n^2|B_{mn}^1| &\leqslant C \frac{F(L)}{\t{F}(L)} \frac{\t{F}(r)n^{2}}{\m}
\biggl\{ \|f\|_{2,\alpha;Q}^{(1-\delta,2-\delta_0)} \mu_m^{-1} \int_0^1 \frac{\t{F}(s)s^{\delta}}{(1+z^2(s))^{\frac 14}} ds \notag\\
&\quad + \|f\|_{4,\alpha;Q}^{(1-\delta,2-\delta_0)} \mu_m^{-1}n^{-2} \int_1^L \frac{\t{F}(s)s^{\delta_0-3}}{(1+z^2(s))^{\frac 14}} ds \biggr\} \notag\\
& \leqslant C \frac{F(L) \t{F}(r)}{\t{F}(L)}
\|f\|_{4,\alpha;Q}^{(1-\delta,2-\delta_0)} n^2 \mu_m^{-\frac 32}
\biggl\{ \int_0^1 \t{F'}(s) \frac{s^{\delta}}{(1+z^2(s))^{\frac 14}} \frac{z(s)}{n\sqrt{1+z^2(s)}} ds \notag\\
& \qquad + n^{-2} \int_1^L \t{F'}(s) \frac{s^{\delta_0-3}}{(1+z^2(s))^{\frac 14}} \frac{z(s)}{n\sqrt{1+z^2(s)}}ds \biggr\} \notag\\
& \leqslant C \frac{F(L) \t{F}(r)}{\t{F}(L)}
  \|f\|_{4,\alpha;Q}^{(1-\delta,2-\delta_0)} n^2 \mu_m^{-\frac 32}
   \biggl\{ n^{-1} \int_0^1 \t{F'}(s) \frac{z(s)}{\sqrt{1+z^2(s)}} ds \notag\\
  & \qquad + n^{-3}\int_1^L \t{F'}(s) \frac{z(s)}{\sqrt{1+z^2(s)}} \biggl(\frac{s}{z(s)}\biggr)^{\frac 12}  s^{\delta_0-\frac 52} ds \biggr\} \notag\\
 & \leqslant C \frac{\t{F}(r)}{\t{F}(L)} \|f\|_{4,\alpha;Q}^{(1-\delta,2-\delta_0)}
 \biggl(n \mu_m^{-\frac 32} F(L) \t{F}(1) + n^{-\frac 32} \mu_m^{-\frac 54} F(L)\t{F}(L) \biggr) \notag\\
& \leqslant C  \|f\|_{4,\alpha;Q}^{(1-\delta,2-\delta_0)}
 \bigl(n \mu_m^{-\frac 32} e^{-n(L-1)} + n^{-\frac 32} \mu_m^{-\frac 54} \bigr) \notag\\
 & \leqslant C  \|f\|_{4,\alpha;Q}^{(1-\delta,2-\delta_0)}
 \bigl(n e^{-n} \mu_m^{-\frac 32} + n^{-\frac 32} \mu_m^{-\frac 54} \bigr).\tag{A.41}
\end{align}

In addition, we obtain for $r\leqslant 1$
\begin{align}
n^2|B_{mn}^2| &\leqslant C\frac{n^2}{\m} \frac{\t{F}(r)}{(1+z^2(r))^{\frac 14}} \|f\|_{4,\alpha;Q}^{(1-\delta,2-\delta_0)}
\biggl\{ n^{-2}\mu_m^{-1} \int_r^1 \frac{F(s)s^{\delta-2}}{(1+z^2(s))^{\frac 14}} ds \notag\\
&\qquad +  n^{-2}\mu_m^{-1} \int_1^L \frac{F(s)s^{\delta_0-3}}{(1+z^2(s))^{\frac 14}} ds \biggr\} \notag\\
&\leqslant C \|f\|_{4,\alpha;Q}^{(1-\delta,2-\delta_0)} \widetilde{F}(r)
n^{-1}\mu_m^{-\frac 32} \biggl\{ \int_r^1 \biggl(-F'(s)\biggr) \frac{z(s)}{\sqrt{1+z^2(s)}}
\biggl(\frac{1}{z(s)}\biggr)^{\frac 12} s^{\delta-2} ds \notag\\
& \qquad + \int_1^L \biggl(-F'(s)\biggr)  \frac{z(s)}{\sqrt{1+z^2(s)}} \biggl(\frac{1}{z(r)}\biggr)^{\frac 12}
s^{\delta_0-3} ds \biggr\} \notag\\
&\leqslant C \|f\|_{4,\alpha;Q}^{(1-\delta,2-\delta_0)} n^{-\frac 32}\mu_m^{-\frac 54}
\biggl(\t{F}(r)F(r)r^{\delta-\frac 52}+ \t{F}(r)F(1)\biggr) \notag\\
&\leqslant C  \|f\|_{4,\alpha;Q}^{(1-\delta,2-\delta_0)} n^{-\frac 32}\mu_m^{-\frac 54} r^{\delta-\frac 52},\tag{A.42}
\end{align}
and for $r \geqslant 1$,
\begin{align}
n^2|B_{mn}^2| &\leqslant n^2 I_{\m}(nr) \int_r^L sK_{\m}(ns)|{f_L}_{mn}^{(1)}(s)| ds \notag\\
& \leqslant C \|f\|_{4,\alpha;Q}^{(1-\delta,2-\delta_0)} n^{-\frac 32}\mu_m^{-\frac 54}.\tag{A.43}
\end{align}

Finally, we estimate $n^2 |B_{mn}^3|$. For $0<r\leqslant 1$, we have
\begin{align}
n^2|B_{mn}^3| &\leqslant n^2 K_{\m}(nr) \biggl(\int_0^{\frac r2} sI_{\m}(ns)|{f_L}_{mn}^{(1)}(s)| d s +\int_{\frac r2}^r sI_{\m}(ns)|{f_L}_{mn}^{(1)}(s)| d s \biggr) \notag\\
&\leqslant C\frac{n^2}{\m} \frac{F(r)}{(1+z^2(r))^{\frac 14}}
\biggl\{\|f\|_{2,\alpha;Q}^{(1-\delta,2-\delta_0)}\mu_m^{-1} \int_0^{\frac r2} \frac{\t{F}(r)}{(1+z^2(r))^{\frac 14}}s^{\delta} ds \notag\\
&\quad + \|f\|_{4,\alpha;Q}^{(1-\delta,2-\delta_0)}n^{-2}\mu_m^{-1} \int_{\frac r2}^r \frac{\t{F}(r)}{(1+z^2(r))^{\frac 14}}s^{\delta-2}ds\biggr\} \notag\\
&\leqslant C\frac{n^2}{\m} F(r) \|f\|_{4,\alpha;Q}^{(1-\delta,2-\delta_0)}
\biggl\{ n^{-1} \mu_m^{-1} \int_0^{\frac r2} \t{F'}(s) \frac{s^{\delta}}{(1+z^2(s))^{\frac 14}} \frac{z(s)}{\sqrt{1+z^2(s)}} ds \notag\\
& \quad + n^{-3} \mu_m^{-1} \int_{\frac r2}^r \t{F'}(s) \frac{s^{\delta-2}}{(1+z^2(s))^{\frac 14}} \frac{z(s)}{\sqrt{1+z^2(s)}} ds \biggr\} \notag\\
& \leqslant C \frac{n^2}{\m} F(r) \|f\|_{4,\alpha;Q}^{(1-\delta,2-\delta_0)}
\biggl(n^{-1} \mu_m^{-1} \int_0^{\frac r2} \t{F'}(s) ds
+ n^{-3} \mu_m^{-1} \int_{\frac r2}^r \t{F'}(s) \frac{s^{\delta-2}}{(z(s))^{\frac 12}} ds \biggr) \notag\\
&\leqslant C \frac{n^2}{\m} F(r) \|f\|_{4,\alpha;Q}^{(1-\delta,2-\delta_0)}
\biggl(n^{-1} \mu_m^{-1}  \t{F}(\frac r2)
+ n^{-\frac 72} \mu_m^{-\frac 34} \t{F}(r) r^{\delta-\frac 52} \biggr) \notag\\
&\leqslant C \|f\|_{4,\alpha;Q}^{(1-\delta,2-\delta_0)}
\bigl(n e^{-\frac{nr}{2}}\mu_m^{-\frac 32}
+ n^{-\frac 32} \mu_m^{-\frac 54} r^{\delta-\frac 52} \bigr). \tag{A.44}
\end{align}

For $r\geqslant 1$, we obtain
\begin{align}
&n^2|B_{mn}^3|\leqslant C\frac{n^2}{\m} \frac{F(r)}{(1+z^2(r))^{\frac 14}}
\biggl\{\|f\|_{2,\alpha;Q}^{(1-\delta,2-\delta_0)}\mu_m^{-1} \int_0^{\frac 12} \frac{\t{F}(r)}{(1+z^2(r))^{\frac 14}}s^{\delta} ds\notag\\
& \quad+ \|f\|_{4,\alpha;Q}^{(1-\delta,2-\delta_0)}n^{-2}\mu_m^{-1}
\biggl( \int_{\frac 12}^1 \frac{\t{F}(r)}{(1+z^2(r))^{\frac 14}}s^{\delta-2} ds + \int_1^r \frac{\widetilde{F}(r)}{(1+z^2(r))^{\frac 14}}s^{\delta_0-3} ds \biggr)\biggr\}\notag\\
& \leqslant C \frac{n^2}{\m} F(r) \|f\|_{4,\alpha;Q}^{(1-\delta,2-\delta_0)}
\biggl\{n^{-1} \mu_m^{-1} \int_0^{\frac 12} \t{F'}(s) ds
+ n^{-3} \mu_m^{-1} \biggl(\int_{\frac 12}^1 \t{F'}(s) \frac{s^{\delta-2}}{(z(s))^{\frac 12}} ds\notag\\
&\qquad +
\int_1^r \t{F'}(s) \frac{s^{\delta_0-3}}{(z(s))^{\frac 12}} ds \biggr) \biggr\}\notag\\
& \leqslant C \frac{n^2}{\m} F(r) \|f\|_{4,\alpha;Q}^{(1-\delta,2-\delta_0)}
\biggl(n^{-1} \mu_m^{-1}  \t{F}(\frac 12)
+ n^{-\frac 72} \mu_m^{-\frac 34} \t{F}(1) + n^{-\frac 72} \mu_m^{-\frac 34} \t{F}(r) \biggr)\notag\\
& \leqslant C \|f\|_{4,\alpha;Q}^{(1-\delta,2-\delta_0)}
\biggl(n e^{-\frac{n}{2}}\mu_m^{-\frac 32}
+ n^{-\frac 32} \mu_m^{-\frac 54}\biggr).\tag{A.45}
\end{align}
Combining (A.41)-(A.45) yields
\begin{equation}
\left\{
\begin{aligned}
&\sum\limits_{m > M}\sum\limits_{n\geqslant 1 \atop 0 < r\leqslant L} n^2|R_{mn}^{(1)}(r)|\leqslant C \|f\|_{4,\alpha;Q}^{(1-\delta,2-\delta_0)}r^{\min\{\delta,\frac 12\}-\frac 52},\qquad r\leqslant 1,\\
&\sum\limits_{m > M}\sum\limits_{n\geqslant 1 \atop 0 < r\leqslant L} n^2|R_{mn}^{(1)}(r)|\leqslant C \|f\|_{4,\alpha;Q}^{(1-\delta,2-\delta_0)},\qquad 1<r\leqslant L.
\end{aligned}
\right.\tag{A.46}
\end{equation}
Therefore, collecting  (A.38)-(A.40) and (A.46) yields (A.37).\qed

\vskip 0.6 true cm

{\bf Acknowledgements.}  {\it Yin Huicheng wishes to express his gratitude to Professor Xin Zhouping, Chinese University
of Hong Kong,  Professor Witt Ingo, University of G\"otingen, and Professor B. W. Schulze, University of Potsdam
for their many fruitful discussions in this problem.}

\end{document}